\documentclass[graybox,envcountchap,sectrefs]{svmono}
\usepackage[utf8]{inputenc}

\usepackage{type1cm}         

\usepackage{makeidx}         
\usepackage{graphicx}        
\usepackage{caption}
\usepackage{multicol}        
\usepackage[bottom]{footmisc}

\usepackage{newtxtext}       %
\usepackage[varg,charter]{newtxmath}

\usepackage{enumitem}       

\usepackage{hyperref}

\makeindex             
\begin{document}

\author{Enrique Zuazua}
\title{Exact Controllability and Stabilization of the Wave Equation}
\subtitle{-- Monograph --}
\maketitle

\frontmatter

%
%
%

\begin{dedication}
Dedicated to the memory of Luiz Adauto Medeiros, and to my friend and colleague Fred Vasconcellos, and, with them, to the whole Brazilian mathematical community.
\end{dedication}

%
%

\foreword

Enrique Zuazua has written his thesis in applied mathematics under my direction and defended it in the Jacques-Louis Lions Laboratory (LJLL) of Paris in 1988. A very important part of his scientific activity has been, from the very beginning, devoted to the behavior of waves in a bounded domain: oscillations, stability, stabilization and control. After the results on damped systems obtained in the first part of his thesis, he was lucky enough to get involved quite early in the application of the HUM method for exact controllability, just devised around 1986 by Jacques-Louis Lions who had in mind some applications to the stabilization of large spatial structures. To stabilize these structures, some thermal effects have to be taken into account, and the real equations become in fact coupled systems of hyperbolic and parabolic equations. Because parabolic problems are not reversible in time, the controllability of semilinear heat equations has become a field of research in itself, based on quite different methods. But the core of the problem remains control and stabilization of systems governed by evolution equations of the second order in time, which is the object of the present monograph.

Exact controllability of hyperbolic systems can be, whenever possible, obtained by essentially three different methods: propagation arguments related to Holmgren's uniqueness theorem, harmonic analysis in connection with Ingham's lemma and its extensions, and the multiplier method. These three methods can be classified in the following way: the propagation method is the most general one, giving almost optimal results whenever applicable, but does not usually provide any clue for numerical treatment, since it does not seem to be robust, except in one space-dimension where the propagation can be understood by means of Riemann's integral formulas. The harmonic method provides very strong results when applicable, but basically precludes any perturbation; in particular semilinear perturbations will never work in more than one dimension, and only in very special situations (nonlinear oscillations and exponential stabilization). The multiplier method requires geometrical conditions which substantially restrict its applicability, but in contrast with the two other techniques it is, in any dimension, a basically robust method from the point of view of operators, which allows at the same time semilinear perturbations and numerical approximation.

The first four chapters of the monograph are devoted to controllability results relying on the multiplier method and the HUM method of J.-L. Lions. The last two chapters are devoted to stabilization by a dissipative mechanism localized either inside the domain, or at the boundary. In both cases, general stabilization results are proven for large classes of damping terms, and the rate of decay of the energy is specified when it is linear, or behaves like a power of the velocity for small velocities in the region where it is effective. Those decay estimates can be viewed as a generalization of Enrique Zuazua's results on energy decay obtained in his thesis, but the proofs require new tools which had to be devised in connection with the techniques of control theory.

This monograph, completed by a substantial list of references, will be quite useful to those researchers who need to have at hand a synthetic overview of both basic and more recent results in the field of control and stabilization of wave equations in bounded domains.

\vspace{\baselineskip}
\begin{flushright}\noindent
May, 2021 \hfill 
\\Alain Haraux
\\CNRS, French National Centre for Scientific Research
\\INSMI, Institut National des Sciences Mathématiques et de leurs Interactions
\end{flushright}
\vspace{\baselineskip}
Previous Foreword (1990). \hfill
\vspace{\baselineskip}
\\The series ``Texts of Mathematical Methods'' plays a crucial role in serving the field of Mathematics, particularly in its intersection with the natural sciences, targeting postgraduate education and research. Each volume within this series aims to succinctly articulate the core principles of a specific area of interest, catering to the training needs of specialists. These volumes take the form of monographs, delving into a single topic or a connected set of significant investigations conducted by a particular individual.
\vspace{\baselineskip}
\\In essence, the ``Texts of Mathematical Methods'' are designed to address specialized subjects or offer ``unconventional'' perspectives on classic topics. They constitute a bridge between passive study and creative comprehension, facilitating a deeper understanding of mathematical theories. It is important to note that these expositions do not constitute exhaustive surveys or concise treatments of recent research findings. Instead, they serve as instructional tools, with authors striving to illuminate the distinctive characteristics of the theory. These texts act as supplementary resources for courses, seminars, or individual studies, presenting definitions, logical interrelationships between topics, and stating main results and corollaries. While some proofs may be omitted, the authors provide indications of the proof's nature, often accompanied by a simplified outline when technical details are replaced by bibliographic references.
\vspace{\baselineskip}
\\To summarize, ``Texts of Mathematical Methods'' serve as a reservoir of material presented in a semi-formal style, serving as a precursor to formal texts or monographs. They hold immediate relevance, whether due to their treatment of practical applications or the fact that the mathematical topics are actively employed or closely tied to applications.
\vspace{\baselineskip}
\\Through a commitment to fast and cost-effective publication,``Texts of Mathematical Methods'' make contemporary and relevant material accessible to practitioners of Mathematics, mathematicians with an interest in applications, and graduate students.
\vspace{\baselineskip}

\noindent {\it Editorial Board} 
\medskip

\noindent Alvercio Moreira Gomes, Gustavo Alberto Perla Menzala, Luiz Adauto Medeiros, Manuel A. Milla Miranda

\vspace{\baselineskip}
This manuscript presents the topics of the lectures given by Prof. Enrique Zuazua on ``Exact Controllability and System Stabilization'', from July 15 to September 15, 1989, at the Instituto de Matemática of the UFRJ.
\vspace{\baselineskip}
\\We believe that, due to the precise and didactic form with which they were written, they will represent a substantial aid to researchers and doctoral students interested in this area of Partial Differential Equations.
\vspace{\baselineskip}
\\We would like to thank Professor Zuazua for accepting our invitation to join the Department of Mathematical Methods as a Visiting Professor and for having made the arduous commitment to write these notes.
\vspace{\baselineskip}
\\We also thank Ms. Miriam Carqueja for the excellent typing work.

\begin{flushright}\noindent
Rio de Janeiro, December 12, 1990\hfill {\it Carlos Frederico Vasconcellos}\\
\end{flushright}

%
%
%

\preface
These Notes originated from a course I delivered at the Institute of Mathematics of the Universidade Federal do Rio de Janeiro, Brazil (UFRJ) in July-September 1989. I had the honor of visiting this institution upon the invitation and hospitality of Prof. Carlos Frederico Vasconcellos, affectionately known as Fred to his colleagues and friends. The dynamic and friendly team led by Prof. Luiz Adauto Medeiros warmly welcomed me.

These Notes were initially published in 1989 in Spanish under the title \emph{``Controlabilidad Exacta y Estabilizaci\'on de la Ecuaci\'on de Ondas''} in the Lecture Notes Series of the Institute. Spanish and Portuguese are sister languages and we believed that mathematicians at the Institute and in Brazil, as well as Portuguese speakers more broadly, could read and use them without language barriers.

I am indebted to Fred, who, with this invitation in 1989, opened up numerous opportunities for me early in my career. This experience allowed me to learn about Brazil and become fascinated by the country, its culture, nature, and people. I can never thank him enough. My gratitude extends to the entire Brazilian mathematical community, which has offered and continues to offer abundant opportunities for foreign mathematicians and has witnessed remarkable qualitative and quantitative growth over the last three decades.

Francesca Bonadei, who with I had the pleasure and honor to cooperate in a number of publishing and editorial efforts in the Mathematical Sciences during the last twelve years,  recently invited me to translate these Notes into English for publication by Springer Nature in the UNITEXT book series, coordinated by Thomas Hempfling and Christoph Baumann, who were kind enough to support the publication. 

Despite the significant evolution of the topic over the last three decades, I believe that the text, with its synthetic presentation of fundamental tools in the field, remains valuable for researchers in the area, especially for younger generations. It is written from the perspective of the young mathematician I was when I authored the Notes, needing to learn many things in the process and, therefore, taking care to develop details often left to the reader or not readily available elsewhere.

These Notes were written one year after completing my PhD at the Universit\'e Pierre et Marie Curie in Paris and drafting the lectures of Professor Jacques-Louis Lions at Coll\`ege de France in the academic year 1986-1987, later published as a book in 1988. Parts of these Notes offer a concise presentation of content developed in more detail in that book, supplemented by work on the decay of dissipative wave equations carried out during my PhD under the supervision of Professor Alain Haraux in Paris.
I cannot conclude without extending my warmest thanks to Alain Haraux and Jacques-Louis Lions, whose significant, generous, and positive influence has greatly impacted my career and life.\vspace{\baselineskip}
\begin{flushright}\noindent
Erlangen, September 2023\hfill {\it Enrique Zuazua}\\
\end{flushright}
\vspace{\baselineskip}

{\bf Previous Preface (1990)}

In these notes, we delve into the methods developed for the study of exact controllability and stabilization of evolution equations, gathering significant contributions from various authors.

To maintain conciseness, our focus centers on the wave equation with Dirichlet-type boundary conditions. However, it is worth noting that the techniques and results presented herein are versatile, readily adaptable and generalizable to a spectrum of models, including the Schr\"odinger equation, various plate models, the system of linear elasticity, and Maxwell's equations.

The book unfolds across six chapters, each serving a distinct purpose. The introductory first chapter sets the stage by presenting and formulating problems related to exact controllability and stabilization. Following this, the second and third chapters  explore exact boundary and internal controllability of the linear wave equation, respectively, employing the HUM method (Hilbert Uniqueness Method) introduced by J.-L. Lions in his book \cite{lions1988controlabilite}. The fourth chapter presents an exact controllability result for the semilinear wave equation with internal control. Moving forward, the fifth chapter offers estimates of the rate of decay of the energy of solutions for wave equations featuring nonlinear dissipation. The concluding sixth chapter focuses on addressing the intricate problem of boundary stabilization for the semilinear wave equation.

Each chapter includes a dedicated section on literature, offering a comprehensive review of relevant works, and concludes with the presentation of intriguing open problems for further exploration. 

My heartfelt gratitude goes to the Institute of Mathematics of the Federal University of Rio de Janeiro (UFRJ), and particularly to C. F. Vasconcellos, whose support not only made my journey to Rio de Janeiro possible but also facilitated the creation of this book.

\vspace{\baselineskip}
\begin{flushright}\noindent
Madrid, December 1990\hfill {\it Enrique Zuazua}\\
\hfill\\
\end{flushright}

%
%
%
%

\extrachap{Acknowledgements}{
The completion of this English version is owed to Darlis Bracho Tudares, the Digital Communications \& Web Specialist of our Chair and Humboldt Professorship (DCN) and the FAU Research Center for Mathematics of Data (MoD) at Friedrich-Alexander-Universit\"at Erlangen-N\"urnberg, Germany. Darlis meticulously handled the English translations and typesettings. I extend my sincere gratitude to Darlis for her invaluable assistance.

I am also indebted to Alain Haraux, Meryem Kafremer, Louis T. Tebou, Zhongjie Han, Jon Asier B\'arcena-Petisco, Borjan Geshkovski, Yubiao Zhang, Ahmet Ozkan Ozer, and Qiong Zhang for their thorough reviews and numerous insightful suggestions, which significantly enhanced the content.

Additionally, I express my warm appreciation to Francesca Bonadei from Springer-Nature. Our collaboration, initiated in 2012, has been consistently fruitful, culminating in the creation of this book.

The preparation of this final version was made possible through the support of various funding sources, including the Alexander von Humboldt-Professorship program, the ModConFlex Marie Curie Action (HORIZON-MSCA-2021-DN-01), the COST Action MAT-DYN-NET, the DFG Transregio 154 Project ``Mathematical Modelling, Simulation and Optimization Using the Example of Gas Networks,'' as well as grants PID2020-112617GB-C22 and TED2021-131390B-I00 from MINECO (Spain). Furthermore, the Madrid Government - UAM Agreement for the Excellence of the University Research Staff within the framework of the V PRICIT (Regional Programme of Research and Technological Innovation) has played a crucial role in supporting this work.

\vspace{\baselineskip}
\begin{flushright}\noindent
Erlangen, December 2023\hfill {\it Enrique Zuazua}\\
\end{flushright}

}

\tableofcontents

\mainmatter

%
%

\chapter{Presentation and Formulation of Controllability and Stabilization Problems}
\label{chapter01}

\abstract{
In this chapter we introduce the problems of controllability and stabilization for the wave equation. With the aim of being pedagogic, and to avoid unnecessary technical difficulties at this stage, both problems are analyzed in the simplified case where the control (in open-loop or in feedback form) acts on the whole domain.}

\section{The Exact Controllability Problem}
\label{sec:I.1}
Let $ \Omega $ be an open and bounded domain of $\mathbb{R}^n , n \geq 1$, with boundary $\Gamma=\partial \Omega$ of class $C^2$.

We consider the wave equation with the Dirichlet-type boundary conditions,
\begin{eqnarray}\label{eqI01}
\begin{cases}
y^{\prime\prime} - \Delta y=h(x,t) & \text{in} \hspace{1cm} Q = \Omega \times (0,T)
\\ y = 0 & \text{in} \hspace{0.9cm} \Sigma = \Gamma \times (0,T)
\\ y(x, 0) = y^0(x), y^\prime(x,0) = y^1(x) & \text{in} \hspace{1cm} \Omega.
\end{cases}
\end{eqnarray}

Here and in the sequel we use the notation $ ^\prime = \partial \cdot / \partial t $ for the partial derivative with respect to time.

Under adequate regularity and compatibility conditions for the initial data $\{y^0(x),y^1(x)\}$ and the source term $h(x, t)$, the system (\ref{eqI01}) admits a unique solution $y = y(x, t)$ in the energy space $C([0,T] ; H^1_0 (\Omega)) \cap C^1 ([0,T] ; L^2(\Omega))$. See for example H. Brezis \cite{brezis1973operateurs}, A. Haraux \cite{haraux1987semi} and J.-L. Lions \cite{lions1988controlabilite}.

The energy of the system (\ref{eqI01}) is represented by
\begin{equation}\label{eqI02}
E(t) = \frac{1}{2}\int\limits_{\Omega} \left[ \left|\nabla y(x,t) \right|^2 + \left| y^{\prime}(x,t) \right|^2 \right] d x.
\end{equation}
In the homogeneous case where $h\equiv0$, multiplying the equation by $y^\prime$ and integrating by parts in $ \Omega $ we get
\begin{equation*}
\frac{dE}{dt} (t) = 0
\end{equation*}
which implies that the energy is conserved along the trajectory of solutions of \eqref{eqI01}.

 The following is an immediate consequence of this energy conservation law: any non-trivial solution of the unforced wave equation ($h \equiv 0$) does not reach the equilibrium state $\{0,0\}$ at any time (except when the initial data is trivial: $\{y^0,y^1\} \equiv \{0,0\}$).

The problem of exact controllability consists precisely of bringing the trajectories to equilibrium in a uniform time, independently of the initial data, through the action of a suitable external force or control (which, in the case under consideration is represented by the third member $h$ of the equation).

More precisely, the problem of exact controllability can be formulated as follows: Analyze the existence of a time $T> 0$ such that for each pair of initial data $ \{ y^0, y^1 \} $ there is a control $h = h (x, t)$ such that the solution $y = y (x, t)$ of (\ref{eqI01}) satisfies
\begin{equation}\label{eqI03}
y(T) = y^\prime(T) = 0.
\end{equation}

This problem is slightly ambiguous since we have not specified neither the time-horizon $[0,T]$ available to control the system, nor the sets or functional spaces of initial data and controls. Obviously, in general, the possibility of obtaining an exact controllability result will depend on how large the time $T$ is, and how the space of the initial data and controls are chosen.

In the previous formulation of the problem we have not imposed any restrictions on the support of the control $h$, which makes its solution trivial, to some extent. Indeed, let $T>0$ and $ \{ y^0, y^1 \} \in H^1_0 (\Omega) \times L^2(\Omega) $ be arbitrary. Let us consider two functions $ a = a(t), b = b(t) \in C^2 ([0,T]) $ such that
\begin{eqnarray*}
a(0) = 1, a^\prime (0) = 0 ; \, a(T) = a^\prime(T) = 0; \, 
b(0) = 0, b^\prime (0) = 1 ; \,  b(T) = b^\prime(T) = 0.
\end{eqnarray*}
Then, the function
$
y(x,t) = a(t)y^0 (x) + b(t) y^1 (x)
$
satisfies
\begin{equation*}
y(0) = y^0, y^\prime (0) = y^1 ; \hspace{0.3cm} y(T) = y^\prime (T) = 0.
\end{equation*}

Therefore, the control
$
h = y^{\prime\prime} - \Delta y = a^{\prime\prime}(t)y^0 + b^{\prime\prime}(t)y^1 - a(t)\Delta y^0 - b(t)\Delta y^1
$
answers the question.

We have shown the following result: \emph{"Given an arbitrary $T>0$, for each initial data pair $ \{ y^0, y^1 \} \in H^1_0 (\Omega) \times L^2(\Omega) $ there is a control $ h \in C([ 0,T] ; H^{-2} (\Omega) ) $ such that the solution of (\ref{eqI01}) satisfies (\ref{eqI03})"}.

In fact, thanks to the linearity and time reversibility of the wave equation, the following more general result holds: Given an arbitrary $T>0$, for each $ \{ y^0, y^1 \}, \{ z^0, z^1 \} \in H^1_0 (\Omega) \times L^2(\Omega)$, there is a control $ h \in C^2([ 0,T] ; H^{-2} (\Omega) ) $ such that the solution of (\ref{eqI01}) satisfies
\begin{equation}\label{eqI04}
y(T) = z^0, y^\prime (T) = z^1.
\end{equation}

Indeed, let $ z = z (x, t) \in C([0,T] ; H^1_0 (\Omega)) \cap C^1 ([0,T] ; L^2(\Omega))$ solve the following system
\begin{eqnarray}\label{eqI05}
\begin{cases}
z^{\prime\prime} - \Delta z = 0 & \text{in} \hspace{1cm} Q
\\ z = 0 & \text{on} \hspace{0.9cm} \Sigma
\\ z (T) = z^0, z^\prime (T) = z^1 & \text{in} \hspace{1cm} \Omega
\end{cases}
\end{eqnarray}
and let us define the new state
\begin{equation}\label{eqI07}
\xi = y - z.
\end{equation}

Then, $y$ satisfies (\ref{eqI01}) if and only if $\xi$ satisfies the following system
\begin{eqnarray}\label{eqI08}
\begin{cases}
\xi^{\prime\prime} - \Delta \xi = h & \text{in} \hspace{0.3cm} Q
\\ 
\xi =0 & \text{in} \hspace{0.3cm} \Sigma
\\ 
\xi (0) = y^0 - z(0), \xi^\prime (0) = y^1 - z^\prime (0)  & \text{in} \hspace{0.3cm} \Omega
\end{cases}
\end{eqnarray}
and (\ref{eqI04}) holds if and only if
\begin{equation}\label{eqI09}
\xi (T) = \xi^\prime (T) = 0.
\end{equation}

According to the previous result, as $ \{ y^0 - z(0), y^1 - z^\prime (0) \} \in H^1_0 (\Omega) \times L^2(\Omega) $, there is a control $ h \in C([ 0,T] ; H^{-2} (\Omega) ) $ such that the solution of (\ref{eqI08}) satisfies (\ref{eqI09}), which implies that the solution of (\ref{eqI01}) satisfies (\ref{eqI04}).

Therefore, showing that every initial state can be driven into equilibrium is equivalent to showing that every initial state can be driven into every final state.

It is important to bear in mind that in the proof of this equivalence we have used the linearity and time-reversibility of the system in a crucial way.

We have seen that the problem of exact controllability with the internal control distributed through the whole domain $ \Omega $ can be treated in a straightforward manner,  and that controllability is achieved in an arbitrarily short time $T > 0$.

However, the systems in which the control acts in a subdomain or on the boundary of the domain $ \Omega $ where waves propagate are much more delicate and they require important further developments, which will be presented along this book.
As we shall see, the methods we shall develop will allow us to show that the regularity of the control above can be significantly improved and that the controllability of finite energy solutions can be achieved with controls in $L^2 (\Omega \times (0,T)  ) $.

\subsection{Localized Internal Control}
\label{subsec:I.1-a}
Let $ \omega \subset \Omega $ be an open and nonempty subdomain of $ \Omega $. Let us denote by $ \chi_\omega $ the characteristic function of $ \omega $, and consider the wave equation
\begin{eqnarray}\label{eqI10}
\begin{cases}
y^{\prime\prime} - \Delta y = h \chi_\omega & \text{in} \hspace{0.3cm} Q
\\ y = 0 & \text{on} \hspace{0.3cm} \Sigma
\\ y (0) = y^0, y^\prime (0) = y^1 & \text{in} \hspace{0.3cm} \Omega.
\end{cases}
\end{eqnarray}

The exact controllability problem is formulated in an analogous way: "Find a time $ T>0 $ such that for each pair of initial data $ \{ y^0, y^1 \} $ there is a control $ h $ such that the solution of (\ref{eqI10}) satisfies (\ref{eqI03})".

In this case, the control $ h $ acts only on the subdomain $ \omega $. It is therefore a problem of exact controllability in which the control is internal and localized.

Due to the finite speed of propagation of waves ($= 1$ in our model), in this setting, exact controllability requires the time $ T $ to be sufficiently large (so that the action of the control on $ \omega $ propagates to the whole $ \Omega $). Furthermore, the smaller $ \omega $ is, the longer the controllability time $T$ will be.

\subsection{Boundary Control}
\label{subsec:I.1-b}
Let $ \Gamma_0 \subset \Gamma $ be a nonempty open subset of $ \Gamma $ and consider the following wave equation with non-homogeneous boundary conditions:
\begin{eqnarray}\label{eqI11}
\begin{cases}
y^{\prime\prime} - \Delta y = 0 \hspace{0.56cm}\hspace{1cm} \text{in} \hspace{0.3cm} Q
\\ y = \begin{cases}
v & \hspace{1cm}\hspace{1cm} \text{on} \hspace{0.3cm} \Sigma_0 = \Gamma_0 \times (0,T)
\\0 & \hspace{1cm}\hspace{1cm} \text{on} \hspace{0.3cm} \Sigma_1 = \Sigma \setminus \Sigma_0
\end{cases}
\\ y(0) = y^0, y^\prime (0) = y^1 \hspace{0.3cm} \text{in} \hspace{0.3cm} \Omega.
\end{cases}
\end{eqnarray}

The problem of exact boundary controllability is formulated as follows: "Find $ T> 0 $ such that for each pair of initial data $ \{ y^0, y^1 \} $ there is a control $ v $ so that the solution of (\ref{eqI11}) satisfies (\ref{eqI03})".

Again, due to the finite speed of propagation, we can only expect to achieve exact boundary controllability if $ T $ is large enough.

As we will see in the next two chapters, both the problem of localized internal control and that of boundary control are complex. We will solve them using the Hilbert Uniqueness Method (HUM) introduced by J.-L. Lions in \cite{lions1986controlabilite}, \cite{lions1988exact},  \cite{lions1988controlabilite}. This method, in addition to being very flexible, has the advantage of systematically providing the optimal control, that is, the control minimizing a certain norm (normally the $L^2$-one, localized in the support of the controls).

Similar questions arise for semilinear wave equations. In Chapter \ref{chapter04}, we will present a method that combines the HUM method and a fixed point argument allowing us to solve the controllability problem when the nonlinearity is globally Lipschitz.

\section{The Stabilization Problem}
\label{sec:I.2}
As we have seen in the previous section, the energy of the solutions of the uncontrolled wave equation (in the absence of control with $h\equiv 0$) is conserved. The stabilization problem consists of forcing the uniform exponential decay of the energy as $ t \rightarrow + \infty $ through the action of a feedback control $ h $.  The control is in feedback or closed-loop form  when its value at each time can be determined by  the solution of the system itself at the same time instant. Its goal is to   induce dissipative effects and enhance the decay of solutions towards the null equilibrium as $ t \rightarrow + \infty $.

More precisely, we look for a function
$
F: \Omega \times \mathbb{R} \times \mathbb{R} \rightarrow \mathbb{R}
$
such that the energy of the solution of
\begin{eqnarray}\label{eqI12}
\begin{cases}
y^{\prime\prime} - \Delta y = F(x, y, y^\prime) & \text{in} \hspace{0.2cm} \Omega \times (0, \infty)
\\ y = 0 & \text{on} \hspace{0.2cm} \Gamma \times (0, \infty)
\\ y(0) = y^0, y^\prime(0) = y^1 & \text{in} \hspace{0.2cm} \Omega
\end{cases}
\end{eqnarray}
decays exponentially and uniformly as $ t \to + \infty $, i.e. there exist constants $ C>1 $ and $ \gamma>0 $ with
\begin{equation}\label{eqI13}
E(t) \leq CE(0) e^{-\gamma t} \hspace{1cm} \forall t >0
\end{equation}
holds for every solution of (\ref{eqI12}).

In order to motivate the type of feedback functions $F$ producing dissipative effects, observe that, formally, the energy dissipation law is as follows (all of this will be rigorously studied in each example):
\begin{equation*}
\frac{dE(t)}{dt} = \int_\Omega F(x, y , y^\prime)y^\prime dx.
\end{equation*}

Therefore, a first candidate is
$
F(x, y , y^\prime) = {-}a(x) y^\prime,
$
 $ a= a(x) $ being a non-negative function. This leads to the following system
\begin{eqnarray}\label{eqI14}
\begin{cases}
y^{\prime\prime} - \Delta y + a(x) y^\prime = 0 & \text{in}  \hspace{0.2cm} \Omega \times (0, \infty)
\\ y = 0 & \text{on} \hspace{0.2cm} \Gamma \times (0, \infty)
\\ y(0) = y^0, y^\prime(0) = y^1 & \text{in}  \hspace{0.2cm} \Omega.
\end{cases}
\end{eqnarray}

In the simplest case where
\begin{equation}\label{eqI15}
\exists a_0 , a_1 >0: a_0 \leq a(x) \leq a_1 \hspace{1cm} \forall x \in \Omega
\end{equation}
i.e. when the damping term is bounded and effective everywhere in the domain $\Omega$, the proof of the exponential decay is very simple. We present here a sketch.

Given $ \varepsilon >0 $, we define the following $\varepsilon-$perturbation of the energy
\begin{equation}\label{eqI16}
E_\varepsilon (t) = E(t) + \varepsilon \int_\Omega yy^\prime dx.
\end{equation}
First, we observe that there exists $ \varepsilon_0 >0 $ such that, if $ \varepsilon \leq \varepsilon_0 $, with $\varepsilon_0 $ small enough, then
\begin{equation}\label{eqI17}
\frac{1}{2} E_\varepsilon (t) \leq E(t) \leq 2E_\varepsilon (t) \hspace{1cm} \forall t>0
\end{equation}
for any solution (simply take $ \varepsilon_0 = \sqrt{\lambda_1}/{2} $ where $ \lambda_1 $ is the first eigenvalue of $ {-}\Delta $ in $ H^1_0 (\Omega)$).

Then,
\begin{align*}
\begin{split}
\frac{dE_\varepsilon}{dt}(t) &= \frac{dE}{dt}(t) + \varepsilon \frac{d}{dt} \left( \int_\Omega yy^\prime dx \right) = {-} \int_\Omega a(x) |y^\prime|^2 dx + \varepsilon \frac{d}{dt} \left( \int_\Omega yy^\prime dx \right).
\end{split}
\end{align*}
Moreover
\begin{align*}
\begin{split}
\frac{d}{dt} \left( \int_\Omega yy^\prime dx \right) & = \int_\Omega |y^\prime|^2 dx + \int_\Omega yy^{\prime\prime} dx
 = \int_\Omega |y^\prime|^2 dx + \int_\Omega y( \Delta y - a(x) y^\prime ) dx
\\ & = \int_\Omega |y^\prime|^2 dx - \int_\Omega |\nabla y|^2 dx - \int_\Omega a(x)yy^\prime dx.
\end{split}
\end{align*}

Therefore,
\begin{equation}\label{eqI18}
\frac{dE_\varepsilon}{dt}(t) = \int_\Omega \left( \varepsilon - a(x) \right) |y^\prime|^2 dx - \varepsilon \int_\Omega |\nabla y|^2 dx - \varepsilon \int_\Omega a(x)yy^\prime dx.
\end{equation}

On the other hand, one can estimate the last term in (\ref{eqI18}) as follows
\begin{align}\label{eqI19}
\begin{split}
\left| \int_\Omega a(x)yy^\prime dx \right| & \leq a_1 \int_\Omega |y| |y^\prime| dx
 \leq a_1 \left( \frac{\lambda_1}{2a_1} \int_\Omega |y|^2 dx + \frac{a_1}{2\lambda_1} \int_\Omega |y^\prime|^2 dx \right)
\\ & \leq \frac{1}{2} \int_\Omega |\nabla y|^2 dx + \frac{a^2_1}{2\lambda_1} \int_\Omega |y^\prime|^2 dx,
\end{split}
\end{align}
and therefore
\begin{equation}\label{eqI20}
\frac{dE_\varepsilon}{dt}(t) \leq \int_\Omega \left( \varepsilon \left(1 + \frac{a^2_1}{2 \lambda_1} \right) - a(x) \right) |y^\prime|^2 dx - \frac{\varepsilon}{2} \int_\Omega |\nabla y|^2 dx.
\end{equation}
Choosing
$
\varepsilon < \varepsilon_1 = (a_0 \lambda_1)/(2 \lambda_1 + a^2_1)
$,
(\ref{eqI20}) reduces to
\begin{equation*}
\frac{dE_\varepsilon}{dt}(t) \leq - \frac{a_0}{2} \int_\Omega |y^\prime|^2 dx - \frac{\varepsilon}{2} \int_\Omega |\nabla y|^2 dx \leq {-} \varepsilon E(t) \leq {-} \frac{\varepsilon}{2} E_\varepsilon (t).
\end{equation*}
Hence
\begin{equation}\label{eqI21}
E_\varepsilon (t) \leq E_\varepsilon (0)e^{-\frac{\varepsilon}{2}t} \hspace{1cm} \forall t >0
\end{equation}
and by (\ref{eqI17}),
\begin{equation}\label{eqI22}
E(t) \leq 4E(0)e^{-\frac{\varepsilon}{2}t} \hspace{1cm} \forall t >0.
\end{equation}
We have shown that: \emph{"If $ a = a(x) $ verifies (\ref{eqI15}), there are constants $ C>1 $ and $ \gamma > 0 $ such that every solution of (\ref{eqI14}) satisfies (\ref{eqI13}).}

The assumption that $a(x)\geq a_0 > 0$ everywhere in $\Omega$
corresponds to the case where the dissipation is effective in the entire domain $ \Omega $. The situation is much more delicate when the dissipation is effective only in a subset $ \omega $ of $ \Omega $ or when it is injected through the boundary conditions.

\subsection{Localized Internal Dissipation}
\label{subsec:I.2-a}
Let $ \omega \subset \Omega $ be a nonempty open subset and suppose that $ a \in L^\infty (\Omega) $ satisfies
\begin{equation}\label{eqI23}
a(x) \geq 0 \hspace{0.5cm} \forall x \in \Omega ;\hspace{0.3cm} a(x) \geq a_0 > 0 \hspace{0.5cm} \forall x \in \omega.
\end{equation}
This case corresponds to the dissipation only effective in $ \omega $.

The stabilization problem is formulated in an analogous way: \emph{Are there constants $ C>1 $ and $ \gamma >0 $ such that the solutions of (\ref{eqI14}) satisfy (\ref{eqI13})?}

As we shall see, when $ \omega$ is a measurable set with positive measure, all solutions tend to zero as time tends to infinity.  But, the exponential convergence requires  suitable and substantial geometric conditions on the set $ \omega $ where the damping is effective. Determining the conditions under which such an exponential decay holds, and developing the analytical methods to prove it, is a challenging problem that we shall address later on.

Similar questions arise for semilinear wave equations of the form
\begin{equation*}
y^{\prime\prime} - \Delta y + f(y) + a(x) y^\prime = 0.
\end{equation*}

We can even consider equations with nonlinear dissipative terms of the type
\begin{equation*}
y^{\prime\prime} - \Delta y + a(x) g( y^\prime )= 0.
\end{equation*}

However, when the damping enters nonlinearly in the system, the uniform exponential decay of the energy cannot be expected, in general. The rate of decay depends typically on the nature of the nonlinearity and, in particular, on its behavior when the velocity is small, i. e. $ y' \sim 0 $, which determines the effective strength of the damping term.

\subsection{Boundary Dissipation}
\label{subsec:I.2-b}
Let $ \Gamma_0 \subset \Gamma $ be a nonempty open subset of the boundary, and $ \Gamma_1 := \Gamma \setminus \Gamma_0 $ be its complement. Let us consider the free wave equation with homogeneous Dirichlet-type boundary condition on $ \Gamma_1 $:
\begin{eqnarray}\label{eqI24}
\begin{cases}
y^{\prime\prime} - \Delta y = 0 & \text{in} \hspace{0.3cm} \Omega \times (0, \infty)
\\
y = 0 & \text{on} \hspace{0.3cm} \Gamma_1 \times (0, \infty).
\end{cases}
\end{eqnarray}

Note that in (\ref{eqI24}) we do not impose so far any boundary condition on $ \Gamma_0 $. The stabilization problem consists, precisely, of defining dissipative boundary conditions on $ \Gamma_0 \times (0, \infty) $, so that the uniform exponential decay of the energy is achieved, that is, (\ref{eqI13}) for some constants $ C>1 $ and $ \gamma > 0 $.

To design suitable dissipative laws, formally, we compute the time-derivative of the energy along the solutions of (\ref{eqI24}) and impose the necessary boundary conditions on $\Gamma_0$ to guarantee its dissipativity.

We have
\begin{equation*}
\frac{dE}{dt}(t) = \int_{\Gamma_0} \frac{\partial y}{\partial v} y^\prime d\sigma,
\end{equation*}
and, therefore, the most natural boundary conditions are of the type
\begin{eqnarray}\label{eqI25}
\frac{\partial y}{\partial v} = -a(x) y^\prime \hspace{0.5cm} \text{over} \hspace{0.5cm} \Gamma_0 \times (0, \infty)
\end{eqnarray}
with $ a = a(x) > 0 \text{ on } \Gamma_0 $.

Having chosen a closed-loop system of the form (\ref{eqI24})-(\ref{eqI25}), the problem is to obtain conditions on the partition $ \{ \Gamma_0, \Gamma_1 \} $ and the function $ a=a(x) $ so that (\ref{eqI13}) holds.

Of course, these questions can also be posed in a more general  framework of semilinear wave equations.
\smallskip

\subsection{Plan}

As we will see in the following chapters, both the problem of exact controllability, thanks to the application of the HUM method, and the problem of stabilization, are solved by obtaining adequate a priori estimates on solutions, allowing to estimate the total energy of the solutions by the energy localized in the region (inside the domain or on its boundary), where the action of the control or dissipative mechanism is injected. These estimates will be referred to as "observability inequalities". This terminology is justified by the fact that observing the energy localized in a subset of the domain or of its boundary suffices to get estimates of the energy everywhere in the domain where waves propagate. To obtain these estimates multiplier techniques will be used, in combination with unique continuation principles.

Chapters \ref{chapter02}, \ref{chapter03} and \ref{chapter04} are devoted to controllability problems, considering first the boundary controllability of the wave equation, then the case where the control acts on  the interior of the domain where waves propagate (on an neighbourhood of its boundary, more precisely), and, finally, the semilinear wave equation.   In  Chapter \ref{chapter05} we will study the problem of stabilization, considering the case of the nonlinear dissipation distributed over the entire domain $ \Omega $, while Chapter \ref{chapter06} is devoted to the case of boundary dissipation. We will not address the problem of localized internal dissipation, for which we refer to A. Haraux \cite{haraux1989remarque} and E. Zuazua \cite{zuazua1990exponential}, \cite{zuazua1991exponential} and the complementary bibliography at the end of these Notes.


%
%

\chapter{Boundary Controllability of the Linear Wave Equation}
\label{chapter02} 

\abstract{
In this chapter we prove the exact boundary controllability for the wave equation, considering first the constant coefficients case, and then  a more general class of space-dependent variable coefficients. The proof is based on the HUM method, that transforms the problem into an observability one, and the multiplier method, which allows to derive the equivalent observability inequalities.
}

\section{Description of the HUM Method}
\label{sec:II.1}
Let $ \Omega $ be a bounded domain of $\mathbb{R}^n, n\geq 1$, with boundary $\Gamma$ of class $C^2$. Let $T>0$ and $\Gamma_0 \subset \Gamma$ be a nonempty open subset of $\Gamma$.

Consider the following wave equation with inhomogeneous boundary conditions:
\begin{eqnarray}\label{eqII01}
\begin{cases}
y^{\prime\prime} - \Delta y = 0  \hspace{1.57cm} \text{in} \hspace{0.4cm} Q = \Omega \times (0, T)
\\ y = \begin{cases}
v & \hspace{2cm}  \text{on} \hspace{0.3cm} \Sigma_0 = \Gamma_0 \times (0, T)
\\ 0 & \hspace{2cm} \text{on} \hspace{0.3cm} \Sigma_1 = ( \Gamma \setminus \Gamma_0) \times (0, T)
\end{cases}
\\ y(0) = y^0, y^\prime(0) = y^1 \hspace{0.35cm} \text{in} \hspace{0.3cm} \Omega.
\end{cases}
\end{eqnarray}

We analyze exact controllability of this system for the set of initial data in $ L^2(\Omega) \times H^{-1}(\Omega) $,  with controls in $ L^2(\Sigma_0) $. The goal is then to tackle the following issue: \emph{Find $ T>0 $ such that for each pair of initial data $ \{ y^0, y^1 \} \in L^2(\Omega) \times H^{-1}(\Omega) $, there is a control $ v \in L^2(\Sigma_0) $ such that the solution of (\ref{eqII01}) satisfies}
\begin{equation}\label{eqII02}
y(T) = y^\prime (T) = 0.
\end{equation}

\label{remark:II.2.1}
\begin{remark}
Due to the finite speed of propagation of waves (which is unit for (\ref{eqII01})), to achieve exact controllability, it is necessary for $T$ to be large enough, that is, $ T > T_0 $, with $ T_0 >0 $ depending on $ \Omega $ and $ \Gamma_0 $.

As we shall see, in the chosen functional framework, exact controllability can only hold when $ \Gamma_0 $ is a sufficiently large subset of $ \Gamma. \quad\quad\hfill \square $
\end{remark}
\begin{remark}
\label{remark:II.2.2}
Recalling the results of Chapter \ref{chapter01}, thanks to the linearity and time-reversibility of the wave equation, the problem can be formulated in the following equivalent way: \emph{Find $ T >0 $ such that for every pair of initial and final data $ \{ y^0, y^1 \}, \{ z^0, z^1 \} \in L^2(\Omega) \times H^{-1}(\Omega) $, there is a control $ v \in L^2(\Sigma_0) $ such that the solution of (\ref{eqII01}) satisfies}
\begin{equation}\label{eqII03}
y(T) = z^0, y^\prime(T) =z^1.
\end{equation}

This type of problem has been studied by numerous authors and the literature is broad, going back to the  celebrated survey article \cite{russell1978control} by  D. L. Russell, back in 1978, where the author made a compilation of the most relevant methods and results available at that time. Although, back then,    many important results were available, a systematic method for treating these problems was lacking. J.-L. Lions \cite{lions1986controlabilite} in 1986 proposed the so-called Hilbert Uniqueness Method (HUM) that allows to reduce the problem of exact controllability to obtaining adequate a priori estimates. This led to a rapid and fruitful development of the field.

Throughout this chapter, we will closely follow Chapter I of the book by J.-L. Lions \cite{lions1988controlabilite}, using the HUM method. All the results presented here, except those corresponding to Section 2.7, in which  variable coefficients are addressed, have been taken from \cite{lions1988controlabilite}.

The adaptation of the HUM method to the exact controllability of (\ref{eqII01}) is as follows:
Given $ \{ \phi^0, \phi^1 \} \in \mathcal{D}(\Omega) \times \mathcal{D}(\Omega) $, we solve the system
\begin{eqnarray}\label{eqII04}
\begin{cases}
\phi^{\prime\prime} - \Delta \phi = 0 & \text{in}  \hspace{0.2cm} Q
\\ \phi = 0 & \text{on} \hspace{0.2cm} \Sigma 
\\ \phi(0) = \phi^0, \phi^\prime(0) = \phi^1& \text{in} \hspace{0.2cm} \Omega,
\end{cases}
\end{eqnarray}
which admits a unique solution in the energy space. Out of the solution $\phi$ of \eqref{eqII04} we consider its normal derivative over the subset of the controlled boundary and solve the backward and controlled wave equation
\begin{eqnarray}\label{eqII05}
\begin{cases}
y^{\prime\prime} - \Delta y = 0 \hspace{1cm} \text{in} \hspace{0.2cm} Q
\\ y(T) = y^\prime(T) = 0 \hspace{0.4cm} \text{in} \hspace{0.2cm} \Omega
\\ y = \begin{cases}
\partial \phi / \partial \nu \hspace{1cm} \text{on} \hspace{0.2cm} \Sigma_0
\\ 0 \hspace{1.77cm} \text{on} \hspace{0.2cm} \Sigma_1.
\end{cases}
\end{cases}
\end{eqnarray}

In (\ref{eqII05}), $ \nu = \nu(x) $ denotes the unit normal vector pointing outside $ \Omega $ at point $ x \in \Gamma \text{ and } \partial \cdot /\partial \nu $ the outward normal derivative.

Problem (\ref{eqII05}), again, admits an unique weak solution whose regularity will be discussed below.

We define the operator
\begin{equation}\label{eqII06}
\Lambda \{ \phi^0, \phi^1 \} = \{ y^\prime (0), {-}y(0) \}.
\end{equation}

Multiplying the first equation (\ref{eqII05}) by $ \theta = \theta(x,t) $, another solution $\theta$ of (\ref{eqII04}) with initial data $ \{ \theta^0, \theta^1 \} $,  i. e.
\begin{eqnarray}\label{eqII04b}
\begin{cases}
\theta^{\prime\prime} - \Delta \theta = 0 & \text{in}  \hspace{0.2cm} Q
\\ \theta= 0 & \text{on} \hspace{0.2cm} \Sigma 
\\ \theta = \theta^0, \theta^\prime(0) = \theta^1& \text{in} \hspace{0.2cm} \Omega,
\end{cases}
\end{eqnarray}
 and integrating by parts, formally for the time being, we get, 
\begin{equation}\label{eqII07}
\langle \Lambda \{ \phi^0, \phi^1 \}, \{ \theta^0, \theta^1 \} \rangle = \int_{\Sigma_0} \frac{\partial \phi}{\partial \nu} \hspace{0.1cm} \frac{\partial \theta}{\partial \nu} d\Sigma
\end{equation}
where $ d\Sigma = d\Gamma dt $ represents the measure on the lateral surface $ \Sigma $ of the cylinder $ Q $.

In particular, when $\theta = \phi$, we have
\begin{equation}\label{eqII08}
\langle \Lambda \{ \phi^0, \phi^1 \}, \{ \phi^0, \phi^1 \} \rangle = \int_{\Sigma_0} \left| \frac{\partial \phi}{\partial \nu} \right|^2 d\Sigma .
\end{equation}

Suppose that $ T $ and $ \Gamma_0 $ are such that the following uniqueness or unique continuation result holds:
\begin{eqnarray}\label{eqII09}
\begin{cases}
\text{If } \phi \text{ solving (\ref{eqII04}) satisfies } \partial \phi / \partial \nu=0 \text{ over } \Sigma_0,
\\ \text{then } \phi \equiv 0 \text{ and therefore } \{ \phi^0, \phi^1 \} = \{ 0,0 \}.
\end{cases}
\end{eqnarray}

In this case, the map
\begin{equation}\label{eqII10}
\left\| \{ \phi^0, \phi^1 \} \right\|_F = \left( \int_{\Sigma_0} \right| \frac{\partial \phi}{\partial \nu} \left|^2 d\Sigma \right)^\frac{1}{2}
\end{equation}
defines a norm on $ \mathcal{D}(\Omega) \times \mathcal{D}(\Omega) $. This allows introducing the following Hilbert space:
\begin{equation}\label{eqII11}
\begin{split}
F = &\text{ the Hilbert space defined as the completion of } 
\\ &\mathcal{D}(\Omega) \times \mathcal{D}(\Omega) \text{ with respect to the norm } \| \cdot \|_F.    
\end{split}
\end{equation}
Thanks to (\ref{eqII07}), the operator $\Lambda$ can be extended to a continuous linear operator from $F$ into its dual $F^\prime$. On the other hand, from (\ref{eqII10}) we deduce that
\begin{equation}\label{eqII12}
\Lambda : F \longrightarrow F^\prime \text{ is an isomorphism.}
\end{equation}

Therefore, given initial data $ \{ y^0, y^1 \} $ such that
\begin{equation}\label{eqII13}
\{ y^1, -y^0 \} \in F^\prime,
\end{equation}
the problem
\begin{equation}\label{eqII14}
 \Lambda \{ \phi^0, \phi^1 \} =  \{ y^1, -y^0 \}
\end{equation}
admits a unique solution $\{ \phi^0, \phi^1 \} \in F$.

But solving (\ref{eqII14}) is equivalent to the fact that $ y = y(x,t) $ solution of (\ref{eqII05}), with $ \phi $ solution of (\ref{eqII04}) and $\{ \phi^0, \phi^1 \}$ solving (\ref{eqII14}), satisfies
\begin{equation}\label{eqII15}
y(0) = y^0; y^\prime(0) = y^1.
\end{equation}

The control sought is therefore given in the following way
$
v = \partial \phi/\partial \nu\big|_{\Sigma_0}
$
where $ \phi = \phi (x,t) $ is the solution of (\ref{eqII04}) corresponding to the data $\{ \phi^0, \phi^1 \} \in F $ satisfying (\ref{eqII14}).

Note that, by construction,
$
v \in L^2 ( \Sigma_0).
$
Therefore, we have shown the following result: {\it "If $\Omega$, $ \Gamma_0 \subset \Gamma $  and $T>0$ are such that the uniqueness result (\ref{eqII09}) is satisfied, for each pair of initial data $ \{ y^0, y^1 \} $ such that $ \{ y^1,-y^0 \} \in F^\prime $, there exists a control $ v \in L^2 ( \Sigma_0) $ such that the solution $y$ of (\ref{eqII01}) satisfies (\ref{eqII02})."}

Thanks to Holmgren's Theorem (see Th.5.5.3, \cite{hormander1976analysis} or section I.8 of \cite{lions1988controlabilite}), it is known that for every Lipschitz domain $\Omega$ and open and non-empty subset $ \Gamma_0 $ of $ \Gamma $, there exists $ T_0 = T_0 (\Gamma_0 , \Omega ) $ such that the uniqueness result above holds true if $ T > T_0 $.

This result we just presented is satisfactory since it applies to a broad class of subsets $ \Gamma_0 $ of the boundary $ \Gamma $. However, it is imprecise since it provides very little information about the space $ F^\prime $ of controlled initial data.
Therefore, the problem that remains to be solved is  the identification of the dual space $ F^\prime $ or, equivalently,  of the space $F$ itself. But since $F$ is the completion of $ \mathcal{D}(\Omega) \times \mathcal{D}(\Omega) $ with respect to the norm $\| \cdot \|_F $, the question is reduced to identifying this norm.
In particular, proving the exact controllability in the optimal space $ L^2(\Omega) \times H^{-1} (\Omega) $, is equivalent to
$
F^\prime = H^{-1} (\Omega) \times L^2(\Omega)
$
or, equivalently, to
$
F = H_0^1 (\Omega) \times L^2(\Omega),
$
which, in turn, is equivalent to the existence of constants $ c>0, C>0 $ such that
\begin{align} \label{eqII16}
\begin{split}
{ c {\left\| \{\phi^0, \phi^1\} \right\|}_{H_0^1 (\Omega) \times L^2(\Omega)} } \leq& \left\| \{\phi^0, \phi^1\} \right\|_F \leq C \left\| \{\phi^0, \phi^1\} \right\|_{H_0^1 (\Omega)\times L^2(\Omega)},
\\ &\quad \quad \quad  \quad \forall \{\phi^0, \phi^1\} \in \mathcal{D}(\Omega) \times \mathcal{D}(\Omega).
\end{split}
\end{align}

In the next two sections, we prove that (\ref{eqII16}) is satisfied if $ \Gamma_0 $ is a sufficiently large subset of $ \Gamma $ (in a sense that we specify later) and $ T>0 $ is large enough (depending on $ \Omega $ and $ \Gamma_0 $).
\end{remark}

\section{The Direct Inequality: Hidden Regularity}
\label{sec:II.2}
We begin by proving a fundamental identity for the solutions of
\begin{eqnarray}\label{eqII17}
\begin{cases}
\theta^{\prime\prime} - \Delta \theta = f \hspace{1cm} & \text{in} \hspace{0.2cm} Q
\\ \theta =0 \hspace{1cm} & \text{on} \hspace{0.2cm} \Sigma
\\ \theta (0) = \theta^0, \theta^\prime(0) = \theta^1 \hspace{1cm} & \text{in} \hspace{0.2cm} \Omega.
\end{cases}
\end{eqnarray}

\begin{lemma}
\label{lemma:II.2.1}
Let $ q = q(x) \text{ be a vector field in } ( C^1(\overline{\Omega}))^n $. For every solution of (\ref{eqII17}) with the initial data $ \{ \theta^0 , \theta^1, f \} \in \mathcal{D}(\Omega) \times\mathcal{D}(\Omega) \times \mathcal{D}(Q) $, the following identity holds
\begin{align} \label{eqII18}
\begin{split}
\frac{1}{2} \int_\Sigma (q \cdot \nu) \left| \frac{\partial \theta}{\partial \nu}  \right|^2 d\Sigma = & \int_\Omega \theta^\prime q \cdot \nabla \theta \bigg|_0^T + \frac{1}{2} \int_Q (\operatorname{div}q)(|\theta^\prime|^2 - |\nabla \theta|^2)
\\ &+  \int_Q \frac{\partial q_k}{\partial x_j} \hspace{0.1cm} \frac{\partial \theta}{\partial x_k} \hspace{0.1cm} \frac{\partial \theta}{\partial x_j} - \int_Q fq \cdot \nabla \theta.
\end{split}
\end{align}
\end{lemma}

\begin{proof}
\label{proof:II.lemma.2.1}
It is enough to multiply the equation (\ref{eqII17}) by $ q \cdot \nabla \theta $ and integrate by parts in $ Q$.
\end{proof}

\begin{remark}
\label{remark:II.2.3}
We denote by $ \cdot $ the Euclidean scalar product in $ \mathbb{R}^n $. In (\ref{eqII18})  the convention for the sum of repeated indices is used, and both the $ (x,t) $ variables and the symbols $ dx $ and $ dt $ in the integration signs, over $ \Omega $ and $ Q $ are omitted.
\end{remark}

The following notation has also been used
\begin{equation*}
    \operatorname{div} q = \text{divergence of } q = \sum_{k=1}^{n} \frac{\partial q_k}{\partial x_k}
\end{equation*}and, finally,
\begin{equation*}
\int_\Omega \theta^\prime q \cdot \nabla \theta \bigg|_0^T = \int_\Omega \theta^\prime (x,T) q(x) \cdot \nabla \theta (x,T) dx - \int_\Omega \theta^\prime (x,0) q(x) \cdot \nabla \theta (x,0) dx. \quad \hfill \square
\end{equation*}
From the identity (\ref{eqII18}) the following estimate is easily obtained.
\begin{theorem}
\label{theorem:II.2.1}
There is a constant $c>0$ such that
\begin{equation}\label{eqII19}
\int_\Sigma \left|\frac{\partial \theta}{\partial \nu} \right|^2 d\Sigma \leq c(1+T) \left({\left\|\theta^0 \right\|}_{H_0^1 (\Omega)}^2 + {\left\|\theta^1 \right\|}_{L^2 (\Omega)}^2 + {\left\| f \right\|}_{L^1(0,T; L^2(\Omega))}^2\right)
\end{equation}
for all $ T>0 $ and all solutions of (\ref{eqII17}) with the initial data
$$
\{ \theta^0, \theta^1, f \} \in H_0^1 (\Omega) \times L^2 (\Omega) \times L^1(0,T; L^2(\Omega)).
$$
\end{theorem}

\begin{proof}
\label{proof:II.theorem.2.1}
Let $ h = h(x) \in C^{1}( \overline{\Omega})^n $ be a vector field such that
$
h = \nu \text{ on } \Gamma.
$
It exists since $ \Omega $ is of class $ C^2 $ (see J.-L. Lions \cite{lions1988controlabilite}, p. 28 for its construction).

Applying the identity (\ref{eqII18}) with $ q = h $, we have
\begin{align}\label{eqII20}
\begin{split}
\frac{1}{2} \int_\Sigma \left|\frac{\partial \theta}{\partial \nu} \right|^2 d\Sigma \leq  C \big\{ &{\left\| \theta^\prime \right\|}_{L^\infty(0, T; L^2(\Omega))}^2 + {\left\|\theta \right\|}_{L^\infty (0, T; H_0^1 (\Omega))}^2
\\ & +  {\left\| \theta^\prime \right\|}_{L^2(Q)}^2 
+  {\left\| \nabla \theta \right\|}_{L^2(Q)}^2 + {\left\|\ f \right\|}_{L^1(0,T; L^2(\Omega))}^2 \big\}
\end{split}
\end{align}
with $ C = C ({\left\| q \right\|}_{W^{1,\infty}(\Omega)}) $.

Moreover,  classical regularity results for the wave equation ensure that 
\begin{equation}\label{eqII21}
 {\left\|(\theta, \theta^\prime) \right\|}_{L^\infty (0, T; H_0^1 (\Omega))\times L^2(\Omega))} \leq
C \left[ {\left\|\theta^0 \right\|}_{H_0^1 (\Omega)} + 
 {\left\|\theta^1 \right\|}_{L^2 (\Omega)} + 
{\left\| f \right\|}_{L^1(0,T; L^2(\Omega))} \right].
\end{equation}

Combining (\ref{eqII20}) and (\ref{eqII21}) gives (\ref{eqII19}) for solutions of (\ref{eqII17}) with regular data. The estimation is extended to solutions with the initial data $ \{  \theta^0, \theta^1, f \} \in H_0^1 (\Omega) \times L^2(\Omega) \times L^1(0,T; L^2(\Omega)) $ by a simple density argument.
\end{proof}

\begin{remark}
\label{remark:II.2.4}
The multiplier $ h \cdot \nabla\theta $ was introduced by F. Rellich in the framework of elliptic equations. The adaptation to the hyperbolic case and (\ref{eqII19}) are due to J.-L. Lions \cite{lions1983controle}. The same method makes it possible to prove an analogous estimate for equations with variable coefficients, that is,
$
\theta^{\prime\prime} - \operatorname{div}(a(x) \nabla \theta) = f
$
with $ a = a(x) \in W^{1, \infty}(\Omega) $ (c.f. \cite{lions1983controle}).

This type of estimates were extended to semilinear equations
\begin{equation*}
\theta^{\prime\prime} - \Delta \theta + g(\theta) = f
\end{equation*}
by J.-L. Lions \cite{lions1987hidden} and M. Milla-Miranda and L. A. Medeiros \cite{milla1988hidden}.
\end{remark}

\begin{remark}
\label{remark:II.2.5}
Inequality (\ref{eqII19}) provides a hidden regularity result. Indeed, the solutions of (\ref{eqII17}) with the initial data $ \{  \theta^0, \theta^1, f \} \in H_0^1 (\Omega) \times L^2(\Omega) \times L^1(0,T; L^2(\Omega)) $ belong to the class
$
\theta \in C ([0,T];H_0^1 (\Omega) ) \cap C^1([0,T]; L^2(\Omega)).
$
But this does not suffice to guarantee 
$
\partial \theta/\partial \nu \big|_\Sigma \in L^2(\Sigma),
$
a fact that we have shown to hold using multiplier arguments. We often refer to this added regularity property, as "hidden regularity", since it is due to the fine interaction of the wave equation satisfied inside the domain and the boundary condition. This shows that finite energy solutions of the wave equation with homogeneous Dirichlet boundary conditions constitute a very specific subclass of functions  in which such added regularity property is fulfilled.
$ \hfill \square $
\end{remark}

\section{The Inverse Inequality: Observability}
\label{sec:II.3}
\begin{figure}[h]
\centering
\captionsetup{justification=centering}
\includegraphics[width=8cm]{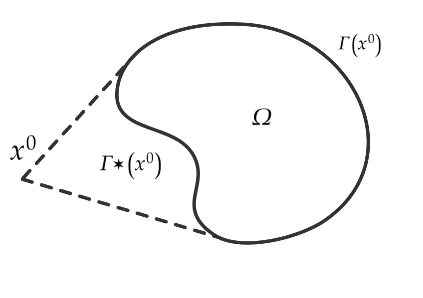}
\caption{The geometric configuration under consideration}
\label{fig:II02.1}
\end{figure}

Given an arbitrary point $ x^0 \in \mathbb{R}^n $, we define $ m(x) = x - x^0 $ and \begin{equation}\label{eqII21.2} R(x^0) = \| m \|_{L^{\infty}(\Omega)}\end{equation}
the radius of the smallest ball centered at $x^0$ and containing $\Omega$. We then split the boundary $ \Gamma $ of $ \Omega $ into two parts, as in Figure \ref{fig:II02.1},
\begin{align*}
\Gamma(x^0) = \{ x \in \Gamma: m(x) \cdot \nu(x) > 0 \}; \hspace{0.2cm}
 \Gamma_{\star}(x^0) = \{ x \in \Gamma: m(x) \cdot \nu(x)  \leq 0 \} = \Gamma \backslash \Gamma(x^0),
\end{align*}
and, similarly, the lateral boundary $\Sigma$:
\begin{equation*}
\Sigma(x^0) = \Gamma(x^0) \times (0,T); \hspace{0.3cm}
\Sigma_{\star}(x^0) = \Gamma_{\star}(x^0) \times (0,T) = \Sigma \backslash \Sigma(x^0).
\end{equation*}

Consider the homogeneous equation
\begin{equation}\label{eqII22}
\begin{cases}
\phi^{\prime \prime} - \Delta \phi = 0 & \text{in} \hspace{0.3cm} Q
\\ \phi = 0 & \text{on} \hspace{0.3cm} \Sigma
\\ \phi(0) = \phi^0, \phi^\prime(0) = \phi^1 & \text{in} \hspace{0.3cm} \Omega.
\end{cases}
\end{equation}
The following result holds.
\label{theo:II.2.2}
\begin{theorem}
\label{theorem:II.2.2}
Let $ T > 2R(x^0) $. Then, for every solution of (\ref{eqII22}) with the initial data $ \{ \phi^0, \phi^1 \} \in H_{0}^1 (\Omega) \times L^2(\Omega) $, the following observability inequality holds
\begin{align}\label{eqII23}
    E_0 = \frac{1}{2} \int_\Omega \left[ \left| \nabla \phi^0 (x) \right|^2 + \left| \phi^1 (x) \right|^2  \right] dx \leq \frac{R(x^0)}{2(T - 2R(x^0))} \int_{\Sigma(x^0)} \left| \frac{\partial \phi}{ \partial \nu} \right|^2 d\Sigma.
\end{align}
\end{theorem}

\begin{remark}
\label{remark:II.2.6}
Inequality (\ref{eqII23}) is often referred to as the "observability inequality" since it guarantees that the total energy of solutions can be estimated (or observed) out of the boundary trace of the normal derivative. These inequalities play also an important role in the theory of inverse problems. $ \hfill \square $
\end{remark}

\begin{proof}
\label{proof:II.remark.2.6}
First, applying (\ref{eqII18}) with $ f = 0, \theta = \phi $ and $ q(x) = m(x) = x - x^0 $ it follows
\begin{align}\label{eqII24}
\begin{split}
    \int_\Omega \phi^\prime m \cdot \nabla\phi \bigg|_0^T + \frac{n}{2} \int_Q \left| \phi^\prime \right|^2 + \left( 1 - \frac{n}{2} \right) \int_Q \left| \nabla \phi \right|^2 
    = & \frac{1}{2} \int_\Sigma (m \cdot \nu) 
    \left| \frac{\partial \phi}{ \partial \nu} \right|^2 d\Sigma
    \\ \leq & \frac{1}{2} \int_{\Sigma(x^0)} (m \cdot \nu) \left| \frac{\partial \phi}{ \partial \nu} \right|^2 d\Sigma
    \\ \leq & \frac{R(x^0)}{2} \int_{\Sigma(x^0)} \left| \frac{\partial \phi}{ \partial \nu} \right|^2 d\Sigma.
\end{split}
\end{align}

On the other hand,
\begin{align}\label{eqII25}
\begin{split}
\frac{n}{2} \int_Q \left| \phi^\prime \right|^2 + \left( 1 - \frac{n}{2} \right) \int_Q \left| \nabla \phi \right|^2 = & \frac{1}{2} \int_Q (\left| \phi^\prime \right|^2 + \left| \nabla \phi \right|^2 ) + \frac{n-1}{2} \int_Q (\left| \phi^\prime \right|^2 - \left| \nabla \phi \right|^2 ) 
\\ =& TE_0 + \frac{n-1}{2} \int_Q (\left| \phi^\prime \right|^2 - \left| \nabla \phi \right|^2 ) 
\end{split}
\end{align}
since the energy
$
E(t) = \frac{1}{2} \int_\Omega \left|  \phi^\prime (x, t) \right|^2 + \left| \nabla \phi (x,t) \right|^2 ) dx
$
is conserved along the trajectory, i.e.
$
E(t) = E_0 
$
for all $ t \in [ 0, T ]$.

Next, multiply the equation by $ \phi $ and integrate by parts to easily obtain the so-called "equipartition of energy" identity
\begin{equation} \label{eqII26}
\int_\Omega \phi^\prime \phi \Big|_0^T = \int_Q (\left| \phi^\prime \right|^2 - \left| \nabla \phi \right|^2 ).
\end{equation}

Now, combining (\ref{eqII24}), (\ref{eqII25}) and (\ref{eqII26}) leads to
\begin{equation} \label{eqII27}
X(t) \Big|_0^T + T E_0 \leq \frac{R(x^0)}{2} \int_{\Sigma(x^0)} \left| \frac{ \partial \phi}{ \partial \nu} \right|^2 d\Sigma
\end{equation}
with
\begin{equation} \label{eqII28}
X(t) = \int_\Omega \phi^\prime (x,t) \left(m(x) \cdot \nabla \phi (x,t) + \frac{n-1}{2} \phi (x,t)\right) dx.
\end{equation}
It is therefore enough to prove that
\begin{equation} \label{eqII29}
\left\| X(\cdot) \right\|_{L^\infty(0,T)} \leq R(x^0) E_0.
\end{equation}

For simplicity, the explicit dependence of $X$ on $t$ is omitted for the rest of the proof. We have
\begin{equation} \label{eqII30}
\left| X \right| \leq \frac{\delta}{2} \int_\Omega \left| \phi^\prime \right|^2 + \frac{1}{2\delta} \int_\Omega \left| m \cdot \nabla \phi + \frac{n-1}{2} \phi \right|^2 , \hspace{0.3cm} \forall \delta > 0.
\end{equation}

Furthermore,
\begin{align}\label{eqII31}
\begin{split}
 \int_\Omega \left| m \cdot \nabla \phi + \frac{n-1}{2} \phi \right|^2 \leq &\int_\Omega \left| m \right|^2 \left| \nabla\phi \right|^2 + 
 \left( \frac{n-1}{2} \right)^2 \int_\Omega \left| \phi \right|^2 
 \\ & +(n-1) \int_\Omega (m \cdot \nabla \phi) \phi.
\end{split}
\end{align} 

Now since $m(x) = x-x^0$
\begin{equation*}
\int_\Omega m \cdot \nabla\phi\phi = \frac{1}{2} \int_\Omega m \cdot \nabla (\phi^2) = {-} \frac{n}{2} \int_\Omega \phi^2
\end{equation*}
and, therefore,
\begin{align}\label{eqII32}
\begin{split}
\int_\Omega \left| m \cdot \nabla \phi + \frac{n-1}{2} \phi \right|^2  &\leq R^2(x^0) \int_\Omega \left| \nabla\phi \right|^2  + \left( \frac{(n-1)^2}{4} - \frac{n(n-1)}{2} \right) \int_\Omega \phi^2
\\ & \leq R^2 (x^0) \int_\Omega \left| \nabla\phi \right|^2.
\end{split}
\end{align}

Finally, from (\ref{eqII30}) and (\ref{eqII32})
\begin{equation*}
    |X| \leq \frac{\delta}{2} \int_\Omega |\phi^\prime|^2 + \frac{R^2 (x^0)}{2\delta} \int_\Omega \left| \nabla\phi \right|^2,
\end{equation*}
from which (\ref{eqII29}) is deduced by choosing $ \delta = R(x^0)$.
\end{proof}

\begin{remark}
\label{remark:II.2.7}
L. F. Ho in \cite{ho1986observabilite} established for the first time an estimate of the type (\ref{eqII23}), for $ T>T_1 $ with some $ T_1 > 2R(x^0) $. Then, J.-L. Lions in \cite{lions1988exact} improved the needed minimal time showing that, if $ T > 2R(x^0) $, the inequality
\begin{equation*}
    E_0 \leq C \int_{\Sigma(x^0)} \left| \frac{\partial \phi}{\partial \nu} \right|^2 d\Sigma
\end{equation*}
holds for a constant $ C>0$ independent of $\phi$. But in \cite{lions1988exact} the constant $C$ was not explicitly given since a "compactness-uniqueness" argument was used in the proof.
The estimate (\ref{eqII23}), as presented here, was proved by V. Komornik in \cite{komornik1987controlabilite}.
$ \hfill \square $
\end{remark}

\section{The Main Controllability Result}
\label{sec:II.4}
The main boundary controllability result for the linear wave equation reads as follows:
\begin{theorem}
\label{theorem:II.2.3}
Let $ x^0 $ be any point of $ \mathbb{R}^n $ and $ T > 2R(x^0) $. Then, for each initial data $ \{ y^0, y^1 \} \in L^2(\Omega) \times H^{-1}(\Omega) $, there is a control $ v \in L^2(\Sigma(x^0)) $ such that the solution of
\begin{equation}\label{eqII33}
    \begin{cases}
    y^{\prime\prime} - \Delta y = 0 \hspace{1.6cm} \text{in} \hspace{0.3cm} Q
    \\
    y = \begin{cases}
    v & \hspace{2cm} \text{on} \hspace{0.3cm} \Sigma(x^0)
    \\ 0 & \hspace{2cm} \text{on} \hspace{0.3cm} \Sigma_\star (x^0)
    \end{cases}
    \\
    y(0) = y^0, y^\prime (0) = y^1 \hspace{0.35cm} \text{in} \hspace{0.3cm} \Omega
    \end{cases}
\end{equation}
satisfies $ y(T) = y^\prime (T) = 0 $.
\end{theorem}

\begin{proof}
\label{proof:II.theorem.2.3}
The conclusion follows as an immediate consequence of the application of the HUM method developed in Section \ref{sec:II.1} and Theorems \ref{theorem:II.2.1} and \ref{theorem:II.2.2}.

Indeed, by Theorem \ref{theorem:II.2.1} with $f=0$ and $\theta = \phi$, and by Theorem \ref{theorem:II.2.2}, it follows that, if $ T>2R(x^0)$, (\ref{eqII16}) holds. This implies that $F = H_0^1(\Omega) \times L^2(\Omega)$ and therefore $ F^\prime = H^{-1}(\Omega) \times L^2(\Omega)$.
\end{proof}

\begin{remark}
\label{remark:II.2.8}
As shown in the next section, the solution $ y = y(x,t) $ of (\ref{eqII33}) in the conditions of Theorem \ref{theorem:II.2.2} belongs to the class
$
C( [0, T]; L^2(\Omega) ) \cap C^1 ([0, T]; H^{-1}(\Omega)).    
$
Therefore, the traces of $y$ and $y^\prime$ for $ t = T$ are well defined.
$ \hfill \square $
\end{remark}

\section{Weak Solutions of the Wave Equation with Inhomogeneous Boundary Conditions}
\label{sec:II.5}
In this section the existence, uniqueness and regularity of the solutions of
\begin{equation}\label{eqII34}
    \begin{cases}
    y^{\prime\prime} - \Delta y = 0 & \text{in} \hspace{0.3cm} Q
    \\ y = v & \text{on} \hspace{0.3cm} \Sigma
    \\ y(0) = y^0, y^\prime (0) = y^1 & \text{in} \hspace{0.3cm} \Omega.
    \end{cases}
\end{equation}
are studied. In particular, the following result is obtained:

\begin{theorem}
\label{theorem:II.2.4}
For each initial data $ \{ y^0, y^1 \} \in L^2(\Omega) \times H^{-1}(\Omega) $ and each boundary condition $ v \in L^2(\Sigma) $ there is a unique solution of (\ref{eqII34}) in the class
\begin{equation}\label{eqII35}
y \in C ([0,T]; L^2(\Omega)) \cap C^1 ([0,T]; H^{-1}(\Omega)).
\end{equation}

In addition, the map $ \{ y^0, y^1, v \} \longrightarrow  \{ y, y^\prime \} $ is linear and continuous from $L^2(\Omega) \times H^{-1}(\Omega) \times L^2(\Sigma) $ into $C([0,T]; L^2(\Omega)) \times C([0,T]; H^{-1}(\Omega))$, that is, for each $T>0$, there is a constant $C(T)>0$ such that
\begin{equation}\label{eqII36}
\| (y, y^\prime) \|_{L^\infty(0, T; L^2(\Omega) \times H^{-1}(\Omega))}  \leq C(T) \{ \| y^0 \|_{L^2(\Omega)} +  \| y^1 \|_{H^{-1}(\Omega)} +  \| v \|_{L^2(\Sigma)} \}.
\end{equation}
\end{theorem}

\begin{proof}
\label{proof:II.theorem.2.4}
In the case of homogeneous boundary conditions, i.e. $ v = 0 $, the result is well known. Due to the linearity of the system, it is enough to consider the case for which $ y^0 = y^1 = 0, v \neq 0 $.
The solutions of
\begin{equation}\label{eqII37}
    \begin{cases}
    y^{\prime\prime} - \Delta y = 0 & \text{in} \hspace{0.3cm} Q
    \\ y = v & \text{in} \hspace{0.3cm} \Sigma
    \\ y(0) = y^\prime (0) = 0 & \text{in} \hspace{0.3cm} \Omega
    \end{cases}
\end{equation}
are defined by the method of transposition as follows.
Let $ \theta = \theta (x,t) $ be a solution of the system
\begin{equation}\label{eqII38}
    \begin{cases}
    \theta^{\prime\prime} - \Delta \theta = f & \text{in} \hspace{0.3cm} Q
    \\ \theta = 0 & \text{on} \hspace{0.3cm} \Sigma
    \\ \theta(T) = \theta^\prime (T) = 0 & \text{in} \hspace{0.3cm} \Omega
    \end{cases}
\end{equation}
with $ f \in L^1(0, T; L^2(\Omega)) $.

Multiplying (formally) equation (\ref{eqII37}) by $\theta$ and integrating by parts, lead to
\begin{equation}\label{eqII39}
\int_Q yf = \int_\Sigma v \frac{\partial\theta}{\partial \nu} d\Sigma.
\end{equation}
We  adopt (\ref{eqII39}) as the formulation of the weak solution of (\ref{eqII37}). That is, $ y = y(x,t) $ is a solution of (\ref{eqII37}) if and only if it satisfies (\ref{eqII39}) for all $ f \in L^1(0, T; L^2(\Omega))$. The estimate (\ref{eqII19}) holds for the solutions of (\ref{eqII38}), due to the reversibility of the wave equation with respect to the time variable. Thus,
\begin{align}\label{eqII40}
\left| \int_\Sigma v \frac{\partial\theta}{\partial \nu} d\Sigma\right| & \leq \| v \|_{L^2(\Sigma)} \left\| \frac{\partial\theta}{\partial \nu} \right\|_{L^2(\Sigma)} \leq C  \| v \|_{L^2(\Sigma)}  \| f \|_{L^1(0, T; L^2(\Omega))} .
\end{align}

Next, it follows from (\ref{eqII40}) that the map
$
    f \to \int_\Sigma v \frac{d\theta}{d\nu} d\Sigma
$
is linear and continuous in $ L^1(0, T; L^2(\Omega)) $. Thus, the existence of a unique solution $ y \in L^\infty(0,T;L^2(\Omega)) $ satisfying (\ref{eqII39}) is deduced, for all $ f \in L^1(0, T; L^2(\Omega)) $. The inequality (\ref{eqII40}) also provides the following estimate
\begin{align}\label{eqII41}
\begin{split}
\| y \|_{L^\infty(0, T; L^2(\Omega))} \leq C  \| v \|_{L^2(\Sigma)}.
\end{split}
\end{align}

To check the time-continuity of the solution $ t \in [0,T] \rightarrow y(t) \in L^2(\Omega) $, first we approximate $ v \in L^2(\Sigma) $ by a sequence $ \{ v_n \} \subset \mathcal{D}((0,T); C^2(\Gamma)) $ (the space of $C^\infty$ and compactly supported functions $(0,T) \rightarrow C^2(\Gamma)$) such that
\begin{equation}\label{eqII42}
    v_n \rightarrow v \hspace{0.3cm} \text{in} \hspace{0.3cm} L^2(\Sigma).
\end{equation}

Let $\{y_n\}$ be the solution of (\ref{eqII37}) associated with $\{v_n\}$. As $v_n$ is regular, $y_n$ is regular as well and, in particular, 
$
 y_n \in C([0,T]; L^2(\Omega)).
$
On the other hand, it follows from (\ref{eqII41}) and (\ref{eqII42}) that
$
    y_n \rightarrow y \in L^\infty(0,T; L^2(\Omega)),
$
and therefore $ y \in C([0,T]; L^2(\Omega)) $.

Now we check that $ y^\prime \in C([0,T]; H^{-1}(\Omega)) $. Given $ y \in \mathcal{D}(Q) $ (the space of $C^\infty$ and compactly supported functions in $Q$) it follows from (\ref{eqII39}), with $ f = g^\prime $, that
\begin{equation}\label{eqII43}
    \int_Q y g^\prime = \int_\Sigma v \frac{\partial \theta}{\partial \nu}
\end{equation}
where $ \theta = \theta (x, t) $ is solution of the system
\begin{eqnarray}\label{eqII44}
\begin{cases}
\theta^{\prime\prime} - \Delta \theta = g^\prime & \text{in} \hspace{0.3cm} Q
\\ \theta = 0 & \text{in} \hspace{0.3cm} \Sigma
\\ \theta(T) = \theta^\prime(T) = 0 & \text{in} \hspace{0.3cm} \Omega.
\end{cases}
\end{eqnarray}
Next, suppose that for all $ T>0 $, there exists $ C(T) > 0 $ such that
\begin{equation}\label{eqII45}
    \left\| \frac{\partial \theta}{\partial \nu} \right\|_{L^2(\Sigma)} \leq C(T) \left\| g \right\|_{L^1(0, T; H_{0}^1(\Omega))}
\end{equation}
for every solution of (\ref{eqII44}).
Under these conditions, again by duality, it follows that
\begin{equation}\label{eqII46}
y^\prime \in L^\infty(0, T; H^{-1}(\Omega))
\end{equation}
and that the following estimate holds
\begin{equation}\label{eqII47}
    \| y^\prime \|_{L^\infty(0, T; H^{-1}(\Omega))}
    \leq C \left\| v \right\|_{L^2(\Sigma)}.
\end{equation}

Finally, by  density,
\begin{equation}\label{eqII48}
y^\prime \in C([0, T]; H^{-1}(\Omega)).
\end{equation}

It is therefore sufficient to prove  estimate (\ref{eqII45}). This is the goal of the next result.
\end{proof}

\begin{lemma}
\label{lemma:II.2.2}
For all $ T>0 $, there is a constant $ C = C(T) > 0 $ such that
\begin{equation}\label{eqII49}
    \left\| \frac{\partial \theta}{\partial \nu} \right\|_{L^2(\Sigma)} \leq C(T) \left\| g \right\|_{L^1(0, T; H_{0}^1(\Omega))}
\end{equation}
for every solution of (\ref{eqII44}).
\end{lemma}

\begin{proof}
\label{proof:II.lemma.2.2}
By density, (\ref{eqII49}) needs to be proved for all $ g \in \mathcal{D}(Q) $.

As in the proof of Theorem \ref{theorem:II.2.1}, apply the identity (\ref{eqII18}) with $ q=h $ where $ h \in (C^1(\overline{\Omega}))^n $ satisfies
$
    h = \nu  $ over $ \Gamma.
$
Then,
\begin{align}\label{eqII50}
    \begin{split}
    \frac{1}{2} \int_\Sigma \left| \frac{\partial \theta}{\partial \nu} \right|^2 d\Sigma = {-\frac{1}{2}} \int_Q (\operatorname{div} h) | \nabla \theta |^2
    + \int_Q \frac{\partial h_j}{\partial x_k} \hspace{0.1cm} \frac{\partial \theta}{\partial x_k} \hspace{0.1cm} \frac{\partial \theta}{\partial x_j} + 
    \int_Q ( \theta^{\prime \prime} - g^\prime )h \cdot \nabla \theta
    \end{split}
\end{align}
where we used the following identity
\begin{equation*}
\int_Q \theta^{\prime \prime} h \cdot \nabla \theta =  \frac{1}{2} \int_Q (\operatorname{div} h) | \theta^\prime |^2 + \int_\Omega \theta^\prime h \cdot \nabla\theta \Big|_0^T .
\end{equation*}

 The proof of (\ref{eqII49}) is reduced to prove the following inequality
\begin{align}\label{eqII54}
    \begin{split}
    \left| {-\frac{1}{2}} \int_Q (\operatorname{div} h) | \nabla \theta |^2 + 
    \int_Q \frac{\partial h_j}{\partial x_k} \hspace{0.1cm} \frac{\partial \theta}{\partial x_k} \hspace{0.1cm} \frac{\partial \theta}{\partial x_j} \right| \leq C \left\| g \right\|^2_{L^1(0, T; H_0^1(\Omega))}
    \end{split}
\end{align}
and  \begin{equation}\label{eqII55}
    \left| \int_Q (\theta^{\prime\prime} - g^\prime) h \cdot \nabla \theta \right|
    \leq C \left\| g \right\|^2_{L^1(0, T; H_0^1(\Omega))}.
\end{equation}

Now, let us observe that $ \theta = \psi^\prime $, where $ \psi = \psi (x, t) $ solves
\begin{eqnarray}\label{eqII51}
\begin{cases}
\psi^{\prime\prime} - \Delta \psi = g & \text{in} \hspace{0.3cm} Q
\\ \psi = 0 & \text{on} \hspace{0.3cm} \Sigma
\\ \psi(T) = \psi^\prime(T) = 0 & \text{in} \hspace{0.3cm} \Omega.
\end{cases}
\end{eqnarray}
Classical regularity results for the wave equation guarantee 
\begin{align}\label{eqII52}
\begin{split}
    \| \psi \|_{L^\infty(0, T; H^2(\Omega))}
   & + \left\| \psi^\prime \right\|_{L^\infty(0, T; H_0^1(\Omega))} \leq C \| g \|_{L^1(0, T; H_0^1(\Omega))}
\end{split}
\end{align}
leading to
\begin{equation}\label{eqII53}
    \left\| \theta \right\|_{L^\infty(0, T; H_0^1(\Omega))} \leq 
    C \left\| g \right\|_{L^1(0, T; H_0^1(\Omega))}.
\end{equation}

This suffices to show \eqref{eqII54}.

To show that \eqref{eqII55} holds as well we observe that
\begin{align*}
\begin{split}
\int_Q (\theta^{\prime\prime} - g^\prime) h \cdot \nabla \theta &= \int_Q (\psi^{\prime\prime\prime} - g^\prime) h \cdot \nabla \psi^\prime 
\\&= {-} \int_Q (\psi^{\prime\prime} - g) h \cdot \nabla \psi^{\prime\prime} + \int_Q (\psi^{\prime\prime} - g) h \cdot \nabla \psi^\prime \bigg|_{0}^T.
\end{split}
\end{align*}

Since $ \psi^{\prime\prime} - g = \Delta \psi $,
\begin{equation}\label{eqII56}
\int_Q (\theta^{\prime\prime} - g^\prime) h \cdot \nabla \theta = \int_\Omega \Delta \psi h \cdot \nabla \psi^\prime \bigg|_{0}^T - \int_Q (\psi^{\prime\prime} - g) h \cdot \nabla \psi^{\prime\prime}.
\end{equation}

It follows from (\ref{eqII52}) that
\begin{equation*}
\left| \int_\Omega \Delta \psi h \cdot \nabla \psi^\prime \bigg |^T_0 \right|  \leq C \left\| g \right\|^2_{L^1(0, T; H_0^1(\Omega))}.
\end{equation*}

Therefore, it is now sufficient to estimate the last term of the identity (\ref{eqII56}):
\begin{align}\label{eqII57}
\begin{split}
\int_Q (\psi^{\prime\prime} - g) h \cdot \nabla \psi^{\prime\prime} =& {-}\frac{1}{2} \int_Q (\operatorname{div} h) |\psi^{\prime\prime}|^2 + \int_Q \operatorname{div} (g h) \psi^{\prime\prime}
\\ =& {-}\frac{1}{2} \int_Q (\operatorname{div} h) ( |\Delta \psi|^2 + |g|^2 + 2\Delta \psi g) + \int_Q \operatorname{div}(g h) (\Delta \psi + g)
\\ =& {-}\frac{1}{2} \int_Q (\operatorname{div} h) ( |\Delta \psi|^2 + 2\Delta \psi g) + \int_Q \operatorname{div}(g h )\Delta \psi
\\ &  -\frac{1}{2} \int_Q (\operatorname{div} h) |g|^2 + \int_Q g \operatorname{div}(g h).
\end{split}
\end{align}

However,
\begin{align*}
\begin{split}
\int_Q g \operatorname{div}(g h) = &\int_Q |g|^2 \operatorname{div} h + \int_Q g \nabla g \cdot h 
\\ =& \int_Q (\operatorname{div} h) |g|^2 +  \frac{1}{2} \int_Q h \cdot \nabla (g^2) = \frac{1}{2} \int_Q (\operatorname{div} h) |g|^2,
\end{split}
\end{align*}
and therefore,
\begin{equation}\label{eqII58}
\int_Q (\psi^{\prime\prime} -g)h \cdot \nabla\psi^{\prime\prime} = {-}\frac{1}{2} \int_Q (\operatorname{div} h) (|\Delta \psi|^2 + 2\Delta \psi g) + \int_Q \operatorname{div}(g h )\Delta \psi.
\end{equation}
Finally, combining (\ref{eqII52}) and (\ref{eqII58}), it follows easily that
\begin{equation*}
\left| \int_Q (\psi^{\prime\prime} -g)h \cdot \nabla\psi^{\prime\prime} \right| \leq C \left\| g \right\|^2_{L^1(0, T; H_0^1(\Omega))},
\end{equation*}
which concludes the proof of the main result.
\end{proof}
Several remarks are in order.
\label{remark:II.2.9}
\begin{remark}
Theorem \ref{theorem:II.2.4} is due to J.-L. Lions \cite{lions1988controlabilite}.

An in-depth presentation of the transposition method and many of its applications to solve Partial Differential Equations under very weak regularity assumptions can be found in J.-L. Lions and E. Magenes \cite{lions1968problemes}.
$ \hfill \square $
\end{remark}

\label{remark:II.2.10}
\begin{remark}
All of the shown results can be extended to the case where $\Omega$ is an open and convex set, without additional assumptions on the regularity of $\Gamma$ (cf. J.-L. Lions \cite{lions1988controlabilite}, Ch. I). This is simply due to the regularity properties ($H^2$ namely) of the solutions of the Laplace equation with homogeneous Dirichlet boundary conditions in such a geometric setting.
$ \hfill \square $
\end{remark}

\section{Some Geometric Considerations}
\label{sec:II.6}

\begin{enumerate}
\item The exact controllability of the wave equation in $ L^2(\Omega) \times H^{-1}(\Omega) $ with controls supported on subsets of the boundary of the form $ \Gamma(x^0) $ has been proved.
It is crucial to note that, in general, $ \Gamma(x^0) $ constitutes a "large" subset of the boundary.

Indeed, in the case where $ \Omega $ is square, i.e. $ \Omega = (0,1) \times (0,1) \subset \mathbb{R}^2 $, for instance, several cases for $x^0$ and $\Gamma(x^0)$ are to be distinguished:
\begin{enumerate}
\item $ \text{If } x^0 \in \Omega, \Gamma(x^0) = \Gamma $
\item $ \text{If } x^0 \in \mathbb{R}^2 \backslash \Omega, \Gamma(x^0) \subset \Gamma = \text{ the union of two or three adjacent sides} $
\end{enumerate}
Therefore, $ \Gamma(x^0) $ contains at least half of the boundary $ \Gamma $.

By fixing the sides $ \Gamma_0 $ on which the control acts and by choosing $x^0$ appropriately, the time of exact controllability $ 2R(x^0) $ can be minimized. If $ \Gamma_0 = \Gamma $, choosing $ x^0 = (1/2, 1/2) $ at the center leads to $ 2R(x^0) $ = diameter of $ \Omega = \sqrt{2}$. On the other hand, if $ \Gamma_0 $ is the union of two consecutive sides, choosing $ x^0 $ on the opposite vertex leads to $ 2R(x^0) = 2 \sqrt{2} $. If $ \Gamma_0 $ is the union of three sides, choosing $ x^0 $ as the midpoint of the fourth side leads to $ 2R(x^0) = \sqrt{5} $.

On the other hand, if $ \Omega = B(0,1) \subset \mathbb{R}^2 $ is the unit disk, it is easy to check that $ \Gamma (x^0) $ contains a half-circle for all $ x^0 \in \mathbb{R}^2 $. Additionally, $ \Gamma (x^0) $ approaches a half-circle if and only if $ |x^0| \rightarrow \infty $, that is, when the controllability time diverges: $ 2R(x^0) \rightarrow \infty $. In the case of controls supported in all $ \Gamma $, the controllability time is minimized by taking $ x^0 = (0,0) $ at the center of the circle, so that $ 2R(x^0) = 2 $.

\item C. Bardos, G. Lebeau and J. Rauch in \cite{bardos1992sharp} considerably improved these results applying microlocal analysis techniques. Their results apply to domains $C^\infty$ and $ \Omega $ and it is roughly the following: {\it "The wave equation is exactly controllable in $ L^2(\Omega) \times H^{-1}(\Omega) $ with controls in $ L^2(\Gamma_0 \times (0,T)) $ if and only if every ray propagating in $ \Omega $ and reflected in its boundary according to the laws of geometric optics, intersects $ \Gamma_0 $ at a non-diffractive point and at time $ t < T $".}

This result guarantees essentially the equivalence between the controllability in a sharp functional setting and the so-called Geometric Control Condition (GCC), which is of a microlocal nature.

This result is optimal since it characterizes the sets $ \Gamma_0 $ for which exact controllability in $ L^2(\Omega) \times H^{-1}(\Omega) $ holds with controls in $ L^2(\Gamma_0 \times (0,T)) $.

By means of this result, Theorem \ref{theorem:II.2.3} is significantly improved, allowing to prove the controllability for a more general class of supports $ \Gamma_0 \subset \Gamma $ of the control and reducing the time needed according the multiplier method. For example, in the case where $ \Omega $ is a disk, exact controllability holds with controls supported in any subset of $ \Gamma $ containing a semicircle. But the time of control that microlocal analysis techniques yield is in general smaller than the one obtained by the method of multipliers. Microlocal tools allow also showing that waves can be controlled acting on three disjoint parts of the boundary that do not contain any half-circle, a result that cannot be achieved by the multiplier technique, as presented above.

\item As mentioned in Section \ref{sec:II.1}, Holmgren's Theorem ensures that given $ \Gamma_0 \subset \Gamma $ any open non-empty subset, there exists $ T_0 = T_0 (\Gamma_0, \Omega) $ such that for all $ T>T_0 $ the uniqueness result is verified (\ref{eqII09}). Applying the HUM method, we deduce that if $ T>T_0 $, the wave equation is exactly controllable with controls in $ L^2(\Gamma_0 \times (0,T)) $ and for initial data $ \{ y^0, y^1 \}$ such that $ \{ y^1, {-}y^0 \}  \in F^\prime$ where $ F = F(\Gamma_0, \Omega, T) $ is a Hilbert space that a priori depends on $ \Gamma_0, \Omega \text{ and } T$.

In Theorem \ref{theorem:II.2.3} we have shown that if $ \Gamma_0 = \Gamma(x^0) $ and $ T>2R(x^0) $, then $ F = H_{0}^1(\Omega) \times L^2(\Omega) $.
On the other hand, the microlocal tools and results in \cite{bardos1992sharp} provide a characterization of the sets $ \Gamma_0 $ and the values of $T$ for which $ F = H_{0}^1(\Omega) \times L^2(\Omega) $.

However, the characterization (or obtaining useful and manageable information) of the space $F$ when these geometric conditions are not fulfilled, i. e. when $F$ is strictly contained in $ H_{0}^1(\Omega) \times L^2(\Omega) $, is to a large extent an open problem. A. Haraux \cite{haraux1988controlabilite} contributed in this direction when the support of the control is a subset of $ \Omega $.  On the other hand,  \cite{allibert1999analytic} provides a sharp description of the nature of the space of controllable data in some particular geometries, like surfaces of revolution. But the issue is still open for general geometries. $ \hfill \square $
\end{enumerate}

\section{The Case of Variable Coefficients in One Space Dimension}
\label{sec:II.7}
We consider the following variable coefficients wave equation
\begin{eqnarray}\label{eqII59}
\begin{cases}
    y^{\prime\prime} - \operatorname{div}(a(x) \nabla y) = 0 \hspace{0.3cm} \text{ in } \hspace{0.3cm} Q
    \\y = \begin{cases}
        v \hspace{2.15cm} & \text{ on } \hspace{0.2cm} \Sigma_0
        \\0 \hspace{2.15cm} & \text{ on } \hspace{0.2cm} \Sigma \backslash \Sigma_0
    \end{cases}
    \\y(0) = y^0, y^\prime (0) = y^1
\end{cases}
\end{eqnarray}
with
\begin{equation}\label{eqII60}
    a = a(x) \in C^1(\overline{\Omega})
\end{equation}
such that
\begin{equation}\label{eqII61}
  \exists  a_0 > 0 : a(x) \geq a_0 \hspace{0.3cm} \forall x \in \Omega.
\end{equation}

The exact controllability of this equation is a problem that is not yet completely well understood. The following partial results are available:
\begin{itemize}
\item If $ \Omega $ is $ C^\infty $ smooth and $ a \in C^\infty (\overline{\Omega}) $, the results in \cite{bardos1992sharp} allow characterizing the subsets $ \Gamma_0 $ of $ \Gamma $ and the values of $ \Gamma $ for which (\ref{eqII59}) is exactly controllable in $ L^2(\Omega) \times H^{-1}(\Omega) $, with controls in $ L^2(\Gamma_0 \times (0, T)) $. 

Note however that bicharacteristic rays are defined according to the symbol of the variable coefficients wave operator. Therefore, the GCC  depends on these coefficients and, in some cases,  there may even exist bicharactetistic rays that remain permanently in the interior of $\Omega$, without never reaching the exterior boundary $\partial\Omega$. In these cases the boundary observability may not hold even if the control acts everywhere on the boundary, and this independently of the length of the time-horizon $T$. We refer to the following comment for further details on this matter.

\item  V. Komornik in \cite{komornik1989exact} observed that the multiplier method used in the previous sections applies, for $ T>0 $ large enough, if $ a= a(x) $ verifies the following additional hypothesis
\begin{equation}\label{eqII62}
    \exists \delta > 0 : (1 - \delta) a(x) - \frac{\operatorname{div} \hspace{0.1cm} a(x) \cdot (x-x^0)}{2} \geq 0  \hspace{0.3cm} \forall x \in \Omega.
\end{equation}

In \cite{macia2002zuazua}  it was shown that  (\ref{eqII62}) is sharp, since when
$
    a = \frac{1}{2} \operatorname{div} a \cdot (x - x_0)
$
there may exist rays of geometric optics trapped inside $\Omega$, that never reach the exterior boundary, making the controllability property impossible.

However, the general case, without structural conditions on the coefficients $a=a(x)$,  cannot be treated by multiplier methods.
The purpose of this section is to prove the exact controllability in the $1-d$ case, without any restriction on the coefficient other than (\ref{eqII60}) and (\ref{eqII61}). 
\end{itemize}
The $1-d$ result we present here shows that, according to our intuition, the characteristic rays cross the whole domain bouncing back and forth on the boundary. This makes  the concentration of waves in the interior  impossible and ensures the boundary controllability in a sufficiently large characteristic time. Note however that, for this to hold, some regularity assumptions on the coefficient need to be imposed, roughly $a \in BV$.

Let therefore $ \Omega = (0,1) \subset \mathbb{R} $ and consider the system $1-d$
\begin{eqnarray}\label{eqII63}
\begin{cases}
    y^{\prime\prime} - (a(x) y_x)_x = 0 \hspace{0.3cm} \text{ in } \hspace{0.3cm} Q = (0,1) \times (0,T)
    \\y(0, t) = v(t) \hspace{1.15cm} \text{ on } \hspace{0.3cm} (0,T)
    \\y(1, t) = 0 \hspace{1.56cm} \text{ on } \hspace{0.3cm} (0,T)
    \\y(0) = y^0, y^\prime (0) = y^1 \hspace{0.1cm} \text{ in } \hspace{0.25cm} (0,1) .
\end{cases}
\end{eqnarray}

The following result is obtained:
\begin{theorem}
\label{theorem:II.2.5}
Suppose $ a \in BV([0,1]) $ is such that (\ref{eqII61}) holds. If $ \hspace{0.1cm} T > 2/\sqrt{a_0} \hspace{0.2cm} $ for all initial data $ \{y^0, y^1 \} \in L^2(0,1) \times H^{-1}(0,1) $, there exists a control $ v = v(t) \in L^2(0,T) $ such that the solution of (\ref{eqII63}) satisfies (\ref{eqII02}).
\end{theorem}

\begin{remark}
\label{remark:II.2.11}
The same result is obtained when the control acts at the tip $ x = 1$, that is, for boundary conditions of the form 
$
    y(0,t) = 0; y(1,t) = v(t)
$ for  $ t \in (0,T)$.

When the control acts simultaneously at both ends of the interval, that is,
\begin{equation*}
    y(0,t) = v(t); y(1,t) = w(t) \hspace{0.3cm} \forall t \in (0,T),
\end{equation*}
the exact controllability holds for $ T > 1/\sqrt{a_0} $  with controls $ v(t), w(t) \in L^2(0,T). $
\end{remark}

\begin{proof}
\label{proof:II.remark.2.11}
We present the proof when $a \in C^1[0,1]$ to avoid technical difficulties. The fact that the same  holds for $ a \in BV(0,1)$ was observed in \cite{fernandez2002zuazua} and can be derived by density.

The HUM method easily applies to equation (\ref{eqII63}). In this way the problem is reduced to obtain the following estimate
\begin{equation}\label{eqII64}
    \left\| \phi^0 \right\|_{H_0^1 (\Omega)}^2 + \left\| \phi^1 \right\|_{L^2 (\Omega)}^2 \leq C \int_{0}^T \left| \phi_x(0,t) \right|^2 dt
\end{equation}
where $\phi$ solves
\begin{eqnarray}\label{eqII65}
\begin{cases}
    \phi^{\prime\prime} - (a(x) \phi_x)_x = 0 \hspace{0.3cm} & \text{in } \hspace{0.2cm} (0,1) \times (0,T)
    \\ \phi(0, t) = \phi(1,t) = 0 \hspace{1.15cm} & \text{on } \hspace{0.2cm} t \in (0,T)
    \\ \phi(0) = \phi^0, \phi^\prime (0) = \phi^1  \hspace{1.15cm}& \text{in } \hspace{0.2cm} (0,1).
\end{cases}
\end{eqnarray}

As we noted in Remark \ref{remark:II.2.11}, the direct inequality is easily extended to systems with variable coefficients in any dimension $n\ge 1$. It holds therefore also in the present $1-d$ setting. Thus, we  focus on the observability inequality. Note however that the techniques we shall use for proving observability can also be employed for the hidden regularity, to get the direct inequality for $BV$ coefficients.

To show the inverse inequality, we define the sidewise energy functional
\begin{equation}\label{eqII66}
F(x) = \frac{1}{2} \int_{\alpha x}^{T-\alpha x} \{ | \phi^\prime (x,t) |^2 + a(x) | \phi_x (x,t) |^2 \} dt
\end{equation}
where $ \alpha = 1/\sqrt{a_0}.$
Note that
\begin{equation}\label{eqII67}
F(0) = \frac{a(0)}{2} \int_{0}^{T} | \phi_x (0,t) |^2 dt.
\end{equation}

Now, proceed computing the derivative of $F$ as follows:
\begin{align}\label{eqII68}
\begin{split}
F^\prime (x) = &\frac{dF(x)}{dx}
\\  = &\int_{\alpha x}^{T-\alpha x} \left[ \phi^\prime (x,t) \phi_x^\prime (x,t)  + \frac{a_x(x)}{2} | \phi_x (x,t) |^2 + a(x) \phi_x (x,t) \phi_{xx} (x,t) \right]dt
\\ & {{-}\frac{\alpha}{2}} \sum_{t=\alpha x, T-\alpha x} \left[ | \phi^\prime (x,t) |^2 + a(x) | \phi_x (x,t) |^2 \right].
\end{split}
\end{align}

Integration by parts on time leads to
\begin{eqnarray*}
\int_{\alpha x}^{T-\alpha x} \phi^\prime (x,t) \phi_x^\prime (x,t) dt = {-} \int_{\alpha x}^{T-\alpha x} \phi^{\prime\prime} (x,t) \phi_x (x,t) dt + \phi^\prime (x,t) \phi_x (x,t) \bigg|^{T-\alpha x}_{\alpha x} \end{eqnarray*}
and we deduce that
\begin{align}\label{eqII69}
\begin{split}
F^\prime (x) & = \int_{\alpha x}^{T-\alpha x} \left[ {-} \phi^{\prime\prime} (x,t) + a(x) \phi_{xx} (x,t) + a_x(x) \phi_x (x,t) \right] \phi_x (x,t) dt
\\ & \quad {{-}\frac{1}{2}} \int_{\alpha x}^{T-\alpha x} a_x(x) \left| \phi_x (x,t) \right|^2 dt + \phi^\prime (x,t) \phi_x (x,t) \bigg|^{T-\alpha x}_{\alpha x}
\\ & \quad {{-}\frac{\alpha}{2}} \sum_{t=\alpha x, T-\alpha x} \left[ | \phi^\prime (x,t) |^2 + a(x) | \phi_x (x,t) |^2 \right]
\\ & \leq {{-}\frac{1}{2}} \int_{\alpha x}^{T-\alpha x} a_x(x) \left| \phi_x (x,t) \right|^2 dt
\end{split}
\end{align}
 since
\begin{equation*}
    \phi^{\prime \prime} - a(x) \phi_{xx} - a_x \phi_x = \phi^{\prime \prime} - (a(x) \phi_x)_x = 0
\end{equation*}
and
\begin{align*}
\left| \phi^\prime (x,t) \phi_x (x,t) \right| &\leq \frac{1}{\sqrt{a_0}} \left( \frac{1}{2} \left| \phi^\prime (x,t) \right|^2 + \frac{a_0}{2} \left| \phi_x (x,t) \right|^2 \right) 
\\ &\leq \frac{\alpha}{2} \left( \left| \phi^\prime (x,t) \right|^2 + a(x) \left| \phi_x (x,t) \right|^2 \right), \hspace{0.3cm} \forall t \in [ 0,T ], \forall x \in (0,1).
\end{align*}

By (\ref{eqII69})
\begin{equation} \label{eqII70}
    F^\prime(x) \leq \frac{\|a_x\|_{L^\infty(0,1)}}{2a_0} \int_{\alpha x}^{T-\alpha x} a(x) \left| \phi_x (x,t) \right|^2 dt \leq  \frac{\|a_x\|_{L^\infty(0,1)}}{a_0} F(x),
\end{equation}
and therefore,
\begin{equation} \label{eqII71}
    F(x) \leq e^{\|a_x\|_{L^\infty(0,1)} x/a_0}   F(0) \leq e^{\|a_x\|_{L^\infty(0,1)}/a_0} F(0), \hspace{0.3cm} \forall x \in (0,1).
\end{equation}

Next, integrating (\ref{eqII71}) with respect to $ x \in (0,1) $ we obtain
\begin{equation} \label{eqII72}
    \int_0^1 F(x) dx \leq e^{\|a_x\|_{L^\infty(0,1)}/a_0} F(0).
\end{equation}

Finally, as $ T> 2 \alpha = 2/\sqrt{a_0} $, observe that
\begin{align*}
 (T - 2 \alpha) E(0) &= \frac{1}{2} \int_{\alpha}^{T-\alpha} \int_0^1 \left[ \left| \phi^\prime (x,t) \right|^2 + a(x) \left| \phi_x (x,t) \right|^2 \right] dx dt 
 \\ &\leq \int_0^1 F(x) dx \leq e^{\|a_x\|_{L^\infty(0,1)}/a_0} F(0).
\end{align*}
This inequality, combined with the  conservation of energy property, 
\begin{equation*}
    E(t) = \frac{1}{2} \int_0^1 \{ \left| \phi^\prime (x,t) \right|^2 + a(x) \left| \phi_x (x,t) \right|^2 \} dx = E(0)
\end{equation*}
yields
\begin{align} \label{eqII73}
\begin{split}
\int_0^1 \left\{ |\phi^1(x)|^2 + a(x) |\phi_x^0(x)|^2 \right\} dx & 
\leq 
\frac{2e^{\|a_x\|_{L^\infty(0,1)}/a_0}} {(T - 2 \alpha)}F(0)
\\ &= \frac{a(0)e^{\|a_x\|_{L^\infty(0,1)}/a_0}} {(T - 2 \alpha)} \int_0^T  |\phi_x(0,t)|^2 dt. \hspace{0.3cm} \hfill \square
\end{split}
\end{align}
\end{proof}
Several remarks are in order.


\begin{remark}
\label{remark:II.2.12}

\begin{itemize}
\item Extending the proof above to $BV$ functions requires a more careful analysis of inequality \eqref{eqII69}. It suffices actually to apply Gronwall's inequality to get a similar estimate, replacing the norm $ \|a_x\|_{L^ \infty (0,1)} $ by $ \|a_x\|_{L^1 (0,1)} $. This yields the estimate for $a \in W^{1,1}$, that, by density, can be extended to any $a \in BV$ (see \cite{fernandez2002zuazua}). 

\item From (\ref{eqII69}) it follows that, when $ a(x) $ is increasing, $F$ is decreasing
and
the constant $C$ in (\ref{eqII64}) does not depend on $ \|a_x\|_{L^ \infty (0,1)} $. 

\item The monotonicity of the coefficient $a=a(x)$ fails in the context of homogenisation with rapidly oscillatory coefficients, a topic treated in  \cite{castro2000zuazua, castro2002zuazua}.
For instance, if $ a(x) \in C^1 (\mathbb{R}) \cap L^\infty (\mathbb{R}) $ is a non constant periodic function (therefore, non-monotonic) that satisfies $ a(x) \geq a_0 > 0 $ for all $ x \in \mathbb{R} $ the observability constant $ C_\varepsilon $ of the estimate (\ref{eqII64}) corresponding to the equation
\begin{equation} \label{eqII74}
    \phi^{\prime\prime} - \left( a\left( \frac{x}{\varepsilon} \right) \phi_x \right)_x = 0
\end{equation}
is of the order of $ C_\varepsilon \sim e^{c/\varepsilon} $ as $ \varepsilon \rightarrow 0. $

In \cite{castro2000zuazua, castro2002zuazua} sharp results were derived, showing that uniform observability results hold in an optimal class of low-frequency solutions. This allows recovering, in the limit as $\varepsilon \to 0$, the observability inequality for the homogenised constant coefficients wave equation.
\end{itemize}
$ \hfill \square $
\end{remark}

\section{Comments}
\label{sec:II.8}
\begin{enumerate}    
\item {\bf Optimal minimal-norm controls.} It follows from Theorem \ref{theorem:II.2.3} that, if $ T>2R(x^0) $, for each $ \{ y^0, y^1 \} \in L^2(\Omega) \times H^{-1}(\Omega) $, the set of admissible controls of (\ref{eqII74}), which is denoted by 
\begin{equation*}
    U_{ad} = \{ v \in L^2(\Sigma (x^0) ) : y \text{ is solution of (\ref{eqII33}) satisfiying (\ref{eqII02})}\}
\end{equation*}
contains infinitely many controls.

Indeed, let be $ \varepsilon > 0 $ such that $ T - \varepsilon > 2R(x^0) $. Given $ w \in L^2(\Gamma (x^0) \times (0, \varepsilon)) $, we solve the system

\begin{equation}\label{eqII75}
    \begin{cases}
    y^{\prime\prime} - \Delta y = 0 \hspace{1.7cm} \text{in} \hspace{0.3cm} \Omega \times (0, \varepsilon)
    \\
    y = \begin{cases}
    w & \hspace{2cm} \text{on} \hspace{0.3cm} \Gamma(x^0) \times (0, \varepsilon)
    \\ 0 & \hspace{2cm} \text{on} \hspace{0.3cm} \Gamma_\star (x^0) \times (0, \varepsilon)
    \end{cases}
    \\
    y(0) = y^0, y^\prime (0) = y^1  \hspace{0.5cm} \text{in} \hspace{0.3cm} \Omega.
    \end{cases}
\end{equation}

It follows from Theorem \ref{theorem:II.2.4} that
\begin{equation} \label{eqII76}
    y^0_\varepsilon = y(\varepsilon) \in L^2(\Omega), y_\varepsilon^1 = y^\prime (\varepsilon) \in H^{-1}(\Omega).
\end{equation}

Since $ T - \varepsilon >2R(x^0) $, Theorem \ref{theorem:II.2.3} ensures the existence of $ u \in L^2( \Gamma(x^0) \times (\varepsilon, T) ) $, such that the solution of the system
\begin{equation}\label{eqII77}
    \begin{cases}
    y^{\prime\prime} - \Delta y = 0 \hspace{1.6cm} \text{in} \hspace{0.3cm} \Omega \times (\varepsilon, T)
    \\
    y = \begin{cases}
    u & \hspace{2cm} \text{on} \hspace{0.3cm} \Gamma(x^0) \times (\varepsilon, T)
    \\ 0 & \hspace{2cm} \text{on} \hspace{0.3cm} \Gamma_\star (x^0) \times (\varepsilon, T)
    \end{cases}
    \\
    y(\varepsilon) = y_\varepsilon^0, y^\prime (\varepsilon) = y_\varepsilon^1 \hspace{0.4cm} \text{in} \hspace{0.3cm} \Omega
    \end{cases}
\end{equation}
satisfies (\ref{eqII02}).

Therefore, the control
\begin{equation}\label{eqII78}
    v = \begin{cases}
    w & \hspace{1cm} \text{in} \hspace{0.3cm} \Gamma(x^0) \times (0, \varepsilon)
    \\ u & \hspace{1cm} \text{in} \hspace{0.3cm} \Gamma (x^0) \times (\varepsilon, T)
    \end{cases}
\end{equation}
belongs to the set of admissible controls $ U_{ad} $, where $ w $ is an arbitrary element of the infinite dimensional space $ L^2(\Gamma(x^0) \times (0, \varepsilon)) $.

Having observed that $ U_{ad} $ contains  infinitely many elements, the following question  arises naturally: { \it Is there an optimal control $ v \in U_{ad} $ of minimal norm so that} 
\begin{equation}\label{eqII79}
    \| v \|^2_{L^2(\Sigma(x^0))} = \min_{u \in U_{ad}}\hspace{0.1cm} \| u \|^2_{L^2(\Sigma(x^0))} \text{?}
\end{equation}

The minimisation problem (\ref{eqII79}) admits a unique solution, which is actually the control obtained in Theorem \ref{theorem:II.2.3} that HUM yields.
In fact,  HUM consistently provides optimal controls of minimal norm in the appropriate functional setting in the various situations in which it is applied (c.f. J.-L. Lions \cite{lions1988controlabilite}, Chapter VIII).
$ \hfill \square $

\item {\bf Other boundary conditions.}  In this chapter we considered only Dirichlet-type boundary conditions. Both the HUM method and the multiplier techniques can be adapted to the Neumann boundary conditions or the Dirichlet-Neumann mixed ones (c.f. J.-L. Lions \cite{lions1988controlabilite}, Chapter III). However, for mixed boundary conditions, the singularities that the solutions present at the boundary points where the boundary condition type changes complicate the analysis considerably. In particular, in that setting, some of the integration by parts formulas being used in the application of multipliers are no longer true. P. Grisvard in \cite{grisvard1987controlabilite305}, \cite{grisvard1989controlabilite} proved that, although the standard identities are not fulfilled, the needed inequalities  for the desired a priori estimates actually hold.

On the other hand, the techniques and results of C. Bardos, G. Lebeau and J. Rauch  \cite{bardos1992sharp} showing that the sharp Geometric Control Condition suffices for controllability can be also adapted to the case of Neumann-type boundary conditions.
$ \hfill \square $

\item {\bf Less regular domains.} Throughout this chapter  the $ C^2 $ regularity of the domain $ \Omega $ (or its convexity) is assumed. These hypotheses can be considerably weakened (see  P. Grisvard \cite{grisvard1987controlabilite} and \cite{grisvard1989controlabilite}).
$ \hfill \square $

\item {\bf More regular controls.} The controllability property holds in $ L^2(\Omega) \times H^{-1}(\Omega) $ with $L^2$-controls. However, the HUM method and multiplier techniques can be adapted to work in other spaces. In particular, it can be shown that if the initial data belongs to $ H_0^1(\Omega) \times L^2(\Omega) $ then the HUM control lies in $ H_0^1(0, T; L^2(\Gamma(x^0)) $. This issue has been  analysed in \cite{ervedoza2010zuazua}, where the ellipticity of the HUM operator has been shown, in the sense that more regular data lead to more regular controls. 

Note however that the analysis in \cite{ervedoza2010zuazua} requires the use of a weight function modulating the controls in time,  making  the effects at extreme points $t=0, T$ vanish. This is essential to prove the ellipticity  of the HUM operator, since it allows to perform integration by parts on the time variable to gain regularity of the controls, without the artefacts generated by the extreme values at $t=0, T$.

Likewise, the exact controllability property in the dual larger space $ H^{-1}(\Omega) \times (H^2(\Omega) \cap H_0^1(\Omega)))^\prime $ can be proved by means of controls in the larger space $ H^{-1} (0,T; L^2(\Gamma(x^0)))  $ (cf. J.-L. Lions \cite{lions1988controlabilite}, Chapter I).
$ \hfill \square $

\item {\bf Model perturbations.}The problem of exact controllability under perturbations such as domain variations, rapidly oscillating coefficients, etc. has been studied by J.-L. Lions in \cite{lions1988controlabilite2}. This is a subject in which there are still plenty interesting open problems. We discuss here some of them.

\item  {\bf Semilinear equations.} The problem of the exact controllability of the semilinear wave equation has been studied in E. Zuazua \cite{zuazua1988controlabilite306}, \cite{zuazua1990exact}, \cite{zuazua1991exact} and solved, in particular,  for globally Lipschitz perturbations. The problem is open when the nonlinear term is superlinear at infinity (see E. Zuazua \cite{zuazua1990controlabilite} for a particular case in $1-d$). 

Note that, in any space dimension and for a broad class of nonlinearities, combining local controllability and stabilization results, one can achieve controllability in a time-horizon whose length depends on the size of the data to be controlled. But whether the control time can be made independent of the data to be controlled is an open problem (see \cite{dehman2003stabilization}).

\item {\bf Perforated domains.} The problem of exact controllability for the linear wave equation in perforated domains has been studied in D. Cioranescu, P. Donato and E. Zuazua \cite{cioranescu1989controlabilite} and \cite{cioranescu1992exact}.$ \hfill \square $

\item {\bf Plate and Schr\"odinger equations.} The HUM method and multiplier techniques can also be adapted to the study of the exact controllability of a vibrating plate equation of the type
\begin{equation}\label{eqII80}
    y^{\prime\prime} + \Delta^2 y = 0 \hspace{0.3cm} \text{in} \hspace{0.3cm} \Omega \times (0,T)
\end{equation}
with various boundary conditions (cf. J.-L. Lions \cite{lions1988controlabilite}, Chapter IV, J. Lagnese and J.-L. Lions \cite{lagnese1988modelling} and the bibliographies  therein).

The distinguished feature of this model is its intrinsic  infinite propagation speed, which, unlike the wave equation, allows exact controllability to hold in an arbitrarily small time.
Analogous results have been proved by E. Machtyngier \cite{machtyngier1990control} for the Schrödinger equation. In G. Lebeau \cite{lebeau1992controle}, it is shown, in a systematic manner, using diadic Fourier decompositions,  that the Geometric Control Condition on $ \Gamma_0 \subset \Gamma $ suffices for the exact controllability of the plate equation (\ref{eqII80}) and the Schrödinger equation, in the corresponding energy spaces, with controls in $ L^2(\Gamma_0 \times (0,T)) $, for $ T>0 $ arbitrarily small, as a consequence of the sharp controllability properties of the wave equation.

Likewise, exact controllability holds for the elasticity system (cf. J.-L. Lions \cite{lions1988controlabilite}, Chapter IV) and for Maxwell's equations (cf. J. Lagnese \cite{lagnese1989exact} and O. Nalin \cite{nalin1989control}).

$ \hfill \square $

\item {\bf Numerical approximation.}  R. Glowinski, C.H. Li and J.-L. Lions \cite{glowinski1990numerical} developed an effective numerical method for studying the exact controllability of the wave equation. We also refer to the complementary bibliography at the end of the book for further references on more recent developments at this respect.
$ \hfill \square $

\end{enumerate}


%
%

\chapter{Internal Controllability of the Linear Wave Equation}
\label{chapter03} 

\abstract{
In this chapter the results of the previous one are adapted for control subdomains in the interior of the domain where waves propagate, focusing on the particular case where the support of the control is placed near the boundary.
}

\section{Description of the HUM method}
\label{sec:III.1}
As mentioned in Chapter \ref{chapter01}, the problem of exact internal controllability of the linear wave equation can be formulated as follows.

Let $\omega$ be a nonempty open subset of $\Omega$ and consider the wave equation
\begin{equation}\label{eqIII01}
    \begin{cases}
    y^{\prime\prime} - \Delta y = h \chi_\omega \hspace{1.6cm} & \text{in} \hspace{0.3cm} Q
    \\
    y = 0 \hspace{0.5cm} & \text{in} \hspace{0.3cm} \Sigma
    \\
    y(0) = y^0, y^\prime (0) = y^1\hspace{0.5cm} & \text{in} \hspace{0.3cm} \Omega,
    \end{cases}
\end{equation}
where $ \chi_\omega $ denotes the characteristic function of $ \omega $.

We are interested in the controllability property in the energy  space $ H_0^1 (\Omega) \times L^2(\Omega) $ with controls in $ L^2( \omega \times (0,T) ) $.

The problem consists of finding $ T>0 $ such that for every pair of initial data  $ \{ y^0, y^1 \} \in H_0^1 (\Omega) \times L^2(\Omega)  $,  there exists a control $ h \in L^2(\omega \times (0,T)) $ such that the solution of (\ref{eqIII01}) satisfies
\begin{equation}\label{eqIII02}
y(T) = y^\prime (T) = 0.
\end{equation}

As in the context of boundary controllability problem, such  result requires that $ \omega $ satisfies some geometric properties.

The HUM method can be adapted to this problem as follows.

Given the $ \{ \phi^0, \phi^1 \} \in \mathcal{D}(\Omega) \times \mathcal{D}(\Omega) $, consider the system
\begin{equation}\label{eqIII03}
    \begin{cases}
    \phi^{\prime\prime} - \Delta \phi = 0 \hspace{1.6cm} & \text{in} \hspace{0.3cm} Q
    \\
    \phi = 0 \hspace{0.5cm} & \text{in} \hspace{0.3cm} \Sigma
    \\
    \phi(0) = \phi^0, \phi^\prime (0) = \phi^1\hspace{0.5cm} & \text{in} \hspace{0.3cm} \Omega
    \end{cases}
\end{equation}
and, then, solve
\begin{equation}\label{eqIII04}
    \begin{cases}
    y^{\prime\prime} - \Delta y = -\phi \chi_\omega \hspace{1.6cm} & \text{in} \hspace{0.3cm} Q
    \\
    y = 0 \hspace{0.5cm} & \text{in} \hspace{0.3cm} \Sigma
    \\
    y(T) = 0, y^\prime (T) = 0 \hspace{0.5cm} & \text{in} \hspace{0.3cm} \Omega
.
    \end{cases}
\end{equation}

Then, define the operator by
\begin{equation}\label{eqIII05}
  \Lambda \{ \phi^0, \phi^1 \} = \{ y^\prime(0), -y(0) \}.
\end{equation}

Multiplying in (\ref{eqIII04}) by $ \phi $ and integrating by parts, leads to
\begin{equation}\label{eqIII06}
  \langle  \Lambda \{ \phi^0, \phi^1 \}, \{ \phi^0, \phi^1 \} \rangle = 
  \int^T_0 \int_\omega |\phi|^2 dx dt.
\end{equation}

Suppose that the following uniqueness result holds
\begin{equation}\label{eqIII07}
    \phi = 0 \hspace{0.2cm} \text{in} \hspace{0.2cm} \omega \times (0,T) 
    \Rightarrow \phi^0 \equiv \phi^1 \equiv 0.
\end{equation}

Then,
\begin{equation}\label{eqIII08}
    \| \{ \phi^0, \phi^1 \} \|_F = \left( \int^T_0 \int_\omega  |\phi|^2 dx dt \right)^ \frac{1}{2}
\end{equation}
defines a norm in $ \mathcal{D}(\Omega) \times \mathcal{D}(\Omega) $, where\\
\begin{equation}\label{eqIII09}
   F = \text{ the Hilbert space completion } \mathcal{D}(\Omega) \times \mathcal{D}(\Omega) \text{ with respect to norm } \|\cdot\|_F.
\end{equation}

It follows from (\ref{eqIII06}) that $ \Lambda $ uniquely extends to an isomorphism from $ F \text{ to } F^\prime $.

Therefore, given $ \{ y^0, y^1 \} $ such that
$
    \{ y^1, -y^0 \} \in F^\prime,
$
there is a unique pair $ \{ \phi^0, \phi^1 \} \in F $ such that
\begin{equation}\label{eqIII10}
    \Lambda \{ \phi^0, \phi^1 \} = \{ y^1, -y^0 \}.
\end{equation}

In particular, if we apply the control $ h = - \phi \in L^2 ( \omega \times (0,T) )$ where $ \phi $ is the solution of (\ref{eqIII03}) with the initial data satisfying (\ref{eqIII10}), the solution of (\ref{eqIII04}) satisfies
$
y(0) = y^0, y^\prime (0) = y^1.
$

Holmgren's Theorem ensures that given any non-empty open subset $ \omega $ of $ \Omega $, there exists $ T_0 ( \omega, \Omega) $ such that if $ T > T_0 ( \omega, \Omega) $ the uniqueness result (\ref{eqIII07}) holds.

Therefore, we have shown that, if $ T > T_0 ( \omega, \Omega) $, the system (\ref{eqIII01}) is controllable in the space of the initial data such that $ \{ y^1, -y^0 \} \in F^\prime $ with the controls $ h \in L^2 ( \omega \times (0,T) ) $.

The problem reduces to obtaining sufficient conditions on $ \omega $ and $ T $ so that
$
F^\prime = L^2(\Omega) \times H^1_0(\Omega)
$
or, equivalently,
$
F = L^2(\Omega) \times H^{-1}(\Omega).
$

It is therefore sufficient to prove the existence of constants $ c>0, C>0 $ such that
\begin{align}\label{eqIII11}
\begin{split}
   c \| \{ \phi^0, \phi^1 \} \|^2_{L^2(\Omega) \times H^{-1}(\Omega)} \leq \int^T_0 \int_\omega |\phi|^2 dx dt \leq 
   C \| \{ \phi^0, \phi^1 \} \|^2_{L^2(\Omega) \times H^{-1}(\Omega)},
   \\ \forall \{ \phi^0, \phi^1 \} \in \mathcal{D}(\Omega) \times \mathcal{D}(\Omega).
\end{split}
\end{align}

It is well known that
$
    \| \phi \|_{L^\infty(0, T; L^2(\Omega) )} \leq 
    C \| \{ \phi^0, \phi^1 \} \|_{L^2(\Omega) \times H^{-1}(\Omega)},
$
from which the upper bound in (\ref{eqIII11}) is deduced.

It remains to prove an inverse observability inequality of the form
\begin{equation}\label{eqIII12}
    \| \{ \phi^0, \phi^1 \} \|_{L^2(\Omega) \times H^{-1}(\Omega)}
    \leq C \int^T_0 \int_\omega |\phi|^2 dx dt,
\end{equation}
which will be the subject of the next section. This will require to impose geometric conditions on the observation subdomain $\omega$.

\section{Inverse Inequality}
\label{sec:III.2}
Let $ x^0 \in \mathbb{R}^n $, and consider the partition of the boundary $ \{ \Gamma(x^0), \Gamma_\star(x^0) \} $ introduced in the previous chapter.

Consider $ \omega \subset \Omega $ to be a neighborhood of $ \overline{\Gamma(x^0)} $, so that, for an open set $ \Theta \subset \mathbb{R}^n $ such that
$
    \overline{\Gamma(x^0)} \subset \Theta,
$
we have
$
    \omega = \Omega \cap \Theta.
$

\begin{figure}
\includegraphics[width=8cm]{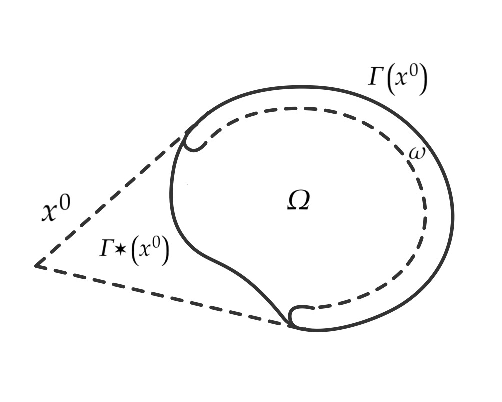}
\centering
\captionof{figure}{Geometric configuration in which the control acts on a neighborhood $\omega$ of $\Gamma(x^0)$.}
\end{figure}

In this geometric setting the following controllability result holds:
\begin{theorem}
\label{theorem:III.3.1}
Let $ \Omega $ be a bounded domain of $ \mathbb{R}^n $ with the boundary $ \Gamma $ of class $ C^2 $. Let $ x^0 \in \mathbb{R}^n $ and $ \omega \subset \Omega $ be a neighborhood of $ \overline{\Gamma(x^0)} $.

If $\hspace{0.2cm} T > T(x^0) = 2R(x^0) $, there exists a constant $ C= C(T)>0 $ such that
\begin{equation} \label{eqIII13}
   \| \phi^0 \|^2_{L^2(\Omega)} + \| \phi^1 \|^2_{H^{-1}(\Omega)} \leq C \int^T_0 \int_\omega |\phi|^2 dx dt
\end{equation}
for any solution of (\ref{eqIII03}) with data $ \{ \phi^0, \phi^1 \} \in L^2(\Omega) \times H^{-1}(\Omega) $.
\end{theorem}

\label{proof:III.theorem.3.1}
\begin{proof}We proceed in several steps.
\newline
\emph{Step 1: Reduction to a higher order inequality.}
Suppose that the following estimate holds
\begin{equation} \label{eqIII14}
   \| \phi^0 \|^2_{H^1_{0}(\Omega)} + \| \phi^1 \|^2_{L^2(\Omega)} \leq C \int^T_0 \int_\omega | \phi^\prime |^2 dx dt.
\end{equation}
It is aimed to show that (\ref{eqIII14}) implies (\ref{eqIII13}).

Given $ \{ \phi^0 , \phi^1 \} \in L^2(\Omega) \times H^{-1}(\Omega) $ define
\begin{equation}\label{eqIII15}
    X \in H^1_{0}(\Omega), \quad \Delta X = \phi^1  \text{ in } \Omega.
\end{equation}

Out of the solution $ \phi $ of (\ref{eqIII03}) with the initial data $ \{ \phi^0 , \phi^1 \} $, define
\begin{equation}\label{eqIII16}
    \psi (t) = \int^t_ 0 \phi (s) ds + X.
\end{equation}

It can be easily seen that $ \psi^\prime = \phi$ and $ \psi = \psi (x, t) $ solves the following system
\begin{equation}\label{eqIII17}
    \begin{cases}
    \psi^{\prime\prime} - \Delta \psi = 0 \hspace{1.6cm} & \text{in} \hspace{0.3cm} Q
    \\
    \psi = 0 \hspace{0.5cm} & \text{in} \hspace{0.3cm} \Sigma
    \\
    \psi(0) = X, \psi^\prime (0) = \phi^0 & \text{in} \hspace{0.3cm} \Omega.
    \end{cases}
\end{equation}

Next, apply the inequality (\ref{eqIII14}) to $ \psi $ to obtain
\begin{align} \label{eqIII18}
\begin{split}
   \|X\|^2_{H^1_{0}(\Omega)} + \| \phi^0 \|^2_{L^2(\Omega)} \leq 
   C \int^T_0 \int_\omega | \psi^\prime |^2 dx dt = C \int^T_0 \int_\omega | \phi |^2 dx dt.
\end{split}
\end{align}

On the other hand, (\ref{eqIII13}) is obtained from (\ref{eqIII18}), thanks to the classical properties of the Dirichlet Laplacian, since
\begin{equation}\label{eqIII19}
    \exists C > 0 : \| \phi^1 \|_{H^{-1}(\Omega)} \leq C \|X\|_{H^1_0}, \quad \forall \phi^1 \in H^{-1}(\Omega).
\end{equation}

Now we are ready to prove (\ref{eqIII14}) in several steps.\\
\label{sec:III.2-c}
\emph{Step 2: Estimation of the normal derivative by a local term.}
By  Theorem \ref{theorem:II.2.2} of Chapter \ref{chapter02}, it is shown that if $ T>T(x^0) $, there exists $ C=C(T)>0 $ such that
\begin{equation}\label{eqIII20}
    E_0 = \frac{1}{2} \int_\Omega \left[ \left| \nabla \phi^0 (x) \right|^2 + \left| \phi^1(x) \right|^2  \right] dx \leq
    C \int ^T_0 \int_{\Gamma(x^0)} \left| \frac{\partial \phi}{\partial \nu} \right|^2 d\Sigma.
\end{equation}

Let $  \varepsilon  > 0 $ be such that $ \varepsilon < (T - T(x^0))/2 $. Then, from (\ref{eqIII02}), the invariance of the wave equation with respect to time-translations, and  the conservation of energy, it can be deduced that
\begin{equation}\label{eqIII21}
    E_0 \leq C \int ^{T-\varepsilon}_\varepsilon \int_{\Gamma(x^0)} \left| \frac{\partial \phi}{\partial \nu} \right|^2 d\Sigma.
\end{equation}

Next, consider a vector field $ q = q(x,t) \in (W^{1,\infty}(Q))^n $ such that
\begin{equation}\label{eqIII22}
\begin{cases}
q(x,t) = \nu(x) & \hspace{0.7cm} \forall(x,t) \in \Gamma (x^0) \times (\varepsilon, T-\varepsilon)
\\
q(x,t)\cdot \nu(x) \geq 0 & \hspace{0.7cm} \forall(x,t) \in \Gamma \times (0,T)
\\
q(x,0) = q(x,T) = 0 & \hspace{0.7cm} \forall x \in \Omega 
\\
q(x,t) = 0 & \hspace{0.7cm} \forall(x,t) \in (\Omega \backslash \omega) \times (0,T).
\end{cases}
\end{equation}

The field $ q $ can be chosen as
$
    q(x,t) = h(x) n (t)
$
where $ n \in C^1([0,T]) $ and $ h = h(x) \in (W^{1, \infty} (\Omega))^n $ satisfy
\begin{equation}\label{eqIII23}
    n(0) = n(T) = 0; \hspace{0.3cm} n(t) = 1 \hspace{0.3cm} \forall t \in (\varepsilon, T- \varepsilon),
\end{equation}
and
\begin{equation}\label{eqIII24}
    \begin{cases}
    h = \nu & \hspace{1.5cm} \text{on} \hspace{0.5cm} \Gamma(x^0)
    \\h \cdot \nu \geq 0 & \hspace{1.5cm} \text{on} \hspace{0.5cm} \Gamma
    \\ h = 0 & \hspace{1.5cm} \text{in} \hspace{0.5cm} \Omega \backslash \omega.
    \end{cases}
\end{equation}
(see J.-L. Lions \cite{lions1988controlabilite}, Chapter I, Section 3.1 for the construction of this vector field).

By the identity (\ref{eqII18}) of Lemma \ref{lemma:II.2.1} in Chapter \ref{chapter02} we obtain,
\begin{align}\label{eqIII25}
\begin{split}
    \frac{1}{2} \int ^{T-\varepsilon}_\varepsilon \int_{\Gamma(x^0)} \left| \frac{\partial \phi}{\partial \nu} \right|^2 d\Sigma \leq & \frac{1}{2} \int ^{T-\varepsilon}_\varepsilon \int_\Gamma (q \cdot \nu) \left| \frac{\partial \phi}{\partial \nu} \right|^2 d\Sigma
    \\ \leq & C \int ^{T}_0 \int_\omega \left( |\phi^\prime|^2 + |\nabla \phi|^2 \right) dx dt
\end{split}
\end{align}
where $C$ is a positive constant depending on $ \| q(x,t) \|_{W^{1, \infty}(Q)}$.

The estimate (\ref{eqIII25}) is valid for all $ T > T(x^0)$. Therefore, if $ \varepsilon > 0 $ is such that $ T - 2\varepsilon > T(x^0) $, by the argument used previously and combining (\ref{eqIII21}) and (\ref{eqIII25}), it can be deduced that
\begin{equation}\label{eqIII26}
    E_0 \leq C \int ^{T-\varepsilon}_\varepsilon \int_\omega \left[ |\phi^\prime|^2 + | \nabla \phi |^2 \right] dx dt.
\end{equation}
\label{sec:III.2-d}
\emph{Step 3: The local $\nabla\phi$-term.}
Let $ \tilde{\omega} $ be a neighborhood of $ \overline{\Gamma(x^0)} $ such that
\begin{equation}\label{eqIII27}
    \Omega \cap \overline{\tilde{\omega}} \subset \omega.
\end{equation}

Since (\ref{eqIII26}) is valid for any neighborhood $ \omega $ of $ \overline{\Gamma(x^0)} $, it also holds for $\tilde{\omega}$ and, therefore
\begin{equation}\label{eqIII28}
    E_0 \leq C \int ^{T-\varepsilon}_\varepsilon \int_{\tilde{\omega}} \left[ |\phi^\prime|^2 + | \nabla \phi |^2 \right] dx dt.
\end{equation}

Let $ p = p(x,t) \in W^{1,\infty}(Q), p \geq 0 $ be such that
\begin{equation}\label{eqIII29}
\begin{cases}
p(x,t) = 1 & \hspace{0.7cm} \forall(x,t) \in \tilde{\omega} \times (\varepsilon, T-\varepsilon)
\\
p(x,t) = 0 & \hspace{0.7cm} \forall(x,t) \in (\Omega \backslash \omega) \times (0,T)
\\
p(x,0) = p(x,T) = 0 & \hspace{0.7cm} \forall x \in \Omega 
\\
|\nabla p|^2/p \in L^\infty (Q).
\end{cases}
\end{equation}
The function $p$ can be chosen as
\begin{equation}\label{eqIII30}
    p(x,t) = \rho^2(x) n(t)
\end{equation}
where $ n \in C^1([0,T]) $, and $ \rho \in W^{1, \infty}(\Omega) $ satisfies
\begin{equation}\label{eqIII31}
    \begin{cases}
    \rho(x) = 1 \hspace{2cm} \forall x \in \tilde{\omega}
    \\
    \rho(x) = 0 \hspace{2cm} \forall x \in \Omega \backslash \omega.
    \end{cases}
\end{equation}

Now, multiplying equation (\ref{eqIII03}) by $ p\phi $ and integrating by parts, the following identity holds immediately
\begin{equation}\label{eqIII32}
    \int_Q p |\nabla \phi|^2 = \int_Q p|\phi^\prime|^2  
    - \int_\Omega p \phi^\prime \phi\bigg|^T_0 
    + \int_Q p^\prime \phi^\prime \phi - \int_Q \nabla p \cdot \nabla \phi \phi.
\end{equation}

Considering (\ref{eqIII29}) and (\ref{eqIII32}) leads to
\begin{equation}\label{eqIII33}
    \int_Q p |\nabla \phi|^2 dx dt \leq \int^T_0 \int_\omega \left[ |\phi^\prime|^2 + |\phi|^2 \right] dx dt
    + \left| \int_Q (\nabla p \cdot \nabla \phi) \phi \right|.
\end{equation}

On the other hand,
\begin{equation}\label{eqIII34}
    \left| \int_Q (\nabla p \cdot \nabla \phi) \phi \right| \leq
    \frac{1}{2} \int_Q p |\nabla \phi|^2 + \frac{1}{2} \int_Q \frac{|\nabla p|^2}{p} |\phi|^2.
\end{equation}

From (\ref{eqIII33}) and (\ref{eqIII34}), using (\ref{eqIII29}), it follows that
\begin{equation}\label{eqIII35}
    \int ^{T-\varepsilon}_\varepsilon \int_\omega | \nabla \phi|^2 \leq C \int^T_0 \int_\omega \left[ |\phi^\prime|^2 + | \phi |^2 \right].
\end{equation}

Finally, combine (\ref{eqIII28}) and (\ref{eqIII35}) to get
\begin{equation}\label{eqIII36}
    E_0 \leq C \int^T_0 \int_{\omega} |\phi^\prime|^2 dx dt + C \int_Q |\phi |^2 dx dt.
\end{equation}
\label{sec:III.2-e}
\emph{Step 4: The lower order $L^2$-term.}
To prove (\ref{eqIII14}) (and therefore the theorem), it suffices  to show the existence of $ C = C(T) > 0 $ such that \begin{equation}\label{eqIII37}
    \int_Q |\phi |^2 dx dt \leq C \int^T_0 \int_\omega |\phi^\prime|^2 dx dt
\end{equation}
holds for every solution of (\ref{eqIII03}).

We proceed arguing by contradiction, using the  compactness-uniqueness argument. If (\ref{eqIII37}) is not satisfied, there exists a sequence of solutions $ \{ \phi_n \} $ of (\ref{eqIII03}) corresponding to the initial data $ \{ \phi^0_n, \phi^1_n \} \subset H^1_0(\Omega) \times L^2(\Omega) $ such that
\begin{equation}\label{eqIII38}
    \int_Q |\phi_n |^2 dx dt = 1 \hspace{0.5cm} \forall n \in \mathbb{N},
\end{equation}
\begin{equation}\label{eqIII39}
    \int^T_0 \int_\omega |\phi^{\prime}_n |^2 dx dt \rightarrow 0,  n \rightarrow \infty.
\end{equation}

From (\ref{eqIII36}), (\ref{eqIII38}) and (\ref{eqIII39}) it follows that $ \{ \phi^0_n, \phi^1_n \}$ is bounded in $ H^1_0(\Omega) \times L^2(\Omega) $, and therefore,
\begin{equation}\label{eqIII40}
    \{ \phi_n \} \hspace{0.3cm} \text{is bounded in} \hspace{0.2cm} H^1(Q).
\end{equation}

From (\ref{eqIII40}) and the compactness of the embedding $ H^1(Q) \subset L^2(Q) $ we have that
\begin{equation}\label{eqIII41}
    \{ \phi_n \} \hspace{0.3cm} \text{is relatively compact in} \hspace{0.2cm}  L^2(Q).
\end{equation}

By extracting a suitable subsequence (which we continue to denote by $ \{ \phi_n \} $ to simplify the notation), the following is deduced
\begin{equation}\label{eqIII42}
    \phi_n \rightarrow \phi \hspace{0.3cm} \text{weakly in} \hspace{0.2cm} H^1(Q),
\end{equation}
\begin{equation}\label{eqIII43}
    \phi_n \rightarrow \phi \hspace{0.3cm} \text{strongly in} \hspace{0.2cm}  L^2(Q),
\end{equation}
where the limit $ \phi = \phi (x,t) $ is a solution of (\ref{eqIII03}). On the other hand, by (\ref{eqIII38}), (\ref{eqIII39}), (\ref{eqIII42}) and (\ref{eqIII43}), $ \phi$ satisfies
\begin{equation}\label{eqIII44}
    \phi^\prime = 0 \hspace{0.3cm} \text{in} \hspace{0.3cm} \omega \times (0,T),
\end{equation}
\begin{equation}\label{eqIII45}
    \| \phi \|_{L^2(Q)} = 1.
\end{equation}

Notice that the function $ \psi = \phi^\prime $ is also a solution of (\ref{eqIII03}), and it satisfies
\begin{equation}\label{eqIII46}
    \psi = 0 \hspace{0.3cm} \text{in} \hspace{0.3cm} \omega \times (0,T).
\end{equation}

As $ T > T(x^0) \geq $ diameter of $ \Omega $, it follows from Holmgren's Uniqueness Theorem that
$
    \psi \equiv 0  
$ in $Q$.
This means that
$
    \phi = \phi (x)
$
is a stationary solution of (\ref{eqIII03}). Hence, this implies that
$
    \phi \equiv 0,
$
which contradicts (\ref{eqIII45}).
\end{proof}

\section{The Main Result}
\label{sec:III.3}
The following result has been proved:
\begin{theorem}
\label{theorem:III.3.2}
Let $ \Omega $ be a bounded domain of $ \mathbb{R}^n $ with the boundary $ \Gamma $ of class $C^2$. Let $ x^0 \in \mathbb{R}^n $ and $ \omega \subset \Omega $ be a neighborhood of $ \overline{\Gamma(x^0)} $.

Suppose that $ T> T(x^0) = 2R(x^0) $.

Then, for each initial data $ \{ y^0, y^1 \} \in H^1_0(\Omega) \times L^2(\Omega) $, there is a control $ h \in L^2( \omega \times (0,T) ) $ such that the solution $y(x,t)$ of (\ref{eqIII01}) satisfies (\ref{eqIII02}).
\end{theorem}

\begin{proof}
\label{proof:III.theorem.3.2}
It is an immediate consequence of the adaptation of the HUM method developed in Section \ref{sec:III.1} and of the observability estimate (\ref{eqIII13}), which implies that the constructed Hilbert space $F$ coincides with $ L^2(\Omega) \times H^{-1}(\Omega) $, by virtue of the equivalence of norms.
\end{proof}

\section{Variable Coefficients in $1-d$ Domains}
\label{sec:III.4}
Consider the following $1-d$ wave equation with space-dependent variable coefficient
\begin{equation}\label{eqIII47}
    \begin{cases}
    y^{\prime\prime} - (a(x) y_x)_x = h \chi_\omega \hspace{0.3cm} & \text{in} \hspace{0.3cm} \Omega \times (0,T)
    \\
    y(0,t) = y(1,t) = 0 \hspace{0.3cm} & \text{for} \hspace{0.3cm} t \in (0,T)
    \\
    y(x,0) = y^0 (x), y^\prime(x,0) = y^1(x) \hspace{0.3cm} & \text{for} \hspace{0.3cm} x \in \Omega

    \end{cases}
\end{equation}
with $ \Omega = (0,1) $ and $ \omega = ( l_1, l_2 ) \subset \Omega. $

Combining the techniques of the previous section together with the ones of Section \ref{sec:II.7}, the following result is obtained.

\begin{theorem}
\label{theorem:III.3.3}
Suppose the coefficient $ a \in BV(0,1) $ and it satisfies
\begin{equation}\label{eqIII48}
    a(x) \geq a_0 > 0, \hspace{0.5cm} \forall x \in (0,1).
\end{equation}

Then, if $ T > 2R $ with $ R = R (l_1, l_2) = \max \left( l_1/\sqrt{a_0}, (1-l_2)/\sqrt{a_0} \right), $ for each initial data pair $ \{ y^0, y^1 \} \in H^1_0(0,1) \times L^2(0,1) $ there exists a control $ h \in L^2( (l_1, l_2) \times (0,T) ) $ such that the solution $ y = y(x,t) $ of (\ref{eqIII47}) satisfies (\ref{eqIII02}).
\end{theorem}

\label{sec:III.4-b}
\begin{proof}
By the HUM method
 the problem is reduced to obtaining the observability estimate \begin{equation}\label{eqIII49}
    \| \phi^0 \|^2_{L^2(\Omega)} + \| \phi^1 \|^2_{H^{-1}(\Omega)}
    \leq C \int^T_0 \int^{l_2}_{l_1} |\phi|^2 dx dt
\end{equation}
for all solutions  $\phi$ of the homogeneous system
\begin{equation}\label{eqIII50}
    \begin{cases}
    \phi^{\prime\prime} - (a(x) \phi_x)_x = 0 \hspace{0.3cm} & \text{in} \hspace{0.3cm} (0,1) \times (0,T)
    \\
    \phi(0,t) = \phi(1,t) = 0 \hspace{0.3cm} & \text{on} \hspace{0.2cm} (0,T)
    \\
    \phi(x,0) = \phi^0 (x), \phi^\prime(x,0) = \phi^1(x) & \text{in} \hspace{0.2cm} (0,1).
    \end{cases}
\end{equation}

Proceeding as in Step 1 of the proof of Theorem \ref{theorem:III.3.1}, it is observed that (\ref{eqIII49}) is a consequence of
\begin{equation}\label{eqIII51}
    \| \phi^0 \|^2_{H^1_0(\Omega)} + \| \phi^1 \|^2_{L^2(\Omega)}
    \leq C \int^T_0 \int^{l_2}_{l_1} |\phi^\prime|^2 dx dt.
\end{equation}

Now, (\ref{eqIII51}) is proved in several steps.\\
\smallskip
\emph{Step 1: A first observability inequality.}\label{sec:III.4-c}
The method of proof of the observability estimate (\ref{eqII73}) of Chapter \ref{chapter02}  allows to show that, for any $ \xi \in [0,1], $ if
\begin{equation}\label{eqIII52}
    T > \frac{2}{\sqrt{a_0}} \xi,
\end{equation}

\begin{equation}\label{eqIII53}
\begin{split}
    \int^\xi_0 \int^{T + \frac{x-\xi}{\sqrt{a_0}}}_{\frac{\xi-x}{\sqrt{a_0}}} 
    \left[ |\phi^\prime|^2 + a(x) |\phi_x|^2 \right] dt dx
    \leq C \int^T_0 \left[ |\phi^\prime (\xi,t)|^2 + a(\xi) |\phi_x(\xi,t)|^2 \right] dt.
\end{split}
\end{equation}

Analogously, if
\begin{equation}\label{eqIII54}
    T > \frac{2}{\sqrt{a_0}} (1-\xi),
\end{equation}
we have
\begin{equation}\label{eqIII55}
\begin{split}
    \int^\xi_0 \int^{T - \frac{x-\xi}{\sqrt{a_0}}}_{\frac{x-\xi}{\sqrt{a_0}}} 
    \left[ |\phi^\prime|^2 + a(x) |\phi_x|^2 \right] dt dx
    \leq C \int^T_0 \left[ |\phi^\prime (\xi,t)|^2 + a(\xi) |\phi_x(\xi,t)|^2 \right] dt.
\end{split}
\end{equation}

Furthermore, the constant $C$ in both estimates (\ref{eqIII53}) and (\ref{eqIII55}) depends continuously on $\xi$ and $T$.

Since $T>2R$, there exists $\varepsilon >0$ such that
\begin{equation}\label{eqIII56}
    T > \frac{2}{\sqrt{a_0}} \max(l_1+\varepsilon, 1 + \varepsilon -l_2).
\end{equation}

From (\ref{eqIII56}), it is deduced that (\ref{eqIII53}) and (\ref{eqIII55}) hold for all $ \xi \in (l_1, l_1 + \varepsilon) $ and $ \xi \in (l_2 - \varepsilon, l_2) $, respectively. Integrating the inequalities (\ref{eqIII53}) and (\ref{eqIII55}) with respect to $ \xi \in (l_1, l_1 + \varepsilon) $ and $ \xi \in (l_2 - \varepsilon, l_2) $), respectively, the following inequalities are immediate
\begin{equation}\label{eqIII57}
    \begin{split}
    \int^{l_1}_0 \int^{T + \frac{(x-l_1)}{\sqrt{a_0}}}_{\frac{l_1-x}{\sqrt{a_0}}} 
    \left[ |\phi^\prime|^2 + a(x) |\phi_x|^2 \right] dt dx
    \leq C \int^T_0 \int^{l_1+\varepsilon}_{l_1} 
    \left[ |\phi^\prime|^2 + a(x) |\phi_x|^2 \right] dx dt,
    \end{split}
\end{equation}
and
\begin{equation}\label{eqIII58}
\begin{split}
    \int^1_{l_2} \int^{T - \frac{(x-l_2)}{\sqrt{a_0}}}_{\frac{x-l_2}{\sqrt{a_0}}} 
    \left[ |\phi^\prime|^2 + a(x) |\phi_x|^2 \right] dt dx
    \leq C \int^T_0 \int^{l_2}_{l_2-\varepsilon} \left[ |\phi^\prime|^2 + a(x) |\phi_x|^2 \right] dx dt.
\end{split}
\end{equation}

Now utilizing (\ref{eqIII57}) and (\ref{eqIII58}), we have
\begin{equation}\label{eqIII59}
    \begin{split}
    & (T-2R) \int^1_0 \left[ |\phi^1(x)|^2 + a(x) |\phi^0_x(x)|^2 \right] dx
    \\ & = \int^{T-R}_{R} \int^1_0 
    \left[ |\phi^\prime|^2 + a(x) |\phi_x|^2 \right] dx dt
    \leq C \int^T_0 \int^{l_2}_{l_1} \left[ |\phi^\prime|^2 + a(x) |\phi_x|^2 \right] dx dt
    \end{split}
\end{equation}
or, equivalently,
\begin{equation}\label{eqIII60}
    \| \phi^0 \|^2_{H^1_0(0,1)} + \| \phi^1 \|^2_{L^2(\Omega)}  \leq C \int^T_0 \int^{l_2}_{l_1} \left[ |\phi^\prime|^2 + a(x) |\phi_x|^2 \right] dx dt.
\end{equation}
\label{sec:III.4-d}
\emph{Step 2: An improved observability inequality.}
The estimate (\ref{eqIII60}) is valid for all $ (l_1, l_2) \subset (0,1) $ as long as $ T > 2R (l_1, l_2) $.

Taking $ \varepsilon >0 $ small enough so that
$
    T - 2\varepsilon > 2R(l_1 - \varepsilon, l_2 - \varepsilon),
$
we have
\begin{equation}\label{eqIII61}
    \begin{split}
    \| \phi^0 \|^2_{H^1_0(0,1)} + \| \phi^1 \|^2_{L^2(0,1)}
    \leq C \int^{T-\varepsilon}_{\varepsilon} \int^{l_2-\varepsilon}_{l_1+\varepsilon}
    \left[ |\phi^\prime|^2 + a(x) |\phi_x|^2 \right] dx dt.
    \end{split}
\end{equation}

The argument of Step 3 of the proof of Theorem \ref{theorem:III.3.1} easily allows to prove that
\begin{equation}\label{eqIII62}
    \begin{split}
    \int^{T-\varepsilon}_{\varepsilon} \int^{l_2-\varepsilon}_{l_1+\varepsilon}
    a(x) |\phi_x|^2 dx dt
    \leq C \int^T_0 \int^{l_2}_{l_1} |\phi^\prime|^2 dx dt + C\|\phi\|^{2}_{L^2(Q)}.
    \end{split}
\end{equation}

Now, combining (\ref{eqIII61}) and (\ref{eqIII62}) leads to
\begin{equation}\label{eqIII63}
    \begin{split}
    \|\phi^0\|^2_{H^1_{0}(0,1)} + \|\phi^1\|^2_{L^2(0,1)} 
    \leq C \int^T_0 \int^{l_2}_{l_1} |\phi^\prime|^2 dx dt + C\|\phi\|^{2}_{L^2(Q)}.
    \end{split}
\end{equation}

Therefore, it is simply enough to prove the inequality
\begin{equation}\label{eqIII64}
    \|\phi\|^{2}_{L^2(Q)} \leq C \int^T_0 \int^{l_2}_{l_1} |\phi^\prime|^2 dx dt.
\end{equation}
\label{sec:III.4-e}
\emph{Step 3. Compactness-uniqueness.}
To prove (\ref{eqIII64}) we proceed by a contradiction argument, as in  the proof of Step 4 of Theorem \ref{theorem:III.3.1}, using the compactness-uniqueness principle.

If (\ref{eqIII64}) does not hold, there exists a solution $ \phi = \phi (x,t) $ of (\ref{eqIII50}) such that
\begin{equation}\label{eqIII65}
    \|\phi\|_{L^2(Q)} = 1,
\end{equation}

\begin{equation}\label{eqIII66}
    \phi^\prime (x,t) = 0, \hspace{0.5cm} \forall (x,t) \in (l_1, l_2) \times (0,T).
\end{equation}

Therefore, $ \psi = \phi^\prime $ is a solution of (\ref{eqIII50}) satisfying that
\begin{equation}\label{eqIII67}
    \psi \in C ([0,T]; L^2(0,1))
\end{equation}
and
\begin{equation}\label{eqIII68}
    \psi = 0 \hspace{0.3cm} \text{in} \hspace{0.3cm} (l_1, l_2) \times (0,T).
\end{equation}

It follows from (\ref{eqIII68}) that
$
    \psi^\prime = 0$ in $ (l_1, l_2) \times (0,T)$.
This, together with (\ref{eqIII67}) and (\ref{eqIII63}), allows to conclude that
$
    \psi(0) \in H^1_0 (0,1)$, $\psi^\prime (0)\in L^2(0,1),
$
and therefore
\begin{equation}\label{eqIII69}
    \psi \in C([0,T]; H^1_0 (0,1)) \cap C^1([0,T]; L^2(0,1)).
\end{equation}

As $\psi$ belongs to the class of finite energy solutions  in (\ref{eqIII69}) (a fact that was not previously known because $\psi = \phi^\prime $ and, therefore, in principle, it is one derivative less regular than $\phi$), the estimates (\ref{eqIII53}) and (\ref{eqIII55}) apply to $ \psi $, and, it follows from (\ref{eqIII68}) that
\begin{equation*}
    \psi(x,t) = 0, \hspace{0.3cm} \forall (x,t) \in (0,1) \times (R,T-R),
\end{equation*}
with $R = R(l_1, l_2)$ as in the statement of Theorem 3.3, and, by the conservation of energy, $ \psi \equiv 0$.

Therefore, $ \phi = \phi (x) $ is a stationary solution of (\ref{eqIII50}) satisfying
\begin{equation*}
    -(a(x) \phi_x)_x = 0 \hspace{0.3cm} \text{in} \hspace{0.3cm} (0,1), \hspace{0.3cm} \phi \in H^1_0(0,1).
\end{equation*}
But then, necessarily, $\phi \equiv 0$, contradicting (\ref{eqIII65}).
\end{proof}

\section{Comments}
\label{sec:III.5}
\begin{enumerate}
\item {\bf Wave equations with variable coefficients.} As seen in Theorems \ref{theorem:III.3.1} and \ref{theorem:III.3.2}, boundary controllability implies controllability with controls with support on a neighborhood of the boundary.

This type of argument can be easily extended to equations with variable coefficients of the type
\begin{equation}\label{eqIII70}
    \begin{cases}
    y^{\prime\prime} - \operatorname{div}(a(x) \nabla y) = h \chi_\omega \hspace{1.6cm} & \text{in} \hspace{0.3cm} Q
    \\
    y = 0 \hspace{0.5cm} & \text{in} \hspace{0.3cm} \Sigma
    \\
    y(0) = y^0, y^\prime (0) = y^1 & \text{in} \hspace{0.3cm} \Omega
    \end{cases}
\end{equation}
with $ a \in W^{1, \infty}(\Omega) $ such that
\begin{equation}\label{eqIII71}
a(x) \geq a_0 > 0 \hspace{0.5cm} \forall x \in \Omega.
\end{equation}
Assume that the coefficient $a=a(x)$ is such that, for some   open subset $ \Gamma_0 \subset \Gamma $ of the boundary and $ T_0 > 0 $, the following estimate holds
\begin{equation}\label{eqIII72}
    \|\phi^0 \|^2_{H^1_0 (\Omega)} + \|\phi^1 \|^2_{L^2(\Omega)}
    \leq C \int_{\Sigma_0} \left| \frac{\partial \phi}{\partial \nu} \right|^2 d\Sigma
\end{equation}
 for all $ T > T_0 $, with $ \Sigma_0 = \Gamma_0 \times (0,T) $, for any solution of
\begin{equation}\label{eqIII73}
    \begin{cases}
    \phi^{\prime\prime} - \operatorname{div}(a(x) \nabla \phi) = 0 \hspace{1.6cm} & \text{in} \hspace{0.3cm} Q
    \\
    \phi = 0 \hspace{0.5cm} & \text{in} \hspace{0.3cm} \Sigma
    \\
    \phi(0) = \phi^0, \phi^\prime (0) = \phi^1 \hspace{0.5cm} & \text{in} \hspace{0.3cm} \Omega.
    \end{cases}
\end{equation}
Recall that such a boundary observability inequality requires some structural conditions on the coefficient $a=a(x)$, assuring that all rays of geometric optics reach the boundary $\Gamma $ and, more specifically, $ \Gamma_0$.

The following result holds:

\begin{theorem}
\label{theorem:III.3.4}
Let $ \Omega $ be a bounded domain of $ \mathbb{R}^n $ with boundary $ \Gamma $ of class $ C^2 $. Let the coefficient $a=a(x)$, $ \Gamma_0 \subset \Gamma $ and $ T_0 > 0 $ be such that for all $ T>T_0 $ the estimate (\ref{eqIII72}) is satisfied.

Suppose that $ \omega $ a neighborhood of $ \overline{\Gamma_0} $ in $ \Omega $. Then, if $ T>T_0 $, for every initial data $ \{ y^0,y^1 \} \in H^1_{0} (\Omega) \times L^2(\Omega) $, there is a control $ h \in L^2(\omega \times (0,T)) $ such that the solution $ y = y (x,t) $ of (\ref{eqIII70}) satisfies (\ref{eqIII02}).
\end{theorem}

\begin{proof}
\label{proof:III.theorem.3.4}
We present a short sketch.

The HUM method reduces the proof of this theorem to obtaining the estimate (\ref{eqIII13}) for the solutions of (\ref{eqIII73}).
As in Step 1 of the proof of Theorem \ref{theorem:III.3.1}, it is observed that proving (\ref{eqIII14}) is enough.
The arguments used for Steps 2 and 3 of the proof of Theorem \ref{theorem:III.3.1} allow us to obtain (\ref{eqIII36}) easily.
Therefore, the problem reduces again to proving (\ref{eqIII37}).
The argument of Step 4 of the proof of Theorem \ref{theorem:III.3.1} does not literally apply in this case since the Holmgren's Theorem is only valid for equations with analytic coefficients.

However, there are two alternative proofs of (\ref{eqIII37}) that do apply in the present case:
\begin{enumerate}
\item The method introduced by V. Komornik \cite{komornik1988methode}, based on the use of non-harmonic Fourier series. This method  applies systematically to evolution equation of the type
$
    \phi^{\prime \prime} + A \phi = 0
$
where $A$ is a self-adjoint elliptic operator of order $ \geq 2 $, and it allows to absorb the lower order terms in estimates of the type (\ref{eqIII36}).

\item The arguments of C. Bardos, G. Lebeau and J. Rauch \cite{bardos1992sharp} reduce (\ref{eqIII37}) to a unique continuation property for the solutions of an elliptic problem of the type
\begin{equation}\label{eqIII74}
    -\operatorname{div}(a(x) \nabla \phi) + \lambda\phi = 0 \hspace{0.3cm} \text{in} \hspace{0.3cm} \Omega
\end{equation}
where $ \lambda \in \mathbb{C} $ and $ \phi = \phi (x) $ a complex valued function.

Then, the question is  reduced to proving that every solution $ \phi \in H^2(\Omega) $  of (\ref{eqIII74}) satisfying
$
    \phi = 0 $ in $\omega$, 
 necessarily vanishes everywhere, which is an immediate consequence of the unique continuation results for elliptic equations obtained through Carleman inequalities (see  \cite{hormander2007analysis}). This is actually a by now classical result on the unique continuation for  the eigenfunctions of elliptic equations, that can be achieved and broadly extended using Carleman inequalities.
\end{enumerate}
\end{proof}

Several remarks and geometric considerations are in order:
\begin{enumerate}
\item {\bf Geometric considerations.} Theorem \ref{theorem:III.3.3} shows that the wave equation in $1-d$ is exactly controllable in $ H^1_{0}(\Omega) \times L^2(\Omega) $ with controls in $ L^2(\omega \times (0,T)) $ for any non-empty open subset $ \omega $ of $ \Omega $.

This is a genuinely $1-d$ result that does not generalize to dimensions $ n>1 $. 

We refer to A. Haraux \cite{haraux1988controlabilite} and J. Lagnese \cite{lagnese1983control} for the analysis of the internal controllability of the constant coefficient wave equation, using Fourier series, when $ \Omega $ is a square and a disk of $ \mathbb{R}^2 $ respectively.
 
 For instance, the spectral analysis in \cite{lagnese1983control} shows that if $ \Omega = B(0,1) \subset \mathbb{R}^n $ and $ \omega = B(0,r) $ with $ r<1 $, the wave equation is not exactly controllable in $ H^1_{0}(\Omega) \times L^2(\Omega) $ with controls $ L^2(\omega \times (0,T)) $ for any $ T>0 $. This is due to the whispering gallery phenomenon, i. e.  the existence of periodic rays of geometric optics  endlessly trapped in the region $ r < |x| < 1 $, without never reaching the control region $ \omega $. In  \cite{lagnese1983control} it is shown that the class of controllable initial data  with $L^2$-controls is smaller than the energy space $ H^1_{0}(\Omega) \times L^2(\Omega) $ but contains, for example, all the data with a finite number of non-zero Fourier coefficients.
 
The results in \cite{haraux1988controlabilite} show that, when $ \overline{\omega} \subset \Omega $, (and, therefore, the GCC is not fulfilled), to control the initial data in $ H^1_{0}(\Omega) \times L^2(\Omega) $, it is necessary to use controls in a class of analytic functionals supported in $ \omega \times (0,T) $. 

However, as we have seen in Section \ref{sec:III.1}, thanks to the HUM method and Holmgren's Uniqueness Theorem, for all non-empty open subset  $ \omega \subset \Omega $ there exists an initial time $ T_0 = T_0(\omega, \Omega) >0 $ such that if $T>T_0$ the exact controllability of (\ref{eqIII01}) with controls in $ L^2(\omega \times (0,T)) $ in a space of initial data such as $ \{ y^1, -y^0 \} \in F^\prime $ holds. The characterization of the space $F^\prime $, or, equivalently, $F$, when $ F \neq L^2(\Omega) \times H^{-1}(\Omega) $, is an open problem to a large extent. We refer to \cite{lebeau1992analytique}, \cite{robbiano1995fonction} and \cite{allibert1999analytic} for some further developments in this context.

As mentioned in Sections \ref{sec:II.6} and \ref{sec:II.7},  estimate (\ref{eqIII72}) is satisfied for  solutions of (\ref{eqIII73}), if $ \Omega $ is $ C^\infty$ smooth, $a \in C^\infty (\overline{\Omega}) $ and the pair $ \{ \Gamma_0, T \} $ is such that the GCC is satisfied, i. e., every ray of geometric optics intersects $ \Gamma_0 $ at a non-diffractive point in a time shorter than $ T $ (cf. \cite{bardos1992sharp}). Note however that for (\ref{eqIII72}) to hold, the coefficient $a(x)$ has to fulfill some geometric conditions ensuring that all rays of geometric optics reach the exterior boundary in a uniform time. In case some of these rays would be trapped inside $\Omega$, without never reaching its exterior boundary,  the boundary observability inequality (\ref{eqIII72}) would not hold and Theorem \ref{theorem:III.3.4} above would be useless. Such pathological examples have been built and discussed in \cite{macia2002zuazua}. Note that for variable coefficient wave equations rays are not straight lines anymore and, in practice, it might be hard to check whether all of them reach $\omega$ in a uniform time.

\item {\bf Minimal norm controls.} The argument of Section \ref{sec:II.7} allows to prove that, if exact controllability in $ H^1_{0}(\Omega) \times L^2(\Omega) $ with controls in $ L^2(\omega \times (0,T)) $ holds, for each $ \{ y^0, y^1 \} \in H^1_{0}(\Omega) \times L^2(\Omega) $ the set of admissible controls, defined as 
$
U_{ad} = \{ h \in L^2(\omega \times (0,T)): y \text{ solution of (\ref{eqIII70}) satisfies (\ref{eqIII02})} \}, 
$
 contains infinitely many elements.

The control $h$ provided by the HUM method is the only optimal control of minimal norm satisfying
\begin{equation*}
    \| h \|^2_{L^2(\omega \times (0,T))} = \min_{g \in U_{ad}} \|g\|^2_{L^2(\omega \times (0,T))}.
\end{equation*}

\item {\bf Optimal control time.} Theorem \ref{theorem:III.3.2} shows that if $ \omega $ is a neighborhood of $ \Gamma(x^0) $ and $ T > T(x^0) $, then the wave equation is exactly controllable in $ H^1_{0}(\Omega) \times L^2(\Omega) $ with controls in $ L^2(\omega \times (0,T)) $.

This result is not optimal with respect to the controllability time. Indeed, thanks to \cite{bardos1992sharp}, we know that the controllability time decreases as the size of $ \omega $ increases. Estimates on the time of control  improve when employing cutting-off arguments and working in the domain $\Omega \backslash \omega$, the complementary set of $\omega$ in $\Omega$.

Note that if $ \omega = \Omega $, exact controllability holds in an arbitrarily small time (cf. J.-L. Lions \cite{lions1988controlabilite}, Section VII 2.2., where the existence of controls $ h \in L^2(Q) $ is proved). Working on the domain $ \Omega \backslash \omega $ allows to prove that the controllability time tends to zero when $\omega$ increasingly converges to the whole $\Omega$ (see \cite{kim1992exact}).

\item {\bf Singular limit towards boundary control.} Consider a sequence of open sets of the form
$
    \omega_\varepsilon = \{ x \in \Omega: d(x, \Gamma(x^0)) < \varepsilon \},
$
where $ d(x, \Gamma(x^0)) $ represents the distance from $x$ in $ \mathbb{R}^n $ to the set $ \Gamma(x^0) $. This is an $\varepsilon$-"width" neighborhood of $ \Gamma(x^0) $ .

Given $ T>T(x^0) $, controls $ h_\varepsilon \in L^2( \omega_\varepsilon \times (0,T)) $ can be constructed so that the system (\ref{eqIII01}) converges to the wave equation with Dirichlet boundary control  (see \cite{fabre1990comportement}, \cite{fabre1990equation}).

\item {\bf Other models.} The method of proof of Theorem \ref{theorem:III.3.1}  can be employed in a widder context and, in particular, for wave equations with  other boundary conditions or even for plate models (see \cite{zuazua1988controlabilite}).

For instance, for the system
\begin{equation}\label{eqIII75}
    \begin{cases}
    y^{\prime\prime} + \Delta^2 y = h \chi_\omega \hspace{1.6cm} & \text{in} \hspace{0.3cm} Q
    \\
    y = \partial y / \partial \nu = 0 \hspace{0.5cm} & \text{in} \hspace{0.3cm} \Sigma
    \\
    y(0) = y^0, y^\prime (0) = y^1 \hspace{0.5cm} & \text{in} \hspace{0.3cm} \Omega,
    \end{cases}
\end{equation}
the following holds: {\it "Let $ \Omega $ be a bounded domain of $ \mathbb{R}^n $ with boundary of class $ C^3 $ and $ \omega $ be a neighborhood of a subset of boundary of type $ \overline{\Gamma(x^0)} $. Then, for all $ T>0 $ and $ \{ y^0, y^1 \} \in H^2_{0}(\Omega) \times L^2(\Omega) $ there exists a control $ h \in L^2(\omega \times (0,T)) $ such that the solution of (\ref{eqIII75}) satisfies (\ref{eqIII02})."} 

This is the analog of the previous result on the wave equation, but, this time, it is valid in an arbitrarily small time.

In \cite{machtyngier1990control} the techniques in this chapter have been adapted  to study the internal controllability of the Schrödinger equation.

\item {\bf Fourier series and plate models.} Note however that this method only applies when the support of the control $\omega$ is a neighborhood of a suitable subset of the boundary. Substantial further work is needed to consider other control subsets $\omega$.  For instance, in \cite{jaffard1988controle} it is shown, by Fourier series methods, that if $ \Omega $ is a square of $ \mathbb{R}^2 $, the plate equation above can be exactly controlled in $ H^2_{0}(\Omega) \times L^2(\Omega) $ by controls in $ L^2(\omega \times (0,T)) $, $ \omega $ being  strictly contained in the subset of $ \Omega $, i.e. $ \overline{\omega} \subset \Omega $. This result differs in an essential way with the previously mentioned ones regarding the wave equation. Indeed, in this case  the controllability property holds even if the support of the controls, $\omega$, does not need to fulfill the geometric control condition. This shows that the GCC is a sufficient  but not necessary condition for the control of the plate equation.
\end{enumerate}
\end{enumerate}


%
%

\chapter{Internal Controllability of the Semilinear Wave Equation}
\label{chapter04} 

\abstract{
In this chapter the exact controllability of the semilinear wave equation, for globally Lipschitz nonlinearities, and with controls acting in the interior  is studied. The problem is solved by means of linearization and fixed point techniques.
}

\section{Description of the HUM method}
\label{sec:IV.1}
Let $ \Omega $ be a bounded domain of $ \mathbb{R}^n $, $n\ge 1$,  with boundary $ \Gamma $ of class $ C^2 $. Let $ \omega \subset \Omega $ be a neighborhood of the whole boundary $ \Gamma =\partial\Omega $.

Consider the semilinear wave equation
\begin{equation}\label{eqIV01}
    \begin{cases}
    y^{\prime\prime} - \Delta y + f(y) = h \chi_\omega \hspace{1.6cm} & \text{in} \hspace{0.3cm} Q
    \\
    y = 0 \hspace{0.5cm} & \text{in} \hspace{0.3cm} \Sigma
    \\
    y(0) = y^0, y^\prime (0) = y^1 \hspace{0.5cm} & \text{in} \hspace{0.3cm} \Omega
    \end{cases}
\end{equation}
where $f$ is a locally Lipschitz function.

The exact controllability problem of (\ref{eqIV01}) is formulated as follows: \emph{Find $ T>0 $ such that for every pair of initial and final data $ \{ y^0, y^1 \} , \{ z^0, z^1 \} \in H^1_0 (\Omega) \times L^2(\Omega) $, there exists a control $ h \in L^2(\omega \times (0,T)) $ such that the solution $ y = y(x,t) $ of (\ref{eqIV01}) satisfies}
\begin{equation} \label{eqIV02}
    y(T) = z^0, y^\prime(T) = z^1.
\end{equation}

\begin{remark}
\label{remark:IV.4.1}
Unlike in the linear context, in this case, the problem is not reduced to considering the particular case $ z^0 = z^1 = 0 $. We have to prove that any initial data can be driven to any final data. Recall that for the linear wave equation the problem  was reduced to consider the null target because of the time-reversibility and linearity of the system. The semilinear wave equation is time-reversible. But, obviously, the needed linearity of the model fails. We need therefore to consider arbitrary initial and final data.
$ \hfill \square $

We will assume that $f$ satisfies the following growth condition at infinity
\begin{equation}\label{eqIV03}
    |f^\prime(s)| \leq C (1 + |s|^p) \hspace{0.5cm} \forall s \in \mathbb{R}
\end{equation}
with
\begin{equation}\label{eqIV04}
    p(n - 2) \leq 2,
\end{equation}
so that (\ref{eqIV01}) admits an unique local finite energy solution, continuous with respect to time with values in $ H^1_0(\Omega) $ and of class $C^1$ with values in $L^2(\Omega)$.

On the other hand, we assume that $ T>T_0 $ where $ T_0 >0 $ is the exact controllability time of the linear equation with $ f \equiv 0 $. Therefore, we are looking for a result ensuring that, when the linear equation is exactly controllable, the same occurs for the semilinear equation with nonlinearities $f$ satisfying the conditions above and in a class to be determined.

By Theorem \ref{theorem:III.3.1} of Chapter \ref{chapter03}, it is known that if $T>T_0$ with
\begin{equation}\label{eqIV05}
    T_0 = \text{diameter of } \Omega,
\end{equation}
the following estimates hold (they are actually equivalent)
\begin{equation}\label{eqIV06}
\|\phi^0\|^2_{L^2(\Omega)} + \|\phi^1\|^2_{H^{-1}(\Omega)} \leq C \int^T_0 \int_\omega |\phi|^2 dx dt,
\end{equation}
and
\begin{equation}\label{eqIV07}
\|\phi^0\|^2_{H_0^1(\Omega)} + \|\phi^1\|^2_{L^2(\Omega)} \leq C \int^T_0 \int_\omega |\phi^\prime|^2 dx dt,
\end{equation}
for the solutions of the homogeneous linear equation
\begin{equation}\label{eqIV08}
    \begin{cases}
    \phi^{\prime\prime} - \Delta \phi = 0 \hspace{1.6cm} & \text{in} \hspace{0.3cm} Q
    \\
    \phi = 0 \hspace{0.5cm} & \text{in} \hspace{0.3cm} \Sigma
    \\
    \phi(0) = \phi^0, \phi^\prime (0) = \phi^1  \hspace{0.5cm} & \text{in} \hspace{0.3cm} \Omega.
    \end{cases}
\end{equation}

As shown in Theorem \ref{theorem:III.3.2} of Chapter \ref{chapter03}, from (\ref{eqIV06}), it follows that the linear wave equation (with $ f \equiv 0 $ in (\ref{eqIV01})), is exactly controllable in $ H^1_0(\Omega) \times L^2(\Omega) $ with the controls in $ L^2( \omega \times (0,T)) $.

We aim to study the exact controllability of (\ref{eqIV01}) in this functional setting of finite energy solutions, under the hypothesis that the estimates (\ref{eqIV06}) and/or (\ref{eqIV07}) are satisfied.

The exact controllability problem for nonlinear equations has been addressed by several authors. L. Markus \cite{markus1965controlability} introduced an Implicit Function Theorem method  for the study of the controllability of systems of ordinary differential equations. Later, this method was adapted and extended to the study of the controllability of nonlinear wave equations (cf. W. C. Chewning \cite{chewning1976controllability}, H. O. Fattorini \cite{fattorini1975local} and D. L. Russell \cite{russell1978control}). This method provides local controllability results, i.e. the initial data in a neighborhood of $ \{ 0,0 \}$ in $ H^1_0(\Omega) \times L^2(\Omega) $ can be driven to $ \{ 0,0 \} ( \text{if } f(0) = 0) $ in time $T$, $T$ being, essentially,  the controllability time of the underlying linear model. It is important to note that these results are of local nature and that this method does not provide global exact controllability results in the sense formulated above. But the method has an advantage of being applicable to a wide class of models and nonlinearities.

Here, with the aim of achieving global results,  we follow a different path, developing and employing fixed point arguments, that yield global exact controllability results for a suitable  class of nonlinearities.
A first fixed point method was introduced in \cite{zuazua1988controlabilite306} and it provided exact controllability for nonlinearities $f$ behaving in an asymptotically linear manner, i. e.
\begin{equation}\label{eqIV09}
    f^\prime \in L^\infty(R),
\end{equation}
and
\begin{equation}\label{eqIV10}
    \exists \lim_{|s| \to \infty} \frac{f(s)}{s}.
\end{equation}


Subsequently, a second fixed point scheme was introduced in \cite{zuazua1991exact} by means of which the controllability in the class of globally Lipschitz nonlinearities was proved, without the assumption (\ref{eqIV10}). The drawback of this second method is that the constructed controls belong to the class $ H^{-\varepsilon} (0,T;L^2(\omega)) $ with $ \varepsilon > 0 $, while the optimal regularity  for the controls of the  finite energy solutions of the linear model is $L^2$. This limitation, as we shall describe below, was overcame in \cite{fu2007exact} combining this fixed point argument and Carleman inequalities. This permits also to cover a slightly larger class of nonlinearities, allowing for a mild logarithmic superlinear growth at infinity.

In the next four sections, the second fixed point method is developed. In Section \ref{sec:IV.6}, for the sake of completeness,  we prove the local controllability of (\ref{eqIV03}) when $f$ is superlinear at infinity. In the last section, dedicated to other related comments, we briefly describe the fist fixed point method applicable in the case where $f$ is asymptotically linear. 
\end{remark}

\section{Description of the Fixed Point Method}
\label{sec:IV.2}
Suppose that $ f \in C^1(\mathbb{R}) $ is globally Lipschitz and define the function $g$ by
\begin{equation}\label{eqIV11}
    g(s) = 
    \begin{cases}
(f(s) - f(0))/s \hspace{1.6cm} & \text{if} \hspace{0.6cm} s \neq 0
    \\
    f^\prime(0) \hspace{1.6cm} & \text{if} \hspace{0.6cm} s = 0.
    \end{cases}
\end{equation}

Given any $ \xi \in L^2(Q) $, we consider the linearized system
\begin{equation}\label{eqIV12}
    \begin{cases}
    y^{\prime\prime} - \Delta y + g(\xi)y= h \chi_\omega - f(0)\hspace{1.6cm} & \text{in} \hspace{0.3cm} Q

    \\
    y = 0\hspace{1.6cm} & \text{on} \hspace{0.3cm} \Sigma

    \\
    y(0) = y^0, y^\prime (0) = y^1 \hspace{1.6cm} & \text{in} \hspace{0.3cm} \Omega.
    \end{cases}
\end{equation}

Since $ f^\prime \in L^\infty(\mathbb{R})$, the potential $g(\xi) \in L^\infty(Q) $ is uniformly bounded by the Lipschitz constant of $f$:
\begin{equation}\label{eqIV13}
\| g(\xi) \|_{L^\infty(Q)} \leq \| f^\prime \|_{L^\infty(\mathbb{R})}, \hspace{0.3cm} \forall \xi \in L^2(Q).
\end{equation}

We fix the initial  and final finite-energy data,  $ \{y^0, y^1\} $ and $ \{z^0, z^1\}  \in H^1_0(\Omega) \times L^2(\Omega) $ respectively.

Suppose that the system (\ref{eqIV12}) is exactly controllable in time $ T>T_0 $. Then, for any $ \xi \in L^2(Q) $, define the HUM control
\begin{equation}\label{eqIV14}
h=h(x, t; \xi) \in L^2(\omega \times (0,T) )
\end{equation}
driving the solution of (\ref{eqIV12}) from $ \{y^0, y^1\} $ to $ \{z^0, z^1\} $. In this way, define the nonlinear map
\begin{equation}\label{eqIV15}
\mathcal{N}(\xi)= y; \hspace{0.3cm} \mathcal{N} : L^2(Q)\to L^2(Q),
\end{equation}
$y=y(x,t;\xi)$ being the controlled solution of (\ref{eqIV12}).

The problem of the exact controllability of (\ref{eqIV01}) is reduced to obtaining a fixed point of the operator $\mathcal{N}$. Indeed, if $ \mathcal{N}(\xi)= \xi = y$, then the solution $y=y(x,t) $ of (\ref{eqIV12}) satisfies (\ref{eqIV01}) and also (\ref{eqIV02}), by construction.

To prove the existence of a fixed point of $\mathcal{N}$, we will apply Schauder's fixed point theorem. Let us remind its statement:
\emph{Let $X$ be a Banach space and $\mathcal{N}:X \to X$ be a continuous and compact operator. Suppose that there exists a convex, closed, bounded and nonempty set  $ C \subset X $, such that}
\begin{equation}\label{eqIV16}
\mathcal{N} (C) \subset C.
\end{equation}

\emph{Then, $\mathcal{N}$ has at least one fixed point in $C$, i.e.}
\begin{equation}\label{eqIV17}
\exists x \in C : \mathcal{N} x = x.
\end{equation}

\begin{remark}
\label{remark:IV.2}
The hypothesis (\ref{eqIV16}) is trivially satisfied if
\begin{equation}\label{eqIV18}
R = \sup_{x \in X} \|\mathcal{N}x\|_X < \infty.
\end{equation}
Indeed, in that case it is enough to choose $ C = \overline{B_R} $, a closed ball with center zero and radius $R$ in $X$.
$ \hfill \square$

It is expected that, for the problem under consideration, (\ref{eqIV18}) will be satisfied with $X = L^2(Q)$. Indeed, (\ref{eqIV13}) ensures that the potential $g(\xi)$ is uniformly bounded in $ L^\infty(Q) $. Therefore, if the control $h$ continuously depends on the potential (in a sense that is specified later), the corresponding family of controls will be uniformly bounded. Finally, for the fixed initial data in (\ref{eqIV12}) and uniformly bounded potentials and controls, solutions will also be uniformly bounded in the energy space
\begin{equation}\label{eqIV19}
C([0,T]; H^1_0(\Omega)) \cap C^1([0,T]; L^2(\Omega))
\end{equation}
and, in particular, in $H^1(Q)$. Then  (\ref{eqIV18}) holds immediately (and therefore (\ref{eqIV16})) but also the compactness of $\mathcal{N}$ in $L^2(Q)$, due to the compactness of the embedding $H^1(Q) \subset L^2(Q) $.

However, as we shall see in Section \ref{sec:IV.4}, significant technical difficulties arise in obtaining the uniform controllability result. Using a variant of the HUM method, working in a different functional framework, it is possible to build controls uniformly bounded in $ H^{-\varepsilon}(0,T;L^2(\omega)) $ with $\varepsilon>0$ arbitrarily small. This allows to conclude the exact controllability of (\ref{eqIV01}) with controls in $ H^{-\varepsilon}(0,T;L^2(\omega)) $.

To achieve this goal, we proceed in several steps. First, in Section \ref{sec:IV.3} we study the exact controllability of (\ref{eqIV12}). Then, the uniform boundedness of the control is proved in Section \ref{sec:IV.4}. Finally, the nonlinear system (\ref{eqIV01}) is studied, concluding the proof.
\end{remark}

\begin{remark}
\label{remark:IV.4.3}
As mentioned above, this result does not provide controls in the sharp class $L^2(\omega \times (0,T)) $, but only in $ H^{-\varepsilon}(0,T;L^2(\omega)) $. The fact that the control can be chosen in $L^2 (\omega \times (0,T) )$ can be established using the techniques in \cite{fu2007exact}, which, actually, allow to treat more general nonlinearities, whose growth at infinity is weaker than $|s| \log^{1/2}|s|.$ The proof in \cite{fu2007exact} uses a fixed point argument, similar as in the proof below, but combined with sharp observability estimates obtained by means of Carleman inequalities. In this way  Theorem \ref{theorem:III.3.4} is improved in three ways.

\begin{itemize}
    \item{The controls can be taken to be in $L^2(0;T;L^2(\omega))$.}
    \item{The support of the control $\omega$ can be a neighborhood of a subset of the boundary of the form $\Gamma(x^0)$.}
    \item{The growth condition on the nonlinearity $f$ can be relaxed so that $ f$ is not needed to be globally Lipschitz, but simply asymptotically smaller than $|s| \log^{1/2} |s|$ as $|s| \rightarrow \infty$.}
\end{itemize}

The proof we present here leads to a slightly weaker result since it employs a purely perturbative analysis, without using ad-hoc Carleman inequalities for the wave equation with potential terms.
\end{remark}

\section{Internal Controllability of the Wave Equation with a Potential}
\label{sec:IV.3}
System (\ref{eqIV12}) is of the following form
\begin{equation}\label{eqIV20}
    \begin{cases}
    y^{\prime\prime} - \Delta y + V(x,t)y = h \chi_\omega - f(0) & \text{in} \hspace{0.3cm} Q
    \\
    y = 0 & \text{in} \hspace{0.3cm} \Sigma
    \\
    y(0) = y^0, y^\prime (0) = y^1 & \text{in} \hspace{0.3cm} \Omega
    \end{cases}
\end{equation}
where $V=V(x,t) (=g(\xi(x,t)))\in L^\infty(Q)$ uniformly bounded, according to  (\ref{eqIV13}).

For a fixed final data $ \{ z^0,z^1 \} \in H^1_0(\Omega) \times L^2(\Omega) $, we solve the system
\begin{equation}\label{eqIV21}
    \begin{cases}
    z^{\prime\prime} - \Delta z + V(x,t)z = - f(0) & \text{in} \hspace{0.3cm} Q
    \\
    z = 0 & \text{in} \hspace{0.3cm} \Sigma
    \\
    z(T) = z^0, z^\prime (T) = z^1 & \text{in} \hspace{0.3cm} \Omega,
    \end{cases}
\end{equation}
with $z \in C([0,T]; H^1_0(\Omega)) \cap C^1([0,T];L^2(\Omega))$ and introduce the new variable
\begin{equation}\label{eqIV22}
p := y - z,
\end{equation}
where the function $p$  satisfies the following system
\begin{equation}\label{eqIV23}
    \begin{cases}
    p^{\prime\prime} - \Delta p + V(x,t)p = h \chi_\omega & \text{in} \hspace{0.3cm} Q
    \\
    p = 0 & \text{in} \hspace{0.3cm} \Sigma
    \\
    p(0) = y^0 - z(0), p^\prime (0) = y^1 - z^\prime(0) & \text{in} \hspace{0.3cm} \Omega.
    \end{cases}
\end{equation}

Proving the exact controllability of (\ref{eqIV20}) in $  H^1_0(\Omega) \times L^2(\Omega) $ with controls in $ L^2(\omega \times (0,T)) $ is equivalent to showing that every solution of (\ref{eqIV23}) can be driven to the equilibrium in time $T$.

Our interest is therefore focused on the study of the following system
\begin{equation}\label{eqIV25}
    \begin{cases}
    y^{\prime\prime} - \Delta y + V(x,t) y = h \chi_\omega & \text{in} \hspace{0.3cm} Q
    \\
    y = 0 & \text{in} \hspace{0.3cm} \Sigma
    \\
    y(0) = y^0, y^\prime (0) = y^1 & \text{in} \hspace{0.3cm} \Omega.
    \end{cases}
\end{equation}
The following result is immediate.
\begin{theorem}
\label{theorem:IV.4.1}
Let $\Omega$ be a bounded domain of $ \mathbb{R}^n $ with the boundary $ \Gamma $ of class $C^2$. Let $\omega \subset \Omega$ be a neighborhood of $\Gamma$. Suppose that $V=V(x,t) \in L^\infty(Q) $ and
\begin{equation*}
T>T_0 = \text{diameter of } \Omega.
\end{equation*}

Then, for each initial data $ \{ y^0,y^1 \} \in H^1_0(\Omega) \times L^2(\Omega) $, there exists a control $  h \in L^2(\omega \times (0,T)) $ such that the solution of (\ref{eqIV25}) satisfies
\begin{equation}\label{eqIV26}
y(T) = y^\prime(T) = 0.
\end{equation}

Moreover, there exists a positive constant $C=C(\|V\|_{L^\infty(Q)}) $ such that the following estimate holds for all $\{ y^0,y^1 \} $:
\begin{equation}\label{eqIV27}
\| h \|_{L^2(\omega \times (0,T))} \leq C \| \{ y^0,y^1 \} \|_{H^1_0(\Omega) \times L^2(\Omega)}.
\end{equation}
\end{theorem}

\begin{proof}
\label{proof:IV.theorem.4.1}
We apply the HUM method, as in the previous chapter. First solve the equation
\begin{equation}\label{eqIV28}
    \begin{cases}
    \phi^{\prime\prime} - \Delta \phi + V(x,t) \phi = 0 & \text{in} \hspace{0.3cm} Q
    \\
    \phi = 0 & \text{in} \hspace{0.3cm} \Sigma
    \\
    \phi(0) = \phi^0, \phi^\prime (0) = \phi^1 & \text{in} \hspace{0.3cm} \Omega
    \end{cases}
\end{equation}
and then solve
\begin{equation}\label{eqIV29}
    \begin{cases}
    y^{\prime\prime} - \Delta y + V(x,t)y = -\phi \chi_\omega & \text{in} \hspace{0.3cm} Q
    \\
    y = 0 & \text{in} \hspace{0.3cm} \Sigma
    \\
    y(T) = y^\prime (T) = 0 & \text{in} \hspace{0.3cm} \Omega.
    \end{cases}
\end{equation}
Next, define the mapping
\begin{equation}\label{eqIV30}
\Lambda \{ \phi^0,\phi^1 \} = \{ y^\prime (0), -y(0) \}.
\end{equation}
Now, the problem reduces down to proving that the map
$
    \Lambda : L^2(\Omega) \times H^{-1}(\Omega) \to L^2(\Omega) \times H^1_0(\Omega) \text{ is an isomorphism.}
$

It is straightforward to check that $\Lambda$ is a continuous linear operator from $L^2(\Omega) \times H^{-1}(\Omega)$ to $L^2(\Omega) \times H^1_0(\Omega)$.

To prove the invertibility,  observe that
\begin{equation}\label{eqIV31}
\langle \Lambda \{ \phi^0,\phi^1 \}, \{ \phi^0,\phi^1 \} \rangle \hspace{0.2cm} = \int_{0}^{T} \int_{\omega} |\phi|^2 \, dxdt.
\end{equation}

Therefore, it is sufficient to prove the estimate
\begin{equation}\label{eqIV32}
\| \phi^0 \|^2_{L^2(\Omega)} + \| \phi^1 \|^2_{H^{-1}(\Omega)} \hspace{0.2cm} \leq C \int_{0}^{T} \int_{\omega} |\phi|^2 \, dxdt,
\end{equation}
which is the main object of the following result.
\end{proof}

\begin{proposition}
\label{propo:IV.4.1}
Suppose that all hypotheses of Theorem \ref{theorem:IV.4.1} are satisfied. Then, there exists a constant $C=C(\|V\|_{L^\infty(Q)}) >0$ such that  (\ref{eqIV32}) holds for every solution of (\ref{eqIV28}).
\end{proposition}

\begin{proof}
\label{proof:IV.propo.4.1}
We know that (\ref{eqIV32}) holds true if $ V \equiv 0 $. Therefore, we proceed by a perturbation argument. Choose $\phi$
\begin{equation}\label{eqIV33}
\phi = \theta + \eta
\end{equation}
such that $\theta$ and $\eta$, respectively, satisfy
\begin{equation}\label{eqIV34}
    \begin{cases}
    \theta^{\prime\prime} - \Delta \theta = 0 & \text{in} \hspace{0.3cm} Q
    \\
    \theta = 0 & \text{in} \hspace{0.3cm} \Sigma
    \\
    \theta(0) = \phi^0,  \theta^\prime (0) = \phi^1 & \text{in} \hspace{0.3cm} \Omega,
    \end{cases}
\end{equation}
and
\begin{equation}\label{eqIV35}
    \begin{cases}
    \eta^{\prime\prime} - \Delta \eta = - V\phi & \text{in} \hspace{0.3cm} Q
    \\
    \eta = 0 & \text{in} \hspace{0.3cm} \Sigma
    \\
    \eta(0) = \eta^\prime (0) = 0 & \text{in} \hspace{0.3cm} \Omega.
    \end{cases}
\end{equation}

As $T>T_0 = $ diameter of $\Omega$, 
\begin{equation}\label{eqIV36}
\| \phi^0 \|^2_{L^2(\Omega)} + \| \phi^1 \|^2_{H^{-1}(\Omega)} \hspace{0.2cm} \leq C \int_{0}^{T} \int_{\omega} |\theta|^2 \, dxdt
 \leq C \int_{0}^{T} \int_{\omega} \left[ |\phi|^2 + |\eta|^2 \right] \, dxdt.
\end{equation}

Therefore, it is enough to prove the estimate
\begin{equation}\label{eqIV37}
\| \eta \|^2_{L^2(Q)} \leq C \int_{0}^{T} \int_{\omega} |\phi|^2 \, dxdt.
\end{equation}

Now, we go by the method of contradiction. If (\ref{eqIV37}) does not hold, there exists a sequence of solutions $ \{ \phi_n \} $ of (\ref{eqIV28}) corresponding to the initial data $ \{ \phi^0_n, \phi^1_n \} $ such that
\begin{equation}\label{eqIV38}
\int_{0}^{T} \int_{\omega} |\phi_n|^2 \, dxdt \to 0
\end{equation}
while the associated sequence $ \{ \eta_n \} $ of solutions of (\ref{eqIV35}) satisfies
\begin{equation}\label{eqIV39}
\| \eta_n \|_{L^2(Q)} = 1.
\end{equation}

Now, combining (\ref{eqIV36}), (\ref{eqIV38}) and (\ref{eqIV39}), it can be deduced that 
$ \{ \phi^0_n, \phi^1_n \} $ is bounded in $ L^2(\Omega) \times H^{-1}(\Omega)$, and therefore,
$
\{ \phi_n \} \text{ is bounded in } L^2(Q).
$

Since the classical estimates for the wave equation (\ref{eqIV35}) ensure that
$
\{ \eta_n \}$ is bounded in $H^1(Q),
$
we obtain that
$
\{ \eta_n \} \text{ is relatively compact in } L^2(Q).
$

Next, extracting a suitable subsequence and passing to the limit we obtain 
$
\phi_n \rightharpoonup \phi \text{ weakly in } L^2(Q),
$
where $ \phi = \phi(x,t) $ is a solution of (\ref{eqIV28}) such that
\begin{equation}\label{eqIV40}
\phi = 0 \text{ in } \omega \times (0,T).
\end{equation}
On the other hand, $ \eta_n \rightarrow \eta $ in $L^2(Q)$ strongly, where $ \eta = \eta(x,t) $ is a solution of (\ref{eqIV35}) with
\begin{equation}\label{eqIV41}
\|\eta\|_{L^2(Q)} = 1.
\end{equation}

By the unique continuation result of \cite{fu2007exact} the only solution of (\ref{eqIV28}) satisfying (\ref{eqIV40}) is the trivial one. But if $ \phi \equiv 0 $, by the uniqueness of the solutions of (\ref{eqIV35}), it follows that $ \eta \equiv 0 $ and this contradicts  (\ref{eqIV41}).

This way, the proofs of Proposition \ref{propo:IV.4.1} and Theorem \ref{theorem:IV.4.1} are complete.
\end{proof}

In the study of the exact controllability of the semilinear wave equation  a different functional framework is needed. The following result holds for this purpose.

\begin{theorem}
\label{theorem:IV.4.2}
Suppose that the hypotheses of Theorem \ref{theorem:IV.4.1} are satisfied. Let $ \varepsilon \in (0,1/2) $.

Then, for each initial data $ \{y^0, y^1\} \in H^{1-\varepsilon}_0(\Omega) \times H^{-\varepsilon}(\Omega) $, there exists a control $ h \in H^{-\varepsilon}(0,T;L^2(\omega))$ such that the solution $y = y(x,t)$ of (\ref{eqIV25}) satisfies (\ref{eqIV26}).

Moreover, the following estimate holds
\begin{equation}\label{eqIV42}
\|h\|_{H^{-\varepsilon}(0,T; L^2(\omega))} \leq C \| \{y^0, y^1 \} \|_{H^{1-\varepsilon}_0 (\Omega) \times H^{-\varepsilon}(\Omega)}
\end{equation}
with $C=C(\|V\|_{L^\infty(Q)}, \varepsilon)$.
\end{theorem}

\begin{proof}
\label{proof:IV.theorem.4.2}
Consider\begin{equation}\label{eqIV43}
    \begin{cases}
    y^{\prime\prime} - \Delta y + V(x,t)y = -I_\varepsilon \phi \chi_\omega & \text{in} \hspace{0.3cm} Q
    \\
    y = 0 & \text{in} \hspace{0.3cm} \Sigma
    \\
    y(T) = y^\prime (T) = 0 & \text{in} \hspace{0.3cm} \Omega
    \end{cases}
\end{equation}
where $I_\varepsilon : H^\varepsilon(0,T; L^2(\omega)) \to  H^{-\varepsilon}(0,T; L^2(\omega)) $ is the canonical duality isomorphism and define the operator 
$
    \Lambda_\varepsilon \{ \phi^0, \phi^1 \} = \{ y^\prime(0), -y(0) \}.
$

The proof now boils down to show that 
\begin{equation*}
    \Lambda_\varepsilon : H^\varepsilon(\Omega) \times H^{-1+\varepsilon}(\Omega) \to H^{-\varepsilon}(\Omega) \times H^{1-\varepsilon}_0 (\Omega)  \text{ is an isomorphism.}
\end{equation*}

Indeed, $ \Lambda_\varepsilon $ is a linear continuous operator. On the other hand, we have
\begin{equation}\label{eqIV44}
\langle \Lambda_\varepsilon \{ \phi^0,\phi^1 \}, \{ \phi^0,\phi^1 \} \rangle = \| \phi \|^2_{H^\varepsilon(0,T; L^2(\omega))}.
\end{equation}

It is therefore enough to prove the following estimate
\begin{equation}\label{eqIV45}
\| \{ \phi^0, \phi^1 \} \|^2_{H^\varepsilon(\omega) \times H^{\varepsilon-1}(\Omega)} \leq C \| \phi \|^2_{H^\varepsilon(0,T;  L^2(\omega))}
\end{equation}
which is the theme of the following result.
\end{proof}

\begin{proposition}
\label{propo:IV.4.2}
Under the conditions of Theorem \ref{theorem:IV.4.1} the estimate (\ref{eqIV45}) holds for all $ \varepsilon \in (0,1], \varepsilon \neq 1/2 $.
\end{proposition}

\begin{proof}
\label{proof:IV.propo.4.2}
The proof splits into two steps. The case $ V \equiv 0 $ is considered in the first case. The second one handles the general situation.

\emph{Step 1. The case of  $ V \equiv 0.$}
As seen in Section \ref{sec:IV.1}, as $ T >$ diameter of $ \Omega $, (\ref{eqIV06}) and (\ref{eqIV07}) hold and therefore (\ref{eqIV45}) holds as well for $ \varepsilon = 0,1 $. The general case (\ref{eqIV45}) with $ \varepsilon \in (0,1) $  is derived by  interpolation.

\emph{Step 2. The general case $ V \in L^\infty(Q). $ }
Decompose the solution of (\ref{eqIV28}) as in (\ref{eqIV33}) and apply the estimate (\ref{eqIV45}), with $ V \equiv 0 $, to obtain
\begin{align}\label{eqIV46}
    \begin{split}
\| \{ \phi^0, \phi^1 \} \|^2_{H^\varepsilon(\omega) \times H^{-1+\varepsilon}(\Omega)} \leq C \left[ \| \phi \|^2_{H^\varepsilon(0,T; L^2(\omega))} + \| \eta \|^2_{H^\varepsilon(0,T; L^2(\omega))} \right].
    \end{split}
\end{align}

However, given that $V$ is uniformly bounded in $L^\infty(Q)$, it holds
\begin{align}\label{eqIV47}
    \begin{split}
    \| \eta \|_{H^\varepsilon(0,T; L^2(\omega))} \leq 
    C \| \eta \|_{W^{1,\infty}(0,T;L^2(\Omega))}
    \leq C \| V\phi \|_{L^2(Q)} \leq C \| \phi \|_{L^2(Q)}.
    \end{split}
\end{align}

Now by (\ref{eqIV46}) and (\ref{eqIV47}), 
\begin{align}\label{eqIV48}
    \begin{split}
\| \{ \phi^0, \phi^1 \} \|^2_{H^\varepsilon(\Omega) \times H^{-1+\varepsilon}(\Omega)}
\leq C \left[ \| \phi \|^2_{H^\varepsilon(0,T; L^2(\omega))} + \| \phi \|^2_{L^2(Q)} \right].
    \end{split}
\end{align}

We need to show that
\begin{equation}\label{eqIV49}
\| \phi \|_{L^2(Q)} \leq C \| \phi \|_{H^\varepsilon(0,T; L^2(\omega))}.
\end{equation}

The rest of the proof goes by the contradiction argument. If (\ref{eqIV49}) is not fulfilled, there exists a sequence $ \{ \phi_k \} $ of solutions of (\ref{eqIV28}) such that
\begin{equation}\label{eqIV50}
\| \phi_k \|_{H^\varepsilon(0,T; L^2(\omega))} \to 0,
\end{equation}
and
\begin{equation}\label{eqIV51}
\| \phi_k \|_{L^2(Q)} = 1.
\end{equation}
From (\ref{eqIV48}), (\ref{eqIV50}) and (\ref{eqIV51}) it follows that
\begin{equation}\label{eqIV52}
\{ \phi^0_k, \phi^1_k \} \text{ is bounded in } H^\varepsilon(\Omega) \times H^{-1+\varepsilon}(\Omega).
\end{equation}

By (\ref{eqIV52}) and by means of classical estimates for the wave equation (\ref{eqIV28}), we obtain that
\begin{equation}\label{eqIV53}
\{ \phi_k \} \text{ is bounded in } H^\varepsilon(Q)
\end{equation}
and therefore, as $ \varepsilon > 0, $

\begin{equation}\label{eqIV54}
\{ \phi_k \} \text{ is relatively compact in } L^2(Q).
\end{equation}

Now, extracting  subsequences, we have
\begin{equation*}
\begin{split}
\phi_k \rightharpoonup \phi \hspace{0.1cm} \text{ weakly in } H^\varepsilon(Q); \hspace{0.3cm} \phi_k \rightarrow \phi \hspace{0.1cm} \text{ strongly in } L^2(Q)
\end{split}
\end{equation*}
where $ \phi = \phi (x,t) $ is a solution of (\ref{eqIV28}) and
\begin{equation}\label{eqIV55}
 \phi = 0 \text{ in } \omega \times (0,T)
\end{equation}
\begin{equation}\label{eqIV56}
 \| \phi \|_{L^2(Q)} = 1.
\end{equation}

The unique continuation result in \cite{fu2007exact}  ensures that the only solution of (\ref{eqIV28}) satisfying (\ref{eqIV55}) is $ \phi \equiv 0, $ which contradicts (\ref{eqIV56}).

The proof of Proposition \ref{propo:IV.4.2} is completed as well as that of Theorem \ref{theorem:IV.4.2}.
\end{proof}

\begin{remark}
\label{remark:IV.4.4}
Note that the estimate (\ref{eqIV45}) has not been shown for $ \varepsilon =0 $ due to the lack of compactness.
$ \hfill \square $
\end{remark}
\section{Uniform Controllability}
\label{sec:IV.4}
We start by proving that estimate  (\ref{eqIV45}) holds uniformly provided the potentials are uniformly bounded in $ L^\infty(Q) $.

\begin{proposition}
\label{propo:IV.4.3}
Suppose that the hypotheses of Theorem \ref{theorem:IV.4.1} are satisfied and $ \varepsilon \in (0,1], \varepsilon \neq 1/2 $. Then, for each $ M>0 $ there is a constant $ C=C_\varepsilon(M) $ such that the inequality (\ref{eqIV45}) holds for every solution of (\ref{eqIV28}) with potential $ V \in L^\infty(Q) $ such that $ \| V \|_{L^\infty(Q)} \leq M $.
\end{proposition}

\begin{proof}
\label{proof:IV.propo.4.3}
We proceed as in the proof of Proposition \ref{propo:IV.4.2}. The estimate (\ref{eqIV48}) is uniform with respect to potentials $V$ such that $ \| V \|_{L^\infty(Q)} \leq M $. Therefore it is enough to prove the existence of $ C=C_\varepsilon(M) $ such that
\begin{equation}\label{eqIV57}
\| \phi \|_{L^2(Q)} \leq C \| \phi \|_{H^\varepsilon(0,T; L^2(\omega))}
\end{equation}
for each solution of (\ref{eqIV28}) with potential $V$ such that $ \| V \|_{L^\infty(Q)} \leq M $.

We argue by contradiction. If (\ref{eqIV57}) is not fulfilled, there is a sequence of potentials $ \{ V_n \} $, uniformly bounded in $ L^\infty(Q) $, and a sequence $ \{ \phi_n \} $ of solutions of the corresponding problems (\ref{eqIV28}) satisfying (\ref{eqIV50}) and (\ref{eqIV51}).

From the uniformity of the estimate (\ref{eqIV48}) with respect to $ \{ V_n \} $, (\ref{eqIV50}) and (\ref{eqIV51}) it follows that
$
\{ \phi^0_n,\phi^1_n \} \text{ is bounded in } H^\varepsilon(\Omega) \times H^{\varepsilon-1}(\Omega)
$
and therefore we have (\ref{eqIV53}).

Extracting sub-sequences, 
$
V_n \rightharpoonup V $ weakly* in $L^\infty(Q)$,
$\phi_n \rightharpoonup \phi$ weakly in $H^\varepsilon(Q)$, and
$
\phi_n \rightarrow \phi$ strongly in $L^2(Q).
$

Passing to the limit in the equation verified by $V_n$ and $\phi_n$ it follows that $ \phi = \phi(x,t) $ is a solution (\ref{eqIV28}) for the limit potential $V$ and that it satisfies (\ref{eqIV55}) and (\ref{eqIV56}).

Thanks to the unique continuation principle in \cite{fu2007exact} we reach again a contradiction. This concludes the proof of Proposition 4.3.
\end{proof}

As an immediate consequence of the uniformity of the estimate (\ref{eqIV45}) the following uniform exact controllability result is obtained.

\begin{theorem}
\label{theorem:IV.4.3}
Suppose that the hypotheses of Theorem \ref{theorem:IV.4.1} are satisfied. Let $ \varepsilon \in (0, 1/2) $.

Then, for all $M>0$ there is a constant $ C = C_\varepsilon(M) $ such that the control $ h \in H^{-\varepsilon}(0,T;L^2(\omega)) $ constructed in Theorem  \ref{theorem:IV.4.2} satisfies
\begin{equation}\label{eqIV58}
\| h \|_{H^{-\varepsilon}(0,T; L^2(\omega))} \leq C \| \{ y^0, y^1 \} \|_{H^{1-\varepsilon}_0(\Omega) \times H^{-\varepsilon}(\Omega) }
\end{equation}
for all potentials $ V \in L^\infty(Q) $ such that
$
{\| V \|_{L^\infty(Q)}} \leq M.
$
\end{theorem}

\begin{proof}
\label{proof:IV.theorem.4.3}
It is an immediate consequence of the construction of the control carried out in the proof of Theorem \ref{theorem:IV.4.2} and of the uniform estimate of Proposition \ref{propo:IV.4.3}.
\end{proof}

\section{The Main Result}
\label{sec:IV.5}
The following exact controllability result for the semilinear equation (\ref{eqIV01}) holds:

\begin{theorem}
\label{theorem:IV.4.4}
Let $ \Omega $ be a bounded domain of $ \mathbb{R}^n $ with boundary $ \Gamma $ of class $ C^2 $. 

Let $x^0 \in \mathbb{R}^n $ and $\omega$ be a neighborhood of $ \overline{\Gamma(x^0)} \text{ in } \Omega $ and $ T>2 \| x - x^0 \|_{L^\infty(\Omega)} $. 

Suppose that $f$ is globally Lipschitz and $ \varepsilon \in (0, 1/2) $.

Then, for each pair of initial and final data $ \{ y^0,y^1 \}, \{ z^0,z^1 \} \in H^{1-\varepsilon}_0(\Omega) \times H^{-\varepsilon}(\Omega) $ there is a control $ h \in  H^{-\varepsilon}(0,T; L^2(\omega)) $ such that the solution $  y = y(x,t) $ of (\ref{eqIV01}) satisfies (\ref{eqIV02}).
\end{theorem}

\begin{remark}
\label{remark:IV.4.5}
Theorem \ref{theorem:IV.4.4} shows, in particular, that, if $ \{ y^0,y^1 \}, \{ z^0,z^1 \} \in H^{1}_0(\Omega) \times L^2(\Omega) $, there is a control $ h_\varepsilon \in H^{-\varepsilon}(0,T; L^2(\omega)) $ with $ \varepsilon > 0 $ arbitrarily small. However, the existence of a control in $ L^2(\omega \times (0,T)) $ is not established.
$ \hfill \square $
\end{remark}

\begin{proof}
\label{proof:IV.remark.4.5}
Suppose first that $ f \in C^1(\mathbb{R}) $ is globally Lipschitz.

We fix initial and final data
$ \{ y^0,y^1 \}, \{ z^0,z^1 \} \in H^{1-\varepsilon}_0(\Omega) \times H^{-\varepsilon}(\Omega) $. We define a nonlinear operator $ \mathcal{N} : L^2(Q) \to L^2(Q) $ such that for each $ \xi \in L^2(Q), \mathcal{N}\xi = y $ where $ y=y(x,t) $ is the solution of
\begin{equation}\label{eqIV59}
    \begin{cases}
    y^{\prime\prime} - \Delta y + g(\xi)y = h \chi_\omega - f(0) & \text{in} \hspace{0.3cm} Q
    \\
    y = 0 & \text{in} \hspace{0.3cm} \Sigma
    \\
    y(0) = y^0, y^\prime (0) = y^1 & \text{in} \hspace{0.3cm} \Omega
    \end{cases}
\end{equation}
which satisfies (\ref{eqIV02}) with $g$ as in (\ref{eqIV11}) and $ h \in H^{-\varepsilon}(0,T; L^2(\omega)) $ is the control obtained by the HUM method as above.

As mentioned in Section \ref{sec:IV.3}, the control $h$ is the one taking the solution of
\begin{equation}\label{eqIV60}
    \begin{cases}
    p^{\prime\prime} - \Delta p + g(\xi)p = h \chi_\omega & \text{in} \hspace{0.3cm} Q
    \\
    p = 0 & \text{in} \hspace{0.3cm} \Sigma
    \\
    p(0) = y^0 - z(0), p^\prime (0) = y^1 - z^\prime(0) & \text{in} \hspace{0.3cm} \Omega
    \end{cases}
\end{equation}
to equilibrium at time $T$, where $z=z(x,t)$ is the solution of
\begin{equation}\label{eqIV61}
    \begin{cases}
    z^{\prime\prime} - \Delta z + g(\xi)z = -f(0) & \text{in} \hspace{0.3cm} Q
    \\
    z = 0 & \text{in} \hspace{0.3cm} \Sigma
    \\
    z(T) = z^0, z^\prime (T) = z^1 & \text{in} \hspace{0.3cm} \Omega
    \end{cases}
\end{equation}
and $y = p+z$.

As $ \| g(\xi) \|_{L^\infty(Q)} \leq \| f^\prime \|_{L^\infty(\mathbb{R})} $  for all $ \xi \in L^2(Q) $, classical estimates of the wave equation ensure that
\begin{equation}\label{eqIV62}
\begin{split}
\| \{ z(0), z^\prime(0) \} \|_{H^{1-\varepsilon}_0(\Omega) \times H^{-\varepsilon}(\Omega)}
\leq C \{ \| \{ z^0,z^1 \} \|_{H^{1-\varepsilon}_0(\Omega) \times H^{-\varepsilon}(\Omega)} + |f(0)| \}.
\end{split}
\end{equation}

On the other hand, thanks to the uniform estimate (\ref{eqIV58}) we deduce that
\begin{equation}\label{eqIV63}
\| h \|_{H^{-\varepsilon}(0,T;L^2(\Omega))} \leq
C \| \{ y^0 - z(0), y^1 - z^\prime(0) \} \|_{H^{1-\varepsilon}(\Omega) \times H^{-\varepsilon}(\Omega)}.
\end{equation}

Combining (\ref{eqIV62}) and (\ref{eqIV63}) it follows that
\begin{equation}\label{eqIV64}
\| h \|_{H^{-\varepsilon}(0,T;L^2(\omega))} \leq C \hspace{0.6cm} \forall \xi \in L^2(Q)
\end{equation}
and therefore
\begin{equation}\label{eqIV65}
\| y \|_{L^\infty(0,T;H^{1-\varepsilon}(\Omega))} + \| y^\prime \|_{L^\infty(0,T;H^{-\varepsilon}(\Omega))} \leq C \hspace{0.4cm} \forall \xi \in L^2(Q),
\end{equation}
with a constant $C>0$ independent of $\xi$.

In this way we see that the image of $\mathcal{N}$ is contained in a bounded set of $ L^\infty(0,T;H^{1-\varepsilon}(\Omega)) \cap W^{1,\infty}(0,T;H^{-\varepsilon}(\Omega)) $. On the other hand, from the continuity of $g$, it follows that $ \mathcal{N}:L^2(Q) \to L^2(Q) $ is continuous. As the inclusion of $ L^\infty(0,T;H^{1-\varepsilon}(\Omega)) \cap W^{1,\infty}(0,T;H^{-\varepsilon}(\Omega)) $ in $L^2(Q)$ is compact we conclude that
$
\mathcal{N}:L^2(Q) \to L^2(Q)$  is continuous and compact
and also
\begin{equation*}
\sup_{\xi \in L^2(Q)} \hspace{0.3cm} \| \mathcal{N} ( \xi ) \|_{L^2(Q)} < \infty.
\end{equation*}

Applying Schauder's Theorem (see Remark \ref{theorem:IV.4.2} above) it follows that $ \mathcal{N} $ admits a fixed point $ \xi = y $ satisfying (\ref{eqIV01}) and (\ref{eqIV02}). The control $h$ we are looking for is therefore the one corresponding to the fixed point $\xi$.

In the general case in which $f$ is globally Lipschitz but not of class $C^1$, we proceed by an approximation argument.

We introduce a regularizing sequence $ \{ f_n \} \subset C^1( \mathbb{R}) $ such that
\begin{equation}\label{eqIV66}
 f_n \rightarrow f \text{ in } L^\infty(\mathbb{R})
\end{equation}
and
\begin{equation}\label{eqIV67}
 \| f^\prime_n \|_{L^\infty(\mathbb{R})} \leq \| f^\prime \|_{ L^\infty(\mathbb{R})}
 \hspace{0.3cm} \forall n \in \mathbb{N}.
\end{equation}

Equation (\ref{eqIV01}) is exactly controllable for each nonlinearity $f_n$. Furthermore, the estimates (for fixed initial and final data) depend on the nonlinearity $f$ solely through its value at 0 and its Lipschitz constant. From (\ref{eqIV64}) and (\ref{eqIV65}) it follows that the sequence of controls $ h_n $ and the sequence $ y_n $ of solutions obtained will be uniformly bounded in their respective spaces $ H^{-\varepsilon}(0,T;L^2(\omega)) $ and $ L^\infty(0,T;H^{1-\varepsilon}_0(\Omega)) \cap W^{1,\infty}(0,T;H^{-\varepsilon}(\Omega)) $

Passing to the limit, we obtain a control $h$ and a solution $y$ of (\ref{eqIV01}) satisfying (\ref{eqIV02}).

This concludes the proof of Theorem \ref{theorem:IV.4.4}.
\end{proof}

\section{Local Controllability in the Superlinear Case}
\label{sec:IV.6}
Suppose that $ f \in W^{1,\infty}_{loc}(\mathbb{R}) $ is such that
\begin{equation}\label{eqIV68}
 f(0) = 0,
\end{equation}
and
\begin{equation}\label{eqIV69}
 \exists C >0, p>1 : |f^\prime(s)| \leq C|s|^{p-1} \hspace{0.4cm} \forall s \in \mathbb{R} \hspace{0.2cm} \text{with } (n-2)p < n.
\end{equation}
The following local controllability result holds:

\begin{theorem}
\label{theorem:IV.4.5}
Let $ \Omega $ be a bounded domain of $ \mathbb{R}^n $ with boundary $ \Gamma $ of class $C^2$. Let $ x^0 \in \mathbb{R}^n, \omega \subset \Omega $ be a neighborhood of $ \overline{\Gamma(x^0)} $ and $ T>2 \| x-x^0 \|_{L^\infty(\Omega)} $. Suppose that $f$ satisfies (\ref{eqIV68}) and (\ref{eqIV69}).

Then, there is $\delta>0$ such that if the initial data $ \{ y^0,y^1 \} \in H^1_0(\Omega) \times L^2(\Omega) $ verifies
\begin{equation}\label{eqIV70}
\| \{ y^0,y^1 \} \|_{H^1_0(\Omega) \times L^2(\Omega)} < \delta
\end{equation}
there exists a control $ h \in L^2(\omega \times (0,T)) $ such that the solution $ y= y(x,t) $ of (\ref{eqIV01}) satisfies
\begin{equation}\label{eqIV71}
y(T) = y^\prime(T) = 0.
\end{equation}
\end{theorem}

\begin{proof}
\label{proof:IV.theorem.4.5}
We use a nonlinear version of the HUM method.
First we solve the problem (\ref{eqIV08}) with data $ \{ \phi^0,\phi^1 \} \in L^2(\Omega) \times H^{-1}(\Omega) $ and then
\begin{equation}\label{eqIV72}
    \begin{cases}
    y^{\prime\prime} - \Delta y + f(y) = -\phi \chi_\omega & \text{in} \hspace{0.3cm} Q
    \\
    y = 0 & \text{in} \hspace{0.3cm} \Sigma
    \\
    y(T) = y^\prime (T) = 0 & \text{in} \hspace{0.3cm} \Omega.
    \end{cases}
\end{equation}

Thanks to the growth condition (\ref{eqIV69}), (\ref{eqIV72}) admits a unique local finite energy solution, i.e. for every initial data pair $ \{ \phi^0,\phi^1 \} \in L^2(\Omega) \times H^{-1}(\Omega) $ there exists
$
    \tau = \tau ( \| \{ \phi^0,\phi^1 \} \|_{L^2(\Omega) \times H^{-1}(\Omega)} ) > 0
$
such that (\ref{eqIV72}) admits a unique solution in
\begin{equation}\label{eqIV73}
    y \in C ( [T - \tau,T]; H^1_0(\Omega)) \cap C^1( [T - \tau,T]; L^2(\Omega) ).
\end{equation}

As we shall see later, there is $ \alpha >0 $ such that if
\begin{equation}\label{eqIV74}
\| \{ \phi^0,\phi^1 \} \|_{L^2(\Omega) \times H^{-1}(\Omega)} < \alpha
\end{equation}
then $ \tau > T $ and, therefore, $y$ is globally defined in $[0,T]$:
\begin{equation}\label{eqIV75}
    y \in C ( [0,T]; H^1_0(\Omega)) \cap C^1( [0,T]; L^2(\Omega) ).
\end{equation}

We define the nonlinear operator
\begin{equation}\label{eqIV76}
     \mu: B_\alpha \subset L^2(\Omega) \times H^{-1}(\Omega) \to L^2(\Omega) \times H^1_0(\Omega)
\end{equation}
through
\begin{equation}\label{eqIV77}
     \mu \{ \phi^0,\phi^1 \} = \{ y^\prime(0), -y(0) \},
\end{equation}
$ B_\alpha $ being the ball of radius $ \alpha $ and center $ \{ 0,0 \} $ in $ L^2(\Omega) \times H^{-1}(\Omega) $.

The proof of Theorem \ref{theorem:IV.4.5} boils down to proving that  (\ref{eqIV70}) holds, if and only if the equation
\begin{equation}\label{eqIV78}
     \mu \{ \phi^0,\phi^1 \} = \{ y^1, -y^0 \}
\end{equation}
admits a solution $ \{ \phi^0,\phi^1 \} \in B_\alpha $.

To solve (\ref{eqIV78}) we decompose the $ \mu $ operator as follows
\begin{equation}\label{eqIV79}
     \mu = \Lambda + \Theta
\end{equation}
where
\begin{equation}\label{eqIV80}
     \Lambda \{ \phi^0,\phi^1 \} = \{ y^\prime_0(0), -y_0(0) \},
\end{equation}
and
\begin{equation}\label{eqIV81}
     \Theta \{ \phi^0,\phi^1 \} = \{ \eta^\prime(0), -\eta(0) \}
\end{equation}
where $y_0 = y_0(x,t)$ and $ \eta= \eta(x,t) $ are solutions of
\begin{equation}\label{eqIV82}
    \begin{cases}
    y^{\prime\prime}_0 - \Delta y_0 = -\phi \chi_\omega & \text{in} \hspace{0.3cm} Q
    \\
    y_0 = 0 & \text{in} \hspace{0.3cm} \Sigma
    \\
    y_0(T) = y^\prime_0 (T) = 0 & \text{in} \hspace{0.3cm} \Omega
    \end{cases}
\end{equation}
and
\begin{equation}\label{eqIV83}
    \begin{cases}
    \eta^{\prime\prime} - \Delta \eta = -f(y) & \text{in} \hspace{0.3cm} Q
    \\
    \eta = 0 & \text{in} \hspace{0.3cm} \Sigma
    \\
    \eta(T) = \eta^\prime (T) = 0 & \text{in} \hspace{0.3cm} \Omega,
    \end{cases}
\end{equation}
respectively.
Obviously $ y = y_0 + \eta $.

From Theorem \ref{theorem:III.3.1} of Chapter \ref{chapter03} we know that
\begin{equation}\label{eqIV84}
     \Lambda: L^2(\Omega) \times H^{-1}(\Omega) \to L^2(\Omega) \times H^1_0(\Omega) \text{ is an isomorphism}.
\end{equation}

Equation (\ref{eqIV78}) can be rewritten as
\begin{equation*}
     \Lambda \{ \phi^0,\phi^1 \} + \Theta \{ \phi^0,\phi^1 \}  = \{ y^1, -y^0 \}
\end{equation*}
or equivalently
\begin{equation}\label{eqIV85}
     \{ \phi^0,\phi^1 \} = \Lambda^{-1} \{ y^1, -y^0 \} - \Lambda^{-1} \Theta \{ \phi^0,\phi^1 \}.
\end{equation}

Therefore, the problem is reduced to proving the existence of a fixed point of the operator
\begin{equation}\label{eqIV86}
\begin{split}
    \mathcal{K}: L^2(\Omega) \times H^{-1}(\Omega) \to L^2(\Omega) \times H^{-1}(\Omega) ; \\ 
    \mathcal{K}\{ \phi^0,\phi^1 \} = \Lambda^{-1} \{ y^1, -y^0 \} - \Lambda^{-1} \Theta \{ \phi^0,\phi^1 \}.
\end{split}
\end{equation}

For this purpose we apply Schauder's fixed point Theorem. Let us check that $\mathcal{K}$ verifies the required hypotheses.

\emph{Step 1. Compactness of $\mathcal{K}$.}
\label{sec:IV.6-a}
It is enough to prove that
\begin{equation}\label{eqIV87}
     \Theta: B_\alpha \subset L^2(\Omega) \times H^{-1}(\Omega) \to L^2(\Omega) \times H^1_0(\Omega)
\end{equation}
is a continuous and compact operator.

Multiplying equation (\ref{eqIV72}) by $y^\prime$ and integrating in $\Omega$ we obtain
\begin{equation}\label{eqIV88}
     \frac{dE(t)}{dt} = - \int_\omega \phi (x,t)y^\prime (x,t)dx - \int_\Omega f(y (x,t)) y^\prime (x,t)dx
\end{equation}
with
\begin{equation}\label{eqIV89}
     E(t) = \frac{1}{2} \int_\Omega \left[ |y^\prime(x,t)|^2 +  | \nabla y(x,t)|^2 \right] dx.
\end{equation}

One has
\begin{equation}\label{eqIV90}
     E(T) = 0
\end{equation}
and
\begin{equation}\label{eqIV91}
    \begin{split}
   \left| \int_\omega \phi (x,t)y^\prime (x,t)dx \right| & \leq \frac{1}{2} \int_\Omega | y^\prime(x,t) |^2 dx + 
   \frac{1}{2} \int_\Omega | \phi (x,t) |^2 dx
   \\ & \leq E(t) + C \| \{ \phi^0, \phi^1 \} \|^2_{L^2(\Omega) \times H^{-1}(\Omega)}        
    \end{split}
\end{equation}
and
\begin{equation}\label{eqIV92}
   \left| \int_\Omega f(y(x,t))y^\prime (x,t)dx \right| \leq 
   E(t) + C \int_\Omega | f(y(x,t)) |^2 dx.
\end{equation}

Combining (\ref{eqIV88}), (\ref{eqIV91}) and (\ref{eqIV92}) it follows
\begin{equation}\label{eqIV93}
   \left| \frac{dE(t)}{dt} \right| \leq 2E(t) + C \| \{ \phi^0, \phi^1 \} \|^2_{L^2(\Omega) \times H^{-1}(\Omega)} + C \| f(y) \|^2_{L^2(\Omega)}.
\end{equation}

From (\ref{eqIV68}) - (\ref{eqIV69}) it follows that
\begin{equation}\label{eqIV94}
   | f(s) | \leq C |s|^p \hspace{1cm} \forall s \in \mathbb{R}
\end{equation}

and therefore, as $ p(n-2) < n, $
\begin{equation}\label{eqIV95}
    \| f(y) \|_{L^2(\Omega)} \leq C \| y \|^p_{L^{2p}(\Omega)} \leq C \| y \|^p_{H^1_0(\Omega)} \leq C [E(t)]^{p/2}.
\end{equation}

Consequently we obtain the differential inequality
\begin{equation}\label{eqIV96}
   \frac{dE(t)}{dt} \leq C \left ( E(t) + [ E(t) ]^{p/2} + \| \{ \phi^0, \phi^1 \} \|^2_{L^2(\Omega) \times H^{-1}(\Omega)} \right )
\end{equation}
leading to
\begin{equation}\label{eqIV97}
   E(t) \leq C \| \{ \phi^0, \phi^1 \} \|^2_{L^2(\Omega) \times H^{-1}(\Omega)} (T - t)e^{C(T-t)} e^{C\int^T_t [E(s)]^{{p/2}-1}ds}.
\end{equation}

From (\ref{eqIV97}) it can be concluded that if $ \alpha >0 $ is small enough and $ \{ \phi^0, \phi^1 \} \in B_\alpha $, then $ E(\cdot) \in L^\infty(0,T) $ and we have the estimate
\begin{equation}\label{eqIV98}
   E(t) \leq C \| \{ \phi^0, \phi^1 \} \|^2_{L^2(\Omega) \times H^{-1}(\Omega)} \hspace{0.6cm} \forall t \in (0,T)
\end{equation}
for some constant $ C = C (\alpha) > 0 $.
From (\ref{eqIV98}) it follows that if $ \{ \phi^0, \phi^1 \} \in B_\alpha $, 
\begin{equation}\label{eqIV99}
y(t) \, \text{ is uniformly bounded in } L^\infty(0,T;H^1_0(\Omega))
\end{equation}
and, as $p(n-2) <n$,
\begin{equation}\label{eqIV100}
\begin{split}
f(y(t)) \text{ is uniformly bounded in } L^\infty(0,T;H^{\varepsilon}(\Omega))
\end{split}
\end{equation}
for some $ \varepsilon > 0$ that only depends on $p$ (cf. J. Simon \cite{simon1979regularite}).

Combining (\ref{eqIV99}) and classical regularity results for the wave equation it follows that
\begin{equation}\label{eqIV101}
\eta \text{ is uniformly bounded in }
\\ L^\infty(0,T;H^{1+\varepsilon}(\Omega)) \cap W^{1,\infty}(0,T;H^\varepsilon(\Omega)).
\end{equation}

In particular $ \Theta ( B_\alpha ) $ is a bounded set of $ H^\varepsilon(\Omega) \times H^{1+\varepsilon}(\Omega) $ and therefore $ \Theta ( B_\alpha ) $  is relatively compact in $ L^2(\Omega) \times H^1_0(\Omega) $.

The compactness of $\mathcal{K}$ is therefore proved.

\label{sec:IV.6-b}
\emph{Step 2. $ \mathcal{K} (B_\beta)\subset B_\beta$.}
Let us see now that there exists $ \beta \in (0,\alpha]  $ such that
\begin{equation}\label{eqIV102}
\| \mathcal{K} \{ \phi^0, \phi^1 \} \|_{L^2(\Omega) \times H^{-1}(\Omega)} \leq \beta \hspace{0.6cm} \forall \{ \phi^0, \phi^1 \} \in B_\beta.
\end{equation}

One has
\begin{equation*}
    \begin{split}
    \| \Theta \{ \phi^0, \phi^1 \} \|_{L^2(\Omega) \times H^{-1}(\Omega)} \leq C \| f(y) \|_{L^1(0,T;L^2(\Omega))}
    & \leq C \| |y|^p \|_{L^1(0,T;L^2(\Omega))}
    \\ & \leq C \| y \|^p_{L^\infty(0,T;L^{2p}(\Omega))}
    \end{split}
\end{equation*}
and, as $ 2p < 2n/(n-2) $, it follows that
\begin{equation*}
\| \Theta \{ \phi^0, \phi^1 \} \|_{L^2(\Omega) \times H^{-1}(\Omega)} \leq C \| E(t) \|^{\frac{p}{2}}_{L^\infty(0,T)}
\end{equation*}
which, combined with (\ref{eqIV98}), provides the estimate
\begin{equation}\label{eqIV103}
\| \Theta \{ \phi^0, \phi^1 \} \|_{L^2(\Omega) \times H^{-1}(\Omega)} \leq C_0 \| \{ \phi^0, \phi^1 \} \|^{p}_{L^2(\Omega) \times H^{-1}(\Omega)}.
\end{equation}

Let
\begin{equation}\label{eqIV104}
 \gamma = \| \Lambda^{-1} \|_{{\cal L}({L^2(\Omega) \times H^1_0(\Omega)} , {L^2(\Omega) \times H^{-1}(\Omega)})}.
\end{equation}

From (\ref{eqIV103}) - (\ref{eqIV104}) it follows
\begin{equation}\label{eqIV105}
\begin{split}
\| \mathcal{K} \{ \phi^0, \phi^1 \} &\|_{L^2(\Omega) \times H^{-1}(\Omega)}\\
&\leq\gamma \left (\| \{ y^0, y^1 \} \|_{H^1_0(\Omega) \times L^2(\Omega)} 
+  C_0 \| \{ \phi^0, \phi^1 \} \|^p_{L^2(\Omega) \times H^{-1}(\Omega)}\right).
\end{split}
\end{equation}
This shows that (\ref{eqIV102}) is satisfied if
$
    \gamma \| \{ y^0,y^1 \} \|_{H^1_0(\Omega) \times L^2(\Omega) } + \gamma C_0 \beta^p \leq \beta
$
and it holds if
\begin{equation}\label{eqIV106}
    \beta < \min \left[ \alpha, \left( {\frac{1}{\gamma C_0}}\right)^{\frac{1}{p-1}} \right] = \sigma
\end{equation}
as long as
\begin{equation}\label{eqIV107}
\| \{ y^0, y^1 \} \|_{H^1_0(\Omega) \times L^2(\Omega)} \leq \frac{ \beta - \gamma C_0 \beta^p }{\gamma}.
\end{equation}

Applying Schauder's Theorem it follows that $ \mathcal{K} $ admits a fixed point in $B_\beta$ if the initial data $  \{ y^0, y^1 \} $ satisfies (\ref{eqIV107}) and this for each $ \beta > 0 $ fulfilling (\ref{eqIV106}).

Therefore the theorem is proved with
$
\delta = (\sigma - \gamma C_0 \sigma^p)/\gamma.
$
\end{proof}

\section{Further Comments}
\label{sec:IV.7}
\begin{enumerate}
\item {\bf Asymptotically linear nonlinearities.}  Under the hypotheses (\ref{eqIV09}) and (\ref{eqIV10}) on the nonlinearity $f$, the following result can be shown:

\begin{theorem}
\label{theorem:IV.4.6}
Let $\Omega$ be a bounded domain of $ \mathbb{R}^n $ with  boundary $\Gamma$ of class $C^2$. Let $x^0 \in \mathbb{R}^n $ and $\omega$ be a neighborhood of $ \overline{\Gamma(x^0)} \text{ in } \Omega $ and $ T>2 \| x - x^0 \|_{L^\infty(\Omega)} $. 

Suppose that $f$ satisfies (\ref{eqIV09}) and (\ref{eqIV10}).

Then, for each pair of initial and final data $ \{ y^0,y^1 \}, \{ z^0,z^1 \} \in H^1_0(\Omega) \times L^2(\Omega) $ there is a control $ h \in L^2(\omega \times (0,T)) $ such that the solution of (\ref{eqIV01}) verifies (\ref{eqIV02}).
\end{theorem}

\begin{proof} 
\label{proof:IV.theorem.4.6}
We present a sketch. One proceeds in two steps.\newline
\emph{Step 1. Linear $f$.}\label{sec:IV.7-b1}
We start by considering the linear case $ f(s) \equiv \alpha s$ with $ \alpha \in \mathbb{R}$.

We solve the equation
\begin{equation}\label{eqIV108}
    \begin{cases}
    \phi^{\prime\prime} - \Delta \phi + \alpha \phi = 0 & \text{in} \hspace{0.3cm} Q
    \\
    \phi = 0 & \text{in} \hspace{0.3cm} \Sigma
    \\
    \phi(0) = \phi^0, \phi^\prime (0) = \phi^1
    \end{cases}
\end{equation}
and then
\begin{equation}\label{eqIV109}
    \begin{cases}
    y^{\prime\prime} - \Delta y + \alpha y = -\phi \chi_\omega & \text{in} \hspace{0.3cm} Q
    \\
    y_0 = 0 & \text{in} \hspace{0.3cm} \Sigma
    \\
    y(T) = y^\prime (T) = 0.
    \end{cases}
\end{equation}

We define the operator
\begin{equation}\label{eqIV110}
 \Lambda_\alpha \{ \phi^0, \phi^1 \} = \{ y^\prime(0), -y(0) \}
\end{equation}
which satisfies
\begin{equation}\label{eqIV111}
\langle \Lambda_\alpha \{ \phi^0, \phi^1 \}, \{ \phi^0, \phi^1 \} \rangle = \int^T_0 \int_\omega |\phi|^2 dx dt.
\end{equation}

The exact controllability of (\ref{eqIV01}) in $ H^1_0(\Omega) \times L^2(\Omega) $ with controls in $ L^2(\omega \times (0,T)) $ is a consequence of the observability inequality
\begin{equation}\label{eqIV112}
\| \{ \phi^0, \phi^1 \} \|^2_{L^2(\Omega) \times H^{-1}(\Omega)}
\leq C_\alpha \int^T_0 \int_\omega |\phi|^2 dx dt.
\end{equation}

Proposition \ref{propo:IV.4.1} provides an estimate when $ \omega $ is a neighborhood of the entire  boundary $\Gamma$. However, in this particular case where the potential is constant, we have (\ref{eqIV112}) if $\omega$ is a neighborhood of $ \overline{\Gamma (x^0)} $.

Indeed, proceeding in the same way as in the proof of Proposition 4.1, it follows
\begin{equation}\label{eqIV113}
\| \phi^0 \|^2_{L^2(\Omega)} + \| \phi^1 \|^2_{H^{-1}(\Omega)}
\leq C \int^T_0 \int_\omega |\phi|^2 dx dt + C \| \eta \|^2_{L^2(Q)}
\end{equation}
with $ \eta = \eta (x,t)$ the solution of
\begin{equation}\label{eqIV114}
    \begin{cases}
    \eta ^{\prime\prime} - \Delta \eta = - \alpha \phi & \text{in} \hspace{0.3cm} Q
    \\
    \eta = 0 & \text{in} \hspace{0.3cm} \Sigma
    \\
    \eta(0) = \eta^\prime (0) = 0 & \text{in} \hspace{0.3cm} \Omega.
    \end{cases}
\end{equation}

We need to show again (\ref{eqIV37}) and we do it by contradiction. If (\ref{eqIV37}) is not satisfied, a solution $\phi$ of (\ref{eqIV108}) is obtained such that
\begin{equation}\label{eqIV115}
\phi = 0 \hspace{0.3cm} \text{ in } \hspace{0.3cm} \omega \times (0,T)
\end{equation}
and so that the associated solution $\eta$  satisfies
\begin{equation}\label{eqIV116}
\| \eta \|_{L^2(Q)} = 1.
\end{equation}

Now, Holmgren's Uniqueness Theorem (which is applicable since the coefficients of (\ref{eqIV108}) are analytic, constant actually) ensures that the only solution of (\ref{eqIV108}) verifying (\ref{eqIV115}) is $ \phi \equiv 0, $ which contradicts (\ref{eqIV116}).

\label{sec:IV.7-b2}
\emph{Step 2. Asymptotically Linear Case.}
Let $f$ be such that the following limit exists
\begin{equation*}
    \alpha = \operatorname{lim}_{|s| \rightarrow \infty} \frac{f(s)}{s}.
\end{equation*}
Set an arbitrary final state $ \{ z^0, z^1 \} \in H^1_0(\Omega) \times L^2(\Omega)  $. We solve first (\ref{eqIV108}) and then
\begin{equation}\label{eqIV117}
    \begin{cases}
    y^{\prime\prime} - \Delta y + f(y) = -\phi \chi_\omega & \text{in} \hspace{0.3cm} Q
    \\
    y = 0 & \text{in} \hspace{0.3cm} \Sigma.
    \\
    y(T) = z^0, y^\prime (T) = z^1 & \text{in} \hspace{0.3cm} \Omega.
    \end{cases}
\end{equation}

We define the nonlinear operator
$
    \mu: L^2(\Omega) \times H^{-1}(\Omega) \to L^2(\Omega) \times H^1_0(\Omega)
$
such that
$
    \mu \{ \phi^0, \phi^1 \} = \{ y^\prime(0), -y(0) \}.
$
The problem is reduced to proving that the equation
\begin{equation}\label{eqIV118}
\mu \{ \phi^0, \phi^1 \} = \{ y^1, -y^0 \}
\end{equation}
admits a solution $ \{ \phi^0, \phi^1 \} \in L^2(\Omega) \times H^{-1}(\Omega) $ for each $ \{ y^0, y^1 \} \in H^1_0(\Omega) \times L^2(\Omega).$

We decompose the solution of (\ref{eqIV117}) in the following way
$
   y = y_0 + z + \eta
$
with $ y_0, z \text{ and } \eta $ being respectively the solutions of
\begin{equation}\label{eqIV119}
    \begin{cases}
    y^{\prime\prime}_0 - \Delta y_0 + \alpha y_0 = -\phi \chi_\omega & \text{in} \hspace{0.3cm} Q
    \\
    y_0 = 0 & \text{in} \hspace{0.3cm} \Sigma
    \\
    y_0(T) = y^\prime_0 (T) = 0 & \text{in} \hspace{0.3cm} \Omega
    \end{cases}
\end{equation}

\begin{equation}\label{eqIV120}
    \begin{cases}
    z^{\prime\prime} - \Delta z + \alpha z = 0 & \text{in} \hspace{0.3cm} Q
    \\
    z = 0 & \text{in} \hspace{0.3cm} \Sigma
    \\
    z(T) = z^0, z^\prime (T) = z^1 & \text{in} \hspace{0.3cm} \Omega
    \end{cases}
\end{equation}

\begin{equation}\label{eqIV121}
    \begin{cases}
    \eta^{\prime\prime} - \Delta \eta + \alpha \eta = -f(y) + \alpha y & \text{in} \hspace{0.3cm} Q
    \\
    \eta = 0 & \text{in} \hspace{0.3cm} \Sigma
    \\
    \eta(T) = \eta^\prime (T) = 0 & \text{in} \hspace{0.3cm} \Omega.
    \end{cases}
\end{equation}

The operator $\mu$ can be rewritten as follows
\begin{equation}\label{eqIV122}
\mu \{ \phi^0, \phi^1 \} = \Lambda_\alpha \{ \phi^0, \phi^1 \} + \{ z^\prime(0), -z(0) \} + \Theta \{ \phi^0, \phi^1 \}
\end{equation}
with
\begin{equation}\label{eqIV123}
\Theta \{ \phi^0, \phi^1 \} = \{ \eta^\prime(0), -\eta(0) \}.
\end{equation}

From Step 1 it follows that
$
    \Lambda_\alpha : L^2(\Omega) \times H^{-1}(\Omega) \to L^2(\Omega) \times H^1_0(\Omega)
$
is an isomorphism.

Equation (\ref{eqIV118}) can then be rewritten in the following way
\begin{equation}\label{eqIV124}
\begin{split}
 \{ \phi^0, \phi^1 \} & = \Lambda^{-1}_\alpha \{ y^1 - z^\prime(0), z(0) - y(0) \} - \Lambda^{-1}_\alpha \Theta \{ \phi^0, \phi^1 \} = \mathcal{K} \{ \phi^0, \phi^1 \}.
\end{split}
\end{equation}

Thanks to (\ref{eqIV09}), the operator
$
    \mathcal{K} : L^2(\Omega) \times H^{-1}(\Omega) \to L^2(\Omega) \times H^{-1}(\Omega)
$
is compact.

On the other hand, using (\ref{eqIV10}) it is proved that for all $\varepsilon > 0$ there exists $C_\varepsilon(M) > 0$ such that 
\begin{equation}\label{eqIV125}
       \| \Theta \{ \phi^0, \phi^1 \} \|_{L^2(\Omega) \times H^{-1}(\Omega)} \leq \varepsilon \| \{ \phi^0, \phi^1 \} \|_{L^2(\Omega) \times H^{-1}(\Omega)} + C_{\varepsilon}.
\end{equation}

From (\ref{eqIV125}) it follows that for data $ \{ y^0, y^1 \}, \{ z^0, z^1 \} \in H^1_0(\Omega) \times L^2(\Omega) $, there exists $ R > 0 $ such that
$
    \| \mathcal{K} \{ \phi^0, \phi^1 \} \|_{L^2(\Omega) \times H^{-1}(\Omega)} \leq R$ for all $\{ \phi^0, \phi^1 \} \in B_R.
$

The proof of the theorem is concluded applying Schauder's Theorem.
\end{proof}


\item {\bf The critical exponent.} The exponent $ p = n/(n-2) $ is excluded in Theorem \ref{theorem:IV.4.5} because of the lack of compactness. However, using Banach's Fixed Point Theorem, it is possible to prove local controllability in this case as well.
$ \hfill \square $

\item {\bf Equations with variable coefficients.} The methods developed in this chapter can also be applied to semilinear wave equations with variable coefficients.
\begin{equation}\label{eqIV126}
    \begin{cases}
    y^{\prime\prime} - \operatorname{div}(a(x) \nabla y ) + f(y) = h \chi_\omega & \text{in} \hspace{0.3cm} Q
    \\
    y = 0 & \text{in} \hspace{0.3cm} \Sigma
    \\
    y(0) = y^0, y^\prime(0) = y^1 & \text{in} \hspace{0.3cm} \Omega
    \end{cases}
\end{equation}
where $ a \in W^{1,\infty}(\Omega) $ satisfies
\begin{equation}\label{eqIV127}
a(x) \geq a_0 > 0 \hspace{0.3cm} \forall x \in \Omega.
\end{equation}
We refer to \cite{fu2007exact} for an in-depth discussion.

Theorem \ref{theorem:IV.4.4} can be generalized to this case if the following conditions are satisfied:

\begin{enumerate}
\item Estimate
\begin{equation}\label{eqIV128}
\| \phi^0 \|^2_{H^1_0(\Omega)} + \| \phi^1 \|^2_{L^2(\Omega)} \leq C \int^T_0 \int_\Gamma \left| \frac{\partial \phi}{\partial v} \right|^2 d\Sigma
\end{equation}
holds for the solutions of
\begin{equation}\label{eqIV129}
    \begin{cases}
    \phi^{\prime\prime} - \operatorname{div}(a(x) \nabla \phi ) = 0 & \text{in} \hspace{0.3cm} Q
    \\
    \phi = 0 & \text{in} \hspace{0.3cm} \Sigma
    \\
    \phi(0) = \phi^0, \phi^\prime(0) = \phi^1 & \text{in} \hspace{0.3cm} \Omega.
    \end{cases}
\end{equation}

\item The following unique continuation result holds:
\begin{equation*}
\left.\begin{aligned}
        \phi^{\prime\prime} - \operatorname{div}(a(x) \nabla \phi)  + V(x,t) \phi = 0 \hspace{0.3cm} & \text{in} \hspace{0.3cm} Q
    \\
    \phi = 0 \hspace{0.3cm} & \text{in} \hspace{0.3cm} \Sigma
    \\
    \phi = 0 \hspace{0.3cm} & \text{in} \hspace{0.3cm} \omega \times (0,T)
    \\
    \phi \in H^1(Q),\hspace{0.3cm} & V \in L^\infty(Q)
    \end{aligned}
    \right\}
    \to
    \phi \equiv 0.
\end{equation*}
\end{enumerate}

The techniques in \cite{fu2007exact} combined with the fixed point arguments presented here ensure (a) and (b) to hold under further geometric restrictions on the coefficients. These added restrictions are to be expected. They ensure the applicability of Carleman inequalities, and, in particular, that all rays of geometric optics reach the control subdomain in a uniform time.

\item {\bf Dimension $n=1.$}  In the $1-d$ case, using the fixed point arguments of this section, together with the observability inequalities derived from sidewise energy estimates, exact controllability can be proved for variable coefficients of bounded variation and nonlinearities growing asymptotically slower than $|s| \log^2 |s|$ as $|s| \rightarrow \infty$ (see \cite{zuazua1990exact}).

\item {\bf An open problem. }  The main open problem in the context of the semilinear wave equation is obtaining results of global controllability with superlinear nonlinearities at infinity. This is however a delicate issue. Indeed, if the nonlinearity leads to blow-up phenomena, typically, blow-up will occur arbitrarily fast and, therefore, due to the finite speed of propagation, controls will not be able to avoid the blow-up, and, consequently, exact controllability will not hold. The results in \cite{zuazua1993exact}, \cite{fu2007exact} guarantee that exact controllability holds for a class of nonlinearities that grow slightly faster than linear, by a logarithmic factor, for which blow-up cannot occur.

However, nothing prevents, in principle, the cubic wave equation as in 
$$
   y^{\prime\prime} - \Delta y  + y^3 = h \chi_\omega
$$
from being controllable. The existing results combine stabilisation and local controllability for small data showing that all initial data may be controlled to the null state, but in a sufficiently large time, depending on the size of the initial data and which diverges as they grow (\cite{dehman}). Determining whether the exact controllability of this cubic wave equation can be achieved in a time-horizon independent of the initial datum is an open problem.
\end{enumerate}


%
%

\chapter{Wave Equation with a Nonlinear Internal Dissipation}
\label{chapter05} 

\abstract{
In this chapter we study the decay rate of the solutions of the wave equation in the presence of nonlinear damping terms distributed and acting everywhere in the interior of the domain. We introduce and use LaSalle\textquotesingle s invariance principle and build suitable Lyapunov functionals, as perturbations of the energy, allowing to get explicit decay rates as a function of the dissipative nonlinearity.}

\section{Problem Formulation}
\label{sec:V.1}
In this chapter we study the decay rate of solutions of the wave equation with nonlinear dissipation.

Let $ \Omega $ be a bounded domain of  $ \mathbb{R}^n $, $n \ge 1$, and $ g: \mathbb{R} \rightarrow \mathbb{R} $ a function such that
\begin{equation}\label{eqV01}
    g \in C(\mathbb{R}); g(0) = 0; \hspace{0.2cm} g \text{ is non-decreasing; }\hspace{0.2cm} g(s) s > 0 \hspace{0.2cm} \text{ if } \hspace{0.2cm} s \neq 0.
\end{equation}

Consider the dissipative wave equation:
\begin{equation}\label{eqV02}
    \begin{cases}
    y^{\prime\prime} - \Delta y + g(y^\prime) = 0 \hspace{0.2cm} & \text{in} \hspace{0.2cm} \Omega \times (0, \infty)
    \\
    y = 0 \hspace{0.2cm} & \text{on} \hspace{0.2cm} \partial \Omega \times (0, \infty)
    \\
    y(0) = y^0 \in H^1_0(\Omega); y^\prime (0) = y^1 \in L^2(\Omega) & \text{in} \hspace{0.2cm} \Omega.
    \end{cases}
\end{equation}

As we shall see in the next section, problem (\ref{eqV02}) admits a single solution in an appropriate weak sense.

The energy of the system is
\begin{equation}\label{eqV03}
    E_y(t) = \frac{1}{2} \int_\Omega \left[ | y^\prime(x,t) |^2 + | \nabla y(x,t) |^2  \right] dx.
\end{equation}

Multiplying equation (\ref{eqV02}) by $ y^\prime $ and integrating by parts, we obtain
\begin{equation}\label{eqV04}
    \frac{dE_y(t)}{dt} = - \int_\Omega g \left( y^\prime(x,t)\right) y^\prime (x,t) dx \leq 0.
\end{equation}

\begin{remark}
\label{remark:V.1}
These calculations are formal, but will be rigorously justified later.
$ \hfill \square $

From (\ref{eqV04}) it follows that the trajectories $ \{ y(t), y^\prime(t) \} $ are bounded in the energy space $ H^1_0(\Omega) \times L^2(\Omega) $, that is,
\begin{equation}\label{eqV05}
    y \in L^\infty(0, \infty; H^1_0(\Omega)) \cap W^{1,\infty}(0,\infty; L^2(\Omega)).
\end{equation}

Also, the energy decreases over time.

On the other hand, it is observed that the only solution of (\ref{eqV02})  in equilibrium is the trivial one, $ y \equiv 0 $, since $g(s)s=0$ only when $s=0$.

Therefore, the following two questions arise:
\begin{enumerate}
\item Determine sufficient  conditions on $ g $ so that
$
E_y(t) \rightarrow 0$ when $t \rightarrow +\infty
$
for every finite energy solution;
\item Can a uniform decay rate for the energy be obtained  as a function of the properties of nonlinearity $g$?
\end{enumerate}

In Section \ref{sec:I.2} we answered to these questions in the case where $g$ is a linear function, obtaining the uniform exponential decay. In this chapter we consider a more general class of nonlinearities $g$.

In Section \ref{sec:V.2} we state a result due to A. Haraux \cite{haraux1987semi} on the existence, uniqueness and regularity of the solutions of (\ref{eqV02}). In Section \ref{sec:V.3} we apply LaSalle's Invariance Principle to give an affirmative answer to the first question in a wide class of nonlinearities. In Section \ref{sec:V.4} we address the second one,  developing a technique based on the construction of suitable Lyapunov functions, allowing to obtain estimates of the order of uniform decay of the energy, under suitable conditions on the behavior of the nonlinearity $ g $ at the origin and at infinity. In the last section we mention some possible extensions of the results of this chapter as well as some open problems.
\end{remark}

\section{The Initial Value Problem}
\label{sec:V.2}

The following result is a consequence of Th. II.2.3.4, Prop. II.2.3.6, Th. II.1.1.1, Prop. II.1.2.2 in \cite{haraux1987semi} (see also \cite{haraux1986nonlinear}).

\begin{theorem}
\label{theorem:V.5.1}
(\cite{haraux1987semi}) Let $ \Omega $ be a bounded domain of $ \mathbb{R}^n $, $n\ge 1$,  with boundary $ \partial \Omega $ of class $ C^2$. Let $ g: \mathbb{R} \rightarrow \mathbb{R} $ be a function satisfying (\ref{eqV01}). The following holds:

\begin{enumerate}[label=\roman*)]
\item  \textbf{Existence}.  For each pair of initial data  $ \{ y^0, y^1 \} \in H^1_0(\Omega) \times L^2(\Omega) $ there is a unique solution of the above system in the following sense:
\begin{equation}\label{eqV06}
    y \in C( [0,\infty); H^1_0(\Omega)) \cap C^1( [0,\infty); L^2(\Omega) )
\end{equation}
\begin{equation}\label{eqV07}
    y(0) = y^0, y^\prime(0) = y^1
\end{equation}
\begin{equation}\label{eqV08}
    g(y^\prime) \in L^{1}_{loc}( (0,\infty); L^1(\Omega) ) \cap H^{-1}_{loc} ( \Omega \times (0,\infty) ) )
\end{equation}
\begin{equation}\label{eqV09}
    y^{\prime\prime} - \Delta y + g(y^\prime) = 0 \hspace{0.5cm} \text{in} \hspace{0.5cm} H^{-1}_{loc} ( \Omega \times (0,\infty) ) 
\end{equation}
satisfying
\begin{equation}\label{eqV11}
\begin{split}
E_y(t) \in W^{1,1}_{loc}(0, \infty) ; \frac{dE_y(t)}{dt} = - \int_\Omega g(y^\prime) y^\prime dx, \,  \, \text{a.e.} \,  \, t \in (0, \infty).    
\end{split}
\end{equation}

\item \textbf{Stability}. For any pair of finite energy solutions $ y $ and  $ \hat{y} $ the following contraction property holds:
\begin{equation}\label{eqV12}
    E_{y - \hat{y}} (t) \leq E_{y - \hat{y}} (0) \hspace{0.3cm} \forall t \geq 0.
\end{equation}

\item \textbf{Regularity}. If the initial data are more regular, namely, $ \{ y^0, y^1 \} \in [ H^2(\Omega) \cap H^1_0(\Omega) ] \times H^1_0(\Omega) $, and $ g(y^1) \in L^2(\Omega) $, then the solution belongs to the class
\begin{equation}\label{eqV13}
    y \in W^{2, \infty} (0, \infty; L^2(\Omega)) \cap W^{1, \infty} (0, \infty; H^1_0(\Omega)) \cap L^\infty (0, \infty; H^2(\Omega)).
\end{equation}
\end{enumerate}
\end{theorem}

\begin{remark}
\label{remark:V.5.2}Several remarks are in order:
\begin{itemize}

\item From now, when referring to a solution of (\ref{eqV02}), it will be understood that it is the only one satisfying (\ref{eqV06}) - (\ref{eqV11}).

\item In (\ref{eqV08}), by $ g(y^\prime) \in H^{-1}_{loc} (\Omega \times (0, \infty)) $ we mean that
for all $ 0 < t_1 < t_2 < \infty , g(y^\prime) \in H^{-1} (\Omega \times (t_1, t_2)).$ 

\item The uniform in time estimates in \eqref{eqV13} can be obtained by means of a higher order energy estimate. Roughly, multiplying the equation satisfied by $y$ by $-\Delta y^\prime$ and integrating by parts, we get
\begin{equation}
 \frac{dE^+_y(t)}{dt} = - \int_\Omega g'(y^\prime (x,t)) |\nabla y^\prime |^2 dx
\end{equation}
where $E^+_y$ is the higher order energy
\begin{equation}\label{eqV14}
    E^+_y(t) = \frac{1}{2} \int_\Omega \left[ | \nabla y^\prime(x,t) |^2 + | \Delta y(x,t) |^2  \right] dx.
\end{equation}
This is so since, by integration by parts,
$$
- \int_\Omega g(y^\prime(x,t)) \Delta y^\prime(x,t) dx = \int_\Omega g^\prime(y^\prime (x,t)) |\nabla y^\prime(x,t) |^2 dx.
$$
Of course, the monotonicity of $g$ plays a key role when assuring the dissipativity of this higher order energy.

$ \hfill \square $
\end{itemize}
\end{remark}

\section{LaSalle's Invariance Principle and its Consequences}
\label{sec:V.3}
The following result responds to question (1) raised in Section \ref{sec:V.1}.

\begin{theorem}
\label{theorem:V.5.2}
Suppose that the hypotheses of Theorem \ref{theorem:V.5.1} are satisfied and that, in addition, \eqref{eqV01} is fulfilled.

Then, for every finite energy solution $ y = y (x,t) $ of (\ref{eqV01}) with initial data $ \{ y^0, y^1 \} \in H^1_0(\Omega) \times L^2(\Omega) $ we have
\begin{equation}\label{eqV15}
    E_y(t) \rightarrow 0, t \rightarrow + \infty.
\end{equation}
\end{theorem}

\begin{proof}
\label{proof:V.theorem.5.2}
We proceed in two steps.
\newline
\label{sec:V.3-c}
\emph{Step 1. Regular Solutions}. Suppose that the initial data are such that
\begin{equation}\label{eqV16}
\{ y^0, y^1 \} \in [ H^2(\Omega) \cap H^1_0(\Omega) ] \times H^1_0(\Omega); g(y^1) \in L^2(\Omega).
\end{equation}

From Theorem \ref{theorem:V.5.1}, 
$
\{ y(t), y^\prime(t)  \}_{t \geq 0} \text{ is bounded in } \left [H^2(\Omega) \cap H^1_0(\Omega)  \right  ]\times H^1_0(\Omega) 
$
and therefore
$
\{ y(t), y^\prime(t)  \}_{t>0} \text{ is relatively compact in } H^1_0(\Omega) \times L^2(\Omega).
$

It is enough to prove that the only  accumulation point of the trajectory as $t\to \infty$ in $ H^1_0(\Omega) \times L^2(\Omega) $ is $ \{ 0,0 \} $.

Let $ t_n \rightarrow + \infty $ be such that
\begin{equation}\label{eqV17}
\{ y(t_n), y^\prime(t_n) \} \rightarrow \{ z^0, z^1 \} \text{ in } H^1_0(\Omega) \times L^2(\Omega).
\end{equation}

Let $ z_n = z_n (x,t) $ be the function
\begin{equation}\label{eqV18}
    z_n(x,t) = y(x, t + t_n).
\end{equation}

Given $ T>0 $,
$
\{ z_n \} \text{ is bounded in } W^{2, \infty} (0, T; L^2(\Omega)) \cap W^{1,\infty} (0,T; H^1_0(\Omega)),
$
and therefore, extracting a subsequence that we continue to denote by $ \{ z_n \} $, it follows that
\begin{equation}\label{eqV19}
    z_n \rightarrow z \hspace{0.3cm} \text{ in } \hspace{0.3cm} C( [0,T]; H^1_0(\Omega) ) \cap C^1( [0,T]; L^2(\Omega) )
\end{equation}
with $ z = z(x,t) $, which, thanks to (\ref{eqV17}), satisfies
\begin{equation}\label{eqV20}
    z(0) = z^0, z^\prime(0) = z^1.
\end{equation}

On the other hand
\begin{equation}\label{eqV21}
    \begin{split}
    {-} \int^T_0 \int_\Omega g(z^\prime_n)z^\prime_n dx dt 
    & = {-} \int^T_0 \int_\Omega g( y^\prime (x, t_n + t) )y^\prime (x, t_n + t) dt
    \\ & = E_y(t_n + T) - E_y(t_n).     
    \end{split}
\end{equation}

Since $ E_y(\cdot) $ is decreasing and non-negative, the limit
$
\lim_{t \rightarrow \infty} E_y(t) = \ell \ge 0
$
exists, and therefore, from (\ref{eqV21}), it follows that
\begin{equation}\label{eqV22}
\lim_{n \rightarrow \infty} \int^T_0 \int_\Omega g(z^\prime_n)z^\prime_n dx dt = 0.
\end{equation}

Since $ g $ is non-decreasing and satisfies (\ref{eqV01}) we have
\begin{equation}\label{eqV23}
\begin{split}
\int^T_0 \int_\Omega |g(z_n^\prime)| dx dt \leq 
\int^T_0 \int_{\Omega \cap \{|z_n^\prime|\geq 1\}} g(z_n^\prime) z_n^\prime dx dt + \int^T_0 \int_{\Omega \cap \{|z_n^\prime|\leq 1\}} |g(z_n^\prime)| dx dt.
\end{split}
\end{equation}

Combining (\ref{eqV22}) and (\ref{eqV23}) it follows that
$
g(z_n^\prime) \to 0$ in $L^1(\Omega \times (0,T)).
$

As $ z_n = z_n (x,t) $ solves
$
z^{\prime\prime}_n - \Delta z_n + g(z^\prime_n) = 0 $ in $\Omega \times (0,T),$
it follows that
\begin{equation}\label{eqV24}
z^{\prime\prime} - \Delta z = 0 \hspace{0.2cm} \text{ in } \hspace{0.2cm} \Omega \times (0,T)
\end{equation}
and, moreover,
$
g(z^\prime) = 0$  in $ \Omega \times (0,T),
$
which, by (\ref{eqV01}), implies
\begin{equation}\label{eqV25}
z^\prime = 0 \hspace{0.2cm} \text{ in } \hspace{0.2cm} \Omega \times (0,T).
\end{equation}

Combining (\ref{eqV24}) and (\ref{eqV25}) we see that $ z = z(x) \in H^1_0(\Omega) $ and also that $ \Delta z = 0 $ in $ \Omega $, which implies $ z \equiv 0 $, and therefore
$
 z^0 = z^1 = 0.
$

This concludes the proof of the Theorem in the case where the initial data are regular.

\label{sec:V.3-d}
\emph{Step 2. General Case}.
Given $ \{ y^0, y^1 \} \in H^1_0(\Omega) \times L^2(\Omega) $, by density, we may construct a sequence of data $ \{ y^0_n, y^1_n \} \in [ H^2(\Omega) \cap H^1_0(\Omega) ] \times H^1_0(\Omega) $ such that $ g(y^1_n) \in L^2(\Omega) $ and
\begin{equation}\label{eqV26}
y^0_n \to y^0 \text{ in } H^1_0(\Omega), y^1_n \to y^1 \text{ in } L^2(\Omega).
\end{equation}

From (\ref{eqV12}) we get
\begin{equation}\label{eqV27}
\begin{split}
E_y(t) & \leq E_{y-y_n}(t) + E_{y_n}(t) \leq E_{y-y_n}(0) + E_{y_n}(t).
\end{split}
\end{equation}

Given $ \varepsilon > 0 $ we set $ n \in \mathbb{N} $ such that
\begin{equation}\label{eqV28}
E_{y-y_n}(0) < \frac{\varepsilon}{2},
\end{equation}
which is possible thanks to (\ref{eqV26}). Next, for fixed $n$, as shown in the previous step, we see that there is $ t_0 > 0 $ such that
\begin{equation}\label{eqV29}
E_{y_n}(t) < \frac{\varepsilon}{2} \hspace{0.6cm} \forall t \geq t_0 .
\end{equation}

Combining (\ref{eqV27}) - (\ref{eqV29}) it follows that $ E_y(t) \rightarrow 0 $ when $ t \rightarrow + \infty $.

This concludes the proof of the theorem.
\end{proof}

\section{Decay Rates}
\label{sec:V.4}
In the following theorem we prove some upper bounds for the decay rate of the energy as a function of the coercivity of $g$ at the origin and its growth at infinity.

\begin{theorem}
\label{theorem:V.5.3}
Suppose that the conditions of Theorem \ref{theorem:V.5.1} are satisfied and that, furthermore:

(a) There exist $ p \ge 1, c>0 $ such that
\begin{equation}\label{eqV30}
g(s) s \geq c |s|^{p+1}, \quad \forall s \in [-1,1]
\end{equation}
\begin{equation}\label{eqV31}
g(s) s \geq c |s|^2, \quad \forall s \in \mathbb{R}: |s| \geq 1;
\end{equation}
(b) There exist $ \lambda > 0, C>0, q>1 $ such that 
\begin{equation}\label{eqV32}
| g(s) | \leq C |s|^\lambda, \quad \forall s \in [-1,1]
\end{equation}
\begin{equation}\label{eqV33}
| g(s) | \leq C |s|^q, \quad \forall s \in \mathbb{R}: |s| \geq 1 \text{ with } (n-2)q \leq (n+2).
\end{equation}

Then,
\begin{enumerate}[label=\roman*)]
\item{If $ \lambda = p = 1 $, there are constants $ M > 1, \gamma > 0 $ independent of the initial data such that
\begin{equation}\label{eqV34}
E(t) \leq M E(0)e^{-\gamma t} \hspace{0.6cm} \forall t \geq 0
\end{equation}
for every solution of (\ref{eqV02});
}

\item{If $ \lambda < 1 $, there is a constant $M$ depending on $E_y(0)$ such that
\begin{equation}\label{eqV35}
E_y(t) \leq 4 \left[ M t + ( E_y(0))^{-\frac{p+1-2\lambda}{2\lambda}} \right]^{\frac{-2\lambda}{(p+1-2\lambda)}} \hspace{0.6cm} \forall t \geq 0;
\end{equation}
}

\item{If $ \lambda \geq 1 $ and $ p>1 $, there is a constant $ M>0 $ depending on $ E_y(0) $ such that
\begin{equation}\label{eqV36}
E_y(t) \leq 4 \left[ M t + ( E_y(0))^{-\frac{(p-1)}{2}} \right]^{\frac{-2}{(p-1)}} \hspace{0.6cm} \forall t \geq 0.
\end{equation}
}
\end{enumerate}

\end{theorem}

\begin{remark}
\label{remark:V.5.4}
Conditions (\ref{eqV30}) and (\ref{eqV31}) ensure the coercivity of $g$ at the origin and at infinity. On the other hand, (\ref{eqV32}) and (\ref{eqV33}) control and limit its growth.

It is important to note that the decay rate depends only on $p$ and $\lambda$, that is, on the behavior of $g$ at the origin.

The case $ \lambda = p = 1 $ contains the linear one and leads to the exponential decay. In the other cases a polynomial decay order is obtained.
$ \hfill \square $
\end{remark}

\begin{remark}
\label{remark:V.5.5}
Note that in dimension $ n=2 $, thanks in the embedding of $H^1_0(\Omega)$ in all $L^p(\Omega)$ spaces, with $1\le p<\infty$,  it is sufficient that (\ref{eqV33}) is satisfied for some $q>1$.

In dimension $n=1$ the hypothesis (\ref{eqV33}) can be avoided. Indeed, thanks to the continuity of the embedding $ H^1_0(\Omega) \subset L^\infty(\Omega) $, it suffices $g$ to be continuous.
$ \hfill \square $
\end{remark}

\begin{remark}
\label{remark:V.5.6}
We are going to use a method based on the construction of a perturbed energy or Lyapunov functional, equivalent to the energy $E_y(\cdot)$, and for which differential inequalities leading to the decay rates stated in the theorem can be obtained.

The estimates (\ref{eqV34}) and (\ref{eqV36}) were proved in  \cite{zuazua1988stability} by generalizing previous results in \cite{haraux1988decay}. In the proof of (\ref{eqV35})  we were inspired by the work of J. Lagnese and G. Leugering \cite{lagnese1991uniform} for the construction of the appropriate perturbed energy.
$ \hfill \square $
\end{remark}

\begin{proof}
\label{proof:V.theorem.5.6}
We consider regular initial data $ \{ y^0, y^1 \} \in \left[H^2(\Omega) \cap H^1_0(\Omega)\right] \times H^1_0(\Omega)$, $g(y^1) \in L^2(\Omega). $
The estimates will be extended to the general case by a simple density argument.
We distinguish the cases $ \lambda \geq 1 $ and $ \lambda < 1 $.
\smallskip

\noindent {\bf Case 1. $ \lambda \geq 1. $}
\label{case.1:proof:V.theorem.5.6}
We define the function
\begin{equation}\label{eqV37}
\phi(t) = [ E_y(t)]^\frac{(p-1)}{2} \int_\Omega y(x,t) y^\prime (x,t) dx.
\end{equation}

In what follows the energy will be denoted by $E(t)$ instead of $E_y(t)$, as we work with an arbitrary but fixed solution $y$.

It holds
\begin{equation}\label{eqV38}
\begin{split}
\frac{d \phi}{dt}(t) &= \frac{p-1}{2} \left[ E(t) \right]^\frac{p-3}{2} E^\prime(t) \int_\Omega y y^\prime dx + \left[ E(t) \right]^\frac{p-1}{2} \int_\Omega \left[ |y^\prime|^2 +y y^{\prime\prime} \right] dx.
\end{split}
\end{equation}

As
$
 \left| \int_\Omega y y^\prime dx \right| \leq CE(t)
$
and
$
 E(t) \leq E(0) 
$ for all $t \geq 0$,
given that $p\ge 1$, it follows that
\begin{equation*}
\frac{p-1}{2} \left[ E(t) \right]^\frac{p-3}{2} E^\prime (t) \int_\Omega y y^\prime dx \leq - C \left[ E(t) \right]^\frac{p-1}{2} E^\prime (t)\leq-C_1 E^\prime (t).
\end{equation*}

On the other hand, using equation (\ref{eqV02}) we see that
\begin{equation*}
\int_\Omega \left[ |y^{\prime}|^2 + y y^{\prime\prime} \right] dx = \int_\Omega \left[ |y^\prime|^2 - |\nabla y|^2 - yg(y^\prime) \right] dx.
\end{equation*}

From (\ref{eqV38}) we deduce
\begin{equation}\label{eqV39}
\phi^\prime(t) \leq -C_1 E^\prime(t) + \left[ E(t) \right]^\frac{p-1}{2} \int_\Omega \left[ |y^\prime|^2 - |\nabla y|^2 - y g(y^\prime) \right] dx.
\end{equation}

On the other hand,
\begin{equation}\label{eqV40}
\begin{split}
\left| \int_\Omega g (y^\prime)y dx \right| & = \left| \langle g(y^\prime), y \rangle_{H^{-1}, H^1_0} \right| \leq \| g (y^\prime) \|_{H^{-1}(\Omega)} \| y \|_{H^1_0(\Omega)}.
\end{split}
\end{equation}
At this level we need the following Lemma.

\begin{lemma}
\label{lemma:V.5.1}
There are constants $ C_2, C_3 > 0 $, independent of the solution, such that
\begin{equation}\label{eqV41}
 \| g(v) \|_{H^{-1}(\Omega)} \leq C_2 \left[ \int_\Omega g(v) v dx \right]^\frac{q}{q+1} + C_3 \| v \|_{L^2(\Omega)}, \hspace{0.6cm} \forall v \in H^1_0(\Omega).
\end{equation}
\end{lemma}

\begin{proof}
\label{proof:V.lemma.5.1}
First, we observe that, as $ (n-2) q \leq n+2 $, then $ L^{(q+1)/q}(\Omega) \subset H^{-1}(\Omega) $ with continuous  embedding.

From (\ref{eqV33}) it follows that
\begin{equation}\label{eqV42}
 |g(s)|^{q+1} \leq C |s|^q |g(s)|^q \hspace{0.6cm} \forall s \in \mathbb{R}: |s| \geq 1.
\end{equation}

From (\ref{eqV32}) and (\ref{eqV42}) we deduce, using H{\"o}lder's inequality, that
\begin{equation}\label{eqV43}
    \begin{split}
    \| g(v) \|_{H^{-1}(\Omega)} & \leq C \| g(v) \|_{L^{(q+1)/q}(\Omega)}
  \leq C \| v \|_{L^{(q+1)/q}(\Omega)} + C \| g(v) v \|^{q/(q+1)}_{L^1(\Omega)}.
    \end{split}
\end{equation}
This leads to (\ref{eqV41}).
\end{proof}
From (\ref{eqV40}) and (\ref{eqV41}) we obtain
\begin{equation}\label{eqV44}
    \begin{split}
    \bigg| \int_\Omega & g(y^\prime)y dx \bigg| \leq C_3 \| y^\prime \|_{L^2(\Omega)} \| y \|_{H^1_0(\Omega)}
    + C_2 \left[ \int_\Omega g(y^\prime) y^\prime dx \right]^{\frac{q}{(q+1)}} \| y \|_{H^1_0(\Omega)}
    \\  \leq &\frac{1}{2} \| y \|^2_{H^1_0(\Omega)} + C_4 \| y^\prime \|^2_{L^2(\Omega)}
     + C_2 \left[ \int_\Omega g(y^\prime) y^\prime dx \right]^{\frac{q}{(q+1)}} \| y \|_{H^1_0(\Omega)}.
    \end{split}
\end{equation}

Combining (\ref{eqV39}) and (\ref{eqV44}) it follows
\begin{equation}\label{eqV45}
    \begin{split}
    \phi^\prime (t)  \leq &{-C_1} E^\prime (t) + (C_4 + 1) |E(t)|^{\frac{p-1}{2}} \| y^\prime (t) \|^2_{L^2(\Omega)}
    - \frac{1}{2} [E(t) ]^{\frac{p-1}{2}} \| y(t) \|^2_{H^1_0(\Omega)}
    \\ & + C_2 [E(t) ]^{\frac{p-1}{2}} \left[ \int_\Omega g(y^\prime) y^\prime dx \right]^{\frac{q}{(q+1)}} \| y(t) \|_{H^1_0(\Omega)}.
    \end{split}
\end{equation}

Applying Young's inequality with exponents $ (q+1)/q $ and $ q+1 $ to the last term of this inequality we obtain
\begin{equation}\label{eqV46}
    \begin{split}
    \phi^\prime (t)  \leq &{-C_1} E^\prime (t) + (C_4 + 1) |E(t)|^{\frac{p-1}{2}} \| y^\prime (t) \|^2_{L^2(\Omega)}
    \\ & - \frac{1}{4} [E(t) ]^{\frac{p-1}{2}} \| y(t) \|^2_{H^1_0(\Omega)}
     + C_5 \int_\Omega g(y^\prime) y^\prime dx 
    \\  \leq &-(C_1 + C_5) E^\prime(t) + C_4 [E(t)]^{\frac{p-1}{2}}  \| y^\prime (t) \|^2_{L^2(\Omega)}
    - \frac{1}{4} [E(t) ]^{\frac{p-1}{2}} \| y(t) \|^2_{H^1_0(\Omega)}
    \end{split}
\end{equation}
where we used the energy dissipation identity
\begin{equation*}
    E^\prime(t) = - \int_\Omega g(y^\prime (x,t)) y^\prime (x,t)dx.
\end{equation*}

We define the perturbed energy
\begin{equation}\label{eqV47}
    E_\varepsilon (t) = ( 1 + \varepsilon (C_1 + C_5) ) E(t) + \varepsilon \phi (t).
\end{equation}

Note that there exists $ \varepsilon_0 > 0 $ such that
\begin{equation}\label{eqV48}
    \frac{1}{2} [ E_\varepsilon (t) ]^{\frac{p+1}{2}} \leq [E(t)]^{\frac{p+1}{2}} \leq 2[ E_\varepsilon(t) ]^{\frac{p+1}{2}} \hspace{0.3cm} \forall t \geq 0, \forall \varepsilon \leq \varepsilon_0.
\end{equation}

When $ p > 1, \varepsilon_0 $ continuously depends on $E(0)$. However, when $ p = 1, \varepsilon_0 $ can be chosen independently of $ E(0) $.

We then distinguish the cases $ p = 1, p > 1 $.
\smallskip

\noindent {$ \text{\bf Case 1.a. } \lambda \geq 1, p = 1.$}
\label{case.1.a:proof:V.theorem.5.6}
It holds that
\begin{equation}\label{eqV49}
    E^\prime_\varepsilon(t) \leq E^\prime(t) + \varepsilon C_4 \| y^\prime(t) \|^2_{L^2(\Omega)} - \frac{\varepsilon}{4} \| y(t) \|^2_{H^1_0(\Omega)}.
\end{equation}
As
\begin{equation*}
    E^\prime(t) = - \int_\Omega g(y^\prime) y^\prime dx \leq -C_6 \| y^{\prime}(t) \|^2_{L^2(\Omega)}
\end{equation*}
taking
$
\varepsilon < C_6/2C_4 = \varepsilon_1
$
we have that
\begin{equation*}
E^\prime_\varepsilon (t) \leq - \frac{C_6}{2} \| y^\prime(t) \|^2_{L^2(\Omega)} - \frac{\varepsilon}{4} \| y(t) \|^2_{H^1_0(\Omega)}
\end{equation*}
and, if
$
\varepsilon < 2 C_6 = \varepsilon_2,
$
we obtain
\begin{equation}\label{eqV50}
    E^\prime_\varepsilon(t) \leq -\frac{\varepsilon}{2} E(t) \leq - \frac{\varepsilon}{4} E_\varepsilon (t).
\end{equation}
From (\ref{eqV50}) we deduce that if
$
\varepsilon < \min (\varepsilon_0, \varepsilon_1, \varepsilon_2)
$
we have
\begin{equation*}
E_\varepsilon (t) \leq E_\varepsilon (0) e^{-\frac{\varepsilon}{4}t} \hspace{0.6cm} \forall t \geq 0.
\end{equation*}
This leads to
\begin{equation}\label{eqV51}
E (t) \leq 4 E(0) e^{-\frac{\varepsilon}{4}t} \hspace{0.6cm} \forall t \geq 0.
\end{equation}

\noindent {$ \text{\bf Case 1.b. } \lambda \geq 1, p > 1. $} \label{case.1.b:proof:V.theorem.5.6}
In this case we have
\begin{equation}\label{eqV52}
\begin{split}
E^\prime_\varepsilon(t) \leq - \int_\Omega g(y^\prime) y^\prime dx - \frac{\varepsilon}{4} ( E(t) )^{\frac{p-1}{2}} \| y(t) \|^2_{H^1_0(\Omega)} & + \varepsilon C_4 ( E(t))^{\frac{p-1}{2}} \| y^\prime(t) \|^2_{L^2(\Omega)}.
\end{split}
\end{equation}

On the other hand, from (\ref{eqV30}) and (\ref{eqV31}) it follows that there exists $ C_6 > 0 $, depending on the solution, such that
\begin{equation}\label{eqV53}
\int_\Omega g(y^\prime) y^\prime dx \geq C_6 \| y^\prime (t) \|^{p+1}_{L^2(\Omega)}.
\end{equation}

Indeed
\begin{equation*}
\begin{split}
\int_\Omega g(y^\prime) y^\prime dx & \geq C \int_{ \{|y^\prime| \leq 1 \} } |y^\prime|^{p+1} dx + C \int_{ \{|y^\prime| \geq 1 \} } |y^\prime|^2 dx  
\\ &\geq \frac{C}{|\Omega|^{\frac{p-1}{2}} } \left[ \int_{ \{|y^\prime| \leq 1 \} } |y^\prime|^2 dx \right]^{\frac{p+1}{2}} + \frac{C}{[2E(0)]^{\frac{p-1}{2}} } \left[ \int_{ \{|y^\prime| \geq 1 \} } |y^\prime|^2 dx \right]^{\frac{p+1}{2}}
\end{split}
\end{equation*}
from which (\ref{eqV53}) follows with
$
C_6 = C  [\max \{ |\Omega|, 2E(0) \} ]^{- \frac{(p-1)}{2}},
$
where $ | \Omega | $ denotes the Lebesgue measure of $ \Omega $ in $ \mathbb{R}^n $.

Combining (\ref{eqV52}) and (\ref{eqV53}) it follows that
\begin{equation}\label{eqV54}
\begin{split}
E^\prime_\varepsilon (t) & \leq -C_6 \| y^\prime(t) \|^{p+1}_{L^2(\Omega)} + \varepsilon C_4 \left [ E(t) \right]^{\frac{p-1}{2}} \| y^\prime(t) \|^2_{L^2(\Omega)} - \frac{\varepsilon}{4} \left [ E(t) \right ]^{\frac{p-1}{2}} \| y(t) \|^2_{H^1_0(\Omega)} 
\\ & \le -C_6 \| y^\prime(t) \|^{p+1}_{L^2(\Omega)} - \frac{\varepsilon}{2} \left [ E(t) \right ]^{\frac{p+1}{2}} + \varepsilon (C_4 + \frac{1}{4}) \left [ E(t) \right ]^{\frac{p-1}{2}} \| y^\prime(t) \|^2_{L^2(\Omega)}.
\end{split}
\end{equation}

However
\begin{equation}\label{eqV55}
\left [ E(t) \right ]^{\frac{p-1}{2}} \| y^\prime \|^2_{L^2(\Omega)} \leq \mu^{\frac{p+1}{2}} \| y^\prime \|^{p+1}_{L^2(\Omega)} + \mu^{\frac{p+1}{p-1}} \left [ E(t) \right ]^{\frac{p+1}{2}}
\end{equation}
for all $ \mu > 0 $. Combining (\ref{eqV54}) and (\ref{eqV55}) with $ \mu > 0 $ small enough so that
\begin{equation*}
    \mu^{\frac{p+1}{p-1}} (C_4 + \frac{1}{4}) \leq \frac{1}{4}
\end{equation*}
we get
\begin{equation*}
    E_\varepsilon^\prime(t) \leq - \frac{\varepsilon}{4} \left [ E(t) \right ]^{\frac{p+1}{2}} - \left[ C_6 - \varepsilon \left( \frac{4C_4 + 1}{4 \mu^{(p+1)/2}} \right) \right] \| y^\prime(t) \|^2_{L^2(\Omega)}
\end{equation*}
and finally, choosing, 
\begin{equation*}
\varepsilon \leq \varepsilon_1 = \frac{4 C_6 \mu^\frac{(p+1)}{2}}{4 C_4 + 1}
\end{equation*}
it follows
\begin{equation}\label{eqV56}
    E^\prime_\varepsilon(t) \leq - \frac{\varepsilon}{4} \left [ E(t) \right ]^{\frac{p+1}{2}} \leq - \frac{\varepsilon}{8} \left [ E_\varepsilon(t) \right ]^{\frac{p+1}{2}}.
\end{equation}

Integrating the inequality (\ref{eqV56}) we obtain
\begin{equation}\label{eqV57}
\begin{split}
    E_\varepsilon(t) & \leq \left[ \frac{\varepsilon (p-1)}{16} t + [E_\varepsilon (0)]^{-\frac{p-1}{2}} \right]^{-\frac{2}{p-1}}
    \leq 2^{\frac{2}{p+1}} \left[ \frac{\varepsilon (p-1)}{2^{(3p + 5)/(p+1)}} t + [E(0)]^{-\frac{p-1}{2}} \right]^{-\frac{2}{p-1}}
\end{split}
\end{equation}
from which (\ref{eqV36}) follows.
\smallskip

\noindent {$ \text{\bf Case 2. } \lambda < 1.$}
\label{case.2:proof:V.theorem.5.6}
We define the function
\begin{equation}\label{eqV58}
    \phi(t) = [E(t)]^{\frac{(p+1-2\lambda)}{2\lambda}} \int_\Omega y(x,t) y^\prime (x,t) dx.
\end{equation}

From (\ref{eqV30}) and (\ref{eqV32}) it follows that $ p \geq \lambda $ and therefore $ p+1-2\lambda > 0 $.

We have
\begin{equation}\label{eqV59}
\begin{split}
\phi^\prime (t) & = \frac{p+1-2\lambda}{2\lambda} \left [ E(t) \right ]^{\frac{p+1-4\lambda}{2\lambda}} \left[ \int_\Omega y y^\prime dx \right] E^\prime (t)
 + \left [ E(t) \right ]^{\frac{p+1-2\lambda}{2\lambda}}  \int_\Omega \left[ |y^\prime|^2 + y y^{\prime\prime} \right] dx 
\\ & \leq -c_1 E^\prime(t) + \left [ E(t) \right ]^{\frac{p+1-2\lambda}{2\lambda}}  \int_\Omega \left[ |y^\prime|^2 - | \bigtriangledown y |^2 + g(y^\prime)y  \right] dx.
\end{split}
\end{equation}

On the other hand
\begin{equation*}
    \begin{split}
    \left| \int_{|y^\prime|\leq 1} g(y^\prime) y dx \right| & \leq c \int_{|y^\prime|\leq 1} |y^\prime|^\lambda |y| dx \leq c \int_{|y^\prime|\leq 1} \left[ g(y^\prime) y^\prime \right]^{\frac{\lambda}{p+1}} |y| dx
    \\ & \leq c \left[ \int_\Omega g(y^\prime) y^\prime dx \right]^{\frac{\lambda}{p+1}} \left[ \int_\Omega |y|^{\frac{(p+1)}{(p+1-\lambda)}} dx \right]^{\frac{(p+1-\lambda)}{(p+1)}}
    \end{split}
\end{equation*}
but, as $ (p+1)/(p+1-\lambda) < 2 $, we get
\begin{equation}\label{eqV60}
    \left| \int_{|y^\prime|\leq 1} g(y^\prime) y dx \right| \leq c \| y \|_{H^1_0(\Omega)} \left[ \int_\Omega g(y^\prime) y^\prime dx \right]^{\frac{\lambda}{(p+1)}} 
  \leq c [E(t)]^\frac{1}{2} \left[ \int_\Omega g(y^\prime) y^\prime dx \right]^{\frac{\lambda}{(p+1)}}.
\end{equation}

Likewise, arguing like in the proof of Lemma \ref{lemma:V.5.1}, we obtain
\begin{equation}\label{eqV61}
\left| \int_{|y^\prime|\geq 1} g(y^\prime) y dx \right| \leq c  \left[ \int_\Omega g(y^\prime) y^\prime dx \right]^{\frac{q}{(q+1)}} \| y(t) \|_{H^1_0(\Omega)}.
\end{equation}
Combining (\ref{eqV59}) - (\ref{eqV61}) we get
\begin{equation*}
    \begin{split}
        \phi^{\prime}(t)  \leq &-c_1 E^\prime (t) + \left [ E(t) \right ]^{\frac{p+1-2\lambda}{2\lambda}} \| y^\prime (t) \|^2_{L^2(\Omega)}
        - \left [ E(t) \right ]^{\frac{p+1-2\lambda}{2\lambda}} \| y(t) \|^2_{H^1_0(\Omega)}
        \\ & + c_2 \left [ E(t) \right ]^{\frac{p+1-\lambda}{2\lambda}} \left[ \int_\Omega g(y^\prime) y^\prime dx \right]^{\frac{\lambda}{(p+1)}}
        \\ & + c_2 \left [ E(t) \right ]^{\frac{p+1-2\lambda}{2\lambda}} \| y(t) \|_{H^1_0(\Omega)} \left[ \int_\Omega g(y^\prime) y^\prime dx \right]^{\frac{q}{(q+1)}}.
    \end{split}
\end{equation*}
Applying Young's inequality with exponents $ q+1 $ and $ (q+1)/q $ in the last term of the previous inequality, we obtain
\begin{equation*}
    \begin{split}
    \phi^{\prime}(t)  \leq &-(c_1 + c_3) E^\prime (t) + \left [ E(t) \right ]^{\frac{p+1-2\lambda}{2\lambda}} \| y^\prime (t) \|^2_{L^2(\Omega)}
    \\ & - \frac{1}{2} \left [ E(t) \right ]^{\frac{p+1-2\lambda}{2\lambda}} \| y(t) \|^2_{H^1_0(\Omega)}
     + c_2 \left [ E(t) \right ]^{\frac{p+1-\lambda}{2\lambda}} \left[ \int_\Omega g(y^\prime) y^\prime dx \right]^{\frac{\lambda}{(p+1)}}
    \\  = &-(c_1 + c_3) E^\prime (t) - \left [ E(t) \right ]^{\frac{p+1}{2\lambda}} + \frac{3}{2} \left [ E(t) \right ]^{\frac{p+1-2\lambda}{2\lambda}} \|y^\prime(t)\|^2_{L^2(\Omega)}
    \\  &+ c_2 \left [ E(t) \right ]^{\frac{p+1-\lambda}{2\lambda}} \left[ \int_\Omega g(y^\prime) y^\prime dx \right]^{\frac{\lambda}{(p+1)}}.
    \end{split}
\end{equation*}
Applying Young's inequality with  exponents $ (p+1)/(p+1-\lambda) $ and $ (p+1)/\lambda $ in the last term, we deduce
\begin{equation*}
\begin{split}
    \phi^{\prime}(t) \leq & -(c_1 + c_3 + c_4) E^\prime (t) - \frac{1}{2} \left [ E(t) \right ]^{\frac{p+1}{2\lambda}} 
    + \frac{3}{2} \left [ E(t) \right ]^{\frac{p+1-2\lambda}{2\lambda}} \| y^\prime(t) \|^2_{L^2(\Omega)}.
\end{split}
\end{equation*}

Applying again Young's inequality in the last term with $ (p+1)/(p+1-2\lambda) $ and $ (p+1)/2\lambda $  we obtain
\begin{equation}\label{eqV62}
\begin{split}
    \phi^{\prime}(t) \leq & -(c_1 + c_3 + c_4) E^\prime (t) - \frac{1}{4} \left [ E(t) \right ]^{\frac{p+1}{2\lambda}} 
     + c_5 \| y^\prime(t) \|^{\frac{p+1}{\lambda}}_{L^2(\Omega)}.
\end{split}
\end{equation}

Finally we observe that
\begin{equation}\label{eqV63}
-E^\prime(t) = \int_\Omega g(y^\prime)y^\prime dx \geq c_6 \| y^\prime (t) \|^{p+1}_{L^2(\Omega)} \geq c_7 \| y^\prime (t) \|^{\frac{p+1}{\lambda}}_{L^2(\Omega)}.
\end{equation}

Combining (\ref{eqV62}) and (\ref{eqV63}) we deduce
\begin{equation}\label{eqV64}
\phi^{\prime}(t) \leq - c_8 E^\prime (t) - \frac{1}{4} \left [ E(t) \right ]^{\frac{p+1}{2\lambda}}
\end{equation}
with $ c_8 = c_1 + c_3 + c_4 + c_5/c_7. $

We define perturbed energy
\begin{equation}\label{eqV65}
E_\varepsilon(t) = (1 + \varepsilon c_8) E(t) + \varepsilon \phi (t).
\end{equation}

We have, by (\ref{eqV64}),
\begin{equation}\label{eqV66}
\begin{split}
E^\prime_\varepsilon(t) & \leq E^\prime (t) - \frac{\varepsilon}{4} \left [ E(t) \right ]^{\frac{p+1}{2\lambda}} 
\leq - \frac{\varepsilon}{4} \left [ E(t) \right ]^{\frac{p+1}{2\lambda}}.
\end{split}
\end{equation}

Also, if $ \varepsilon> 0 $ we have
\begin{equation}\label{eqV67}
\frac{1}{2} \left [ E_\varepsilon(t) \right ]^{\frac{p+1}{2\lambda}} \leq \left [ E(t) \right ]^{\frac{p+1}{2\lambda}} \leq 2 \left [ E_\varepsilon(t) \right ]^{\frac{p+1}{2\lambda}}.
\end{equation}

Combining (\ref{eqV66}) and (\ref{eqV67}) gives (\ref{eqV35}) as in case 1.b.

In this way the proof is completed for regular initial data. Estimates (\ref{eqV34}) and (\ref{eqV36}) extend to arbitrary initial data of finite energy $ \{ y^0, y^1 \} \in H^1_0(\Omega) \times L^2(\Omega) $ by density, since the constants of the estimates above continuously depend on $ E_y (0) $.
\end{proof}

\section{Comments} \label{sec:V.5}
\begin{enumerate}
\item {\bf Non-monotonic nonlinearities.} In the proof of Theorem \ref{theorem:V.5.2}, the fact that $g$ is non-decreasing is crucial since it allows us to ensure the pre-compactness of the trajectories for regular initial data.

It is not known if the conclusion of Theorem \ref{theorem:V.5.2} is true if $ g \in C(\mathbb{R}) $ satisfies
\begin{equation}\label{eqV68}
g(s) s > 0 \hspace{0.2cm} \text{ if } \hspace{0.2cm} s \neq 0.
\end{equation}
but it is not necessarily non-decreasing. The reason for this is that, in the absence of monotonicity of $g$, the higher order energy bound cannot be achieved. Thus, we cannot guarantee the relative compactness of trajectories in the energy space.

To handle this case, M. Slemrod \cite{slemrod1989weak} developed a "relaxed" version of the invariance principle showing  that when
\begin{equation}\label{eqV69}
    |g(s)| < C|s| \hspace{0.6cm} \forall s \in \mathbb{R}: |s|>1
\end{equation}
and  (\ref{eqV68}) is satisfied, then every trajectory of (\ref{eqV02}) satisfies, in the weak topology,
\begin{equation}\label{eqV70}
    \{ y(t), y^\prime(t) \} \rightharpoonup \{0,0\} \text{ weakly in } H^1_0(\Omega) \times L^2(\Omega) \text{ as }  {t \to +\infty}.
\end{equation}

\item {\bf Localized damping.} Theorem \ref{theorem:V.5.2} can be broadly generalized in several directions. In particular, we can consider wave equations of the form
\begin{equation}\label{eqV71}
    \begin{cases}
    y^{\prime\prime} - \Delta y + a(x) \beta (y^\prime) = 0 \hspace{0.2cm} & \text{in} \hspace{0.2cm} \Omega \times (0, \infty)
    \\
    y = 0 \hspace{0.2cm} & \text{on} \hspace{0.2cm} \partial \Omega \times (0, \infty)
    \\
    y(0) = y^0 \in H^1_0(\Omega); y^\prime (0) = y^1 \in L^2(\Omega) & \text{in} \hspace{0.2cm}  \Omega 

    \end{cases}
\end{equation}
with $ a = a(x) \in L^\infty (\Omega), a \geq 0, $ such that
\begin{equation}\label{eqV72}
a \geq a_0 > 0 \text{ in } \omega \subset \Omega
\end{equation}
where $\omega$ is a measurable subset with positive measure and $\beta$ is a non-decreasing but discontinuous function satisfying
\begin{equation}\label{eqV73}
    \beta (s) s> 0 \hspace{0.6cm} \forall s \in \mathbb{R}_0.
\end{equation}

A. Haraux \cite{haraux1985stabilization} showed that, under these conditions, the trajectories $ \{ y(t), y^\prime(t) \} $ of (\ref{eqV71}) strongly converge to an equilibrium $ \{ z, 0 \} $ in $ H^1_0(\Omega) \times L^2(\Omega) $, where $ z = z(x) $ is a solution of
\begin{equation}\label{eqV74}
    -\Delta z \in -a(x) \beta(0) \hspace{0.3cm} \text{ in } \hspace{0.3cm} \Omega; \hspace{0.3cm} z \in H^1_0(\Omega).
\end{equation}

The symbol $ \in $ in (\ref{eqV74}) indicates that $ \beta $ can be a multivalued function.

Another consequence of the results of \cite{haraux1985stabilization} is that, in order to have
$
    E_y(t) \to 0
$
for any finite energy solution of (\ref{eqV02}), it is sufficient that $ g \in C(\mathbb{R}) $ is non-decreasing and satisfies
$
    0 \notin \operatorname{int} \{ g^{-1}(0) \}.
$

Note that, under the hypothesis (\ref{eqV01}) of the theorem, $ \{ g^{-1}(0)\} = \{ 0 \} $,  implying $ \operatorname{int} \{ g^{-1}(0) \} = \emptyset. \hfill \square$

\item {\bf Optimality of the decay rates.} In Theorem \ref{theorem:V.5.3} we have shown upper bounds on the decay rates. These estimates are probably optimal but, to a large extent, this is an open subject.$ \hfill \square $

\item {\bf More general nonlinearities.} As mentioned in Remark \ref{remark:V.5.4}, the decay rate obtained depends solely on the behavior of $ g $ at the origin. However, we need the conditions (\ref{eqV31}) and (\ref{eqV33}) that limit the behavior of $ g $ at infinity as well.

It would be interesting to know if the estimates (\ref{eqV34}) - (\ref{eqV36}) are valid without hypotheses (\ref{eqV31}) and (\ref{eqV33}), that is, assuming only (\ref{eqV30}) and (\ref{eqV32}).

 We refer to \cite{haraux2023energy} for recent results on relaxing the growth conditions on the nonlinear damping g.

 The hypotheses (\ref{eqV30}) and (\ref{eqV32}) ensure that, near the origin, the graph of $ g $ is contained in the region limited by two power-like functions.

The method of proof of Theorem \ref{theorem:V.5.3} can be adopted to obtain decay rates when $g$ degenerates or grows faster than any function of the form $ |s|^{r-1} s $ with $ r > 0,$ \cite{begout2007interpolation}. The proof is inspired in the same principles but it requires to replace H{\"o}lder's inequality by Jensen's one to deal with a broader class of nonlinearities. 

\item {\bf Generalisation to an abstract setting.}  The results of this chapter generalize to abstract equations of the type
$
    y^{\prime \prime} +  Ay+ g(y^\prime) = 0
$
where $A$ is an elliptic operator of order greater than or equal than two (cf. \cite{haraux1987semi}, \cite{haraux1988decay} and  \cite{zuazua1988stability}). In particular, analogous results are obtained for the wave equation with Neumann boundary conditions
\begin{equation}\label{eqV75}
    \begin{cases}
    y^{\prime\prime} - \Delta y + g(y^\prime) = 0 \hspace{0.2cm} & \text{in} \hspace{0.2cm} \Omega \times (0, \infty)
    \\
    \frac{\partial y}{ \partial \nu} = 0 \hspace{0.2cm} & \text{on} \hspace{0.2cm} \partial \Omega \times (0, \infty)
    \\
    y(0) = y^0 \in H^1(\Omega), y^\prime (0) = y^1 \in L^2(\Omega) \hspace{0.2cm} & \text{in} \hspace{0.2cm} \Omega 
    \end{cases}
\end{equation}
and for the vibrating plate equation
\begin{equation}\label{eqV76}
    \begin{cases}
    y^{\prime\prime} - \Delta^2 y + g(y^\prime) = 0 \hspace{0.2cm} & \text{in} \hspace{0.2cm} \Omega \times (0, \infty)
    \\
    y = \frac{\partial y}{ \partial \nu} = 0 \hspace{0.2cm} & \text{on} \hspace{0.2cm} \partial \Omega \times (0, \infty)
    \\
    y(0) = y^0 \in H^2_0(\Omega); y^\prime (0) = y^1 \in L^2(\Omega) \hspace{0.2cm} & \text{in} \hspace{0.2cm} \Omega . 
    \end{cases}
\end{equation}
$ \hfill \square $

\item {\bf Contractivity estimates.} The methods in this chapter  allow also to obtain estimates for the decay rate of the distance between two different solutions of a non-autonomous equation of the form (cf. \cite{haraux1987semi}, \cite{haraux1988decay} and \cite{zuazua1988stability})
\begin{equation}\label{eqV77}
    \begin{cases}
    y^{\prime\prime} - \Delta y + g(y^\prime) = h(x,t) \hspace{0.2cm} & \text{in} \hspace{0.2cm} \Omega \times (0, \infty)
    \\
    y = 0 \hspace{0.2cm} & \text{on} \hspace{0.2cm} \partial \Omega \times (0, \infty)
    \\
    y(0) = y^0 \in H^1_0(\Omega); y^\prime (0) = y^1 \in L^2(\Omega) \hspace{0.2cm} & \text{in} \hspace{0.2cm} \Omega. 
    \end{cases}
\end{equation}
$\hfill \square $
 
\item {\bf Semilinear perturbations.} The same techniques allow obtaining estimates of the decay rate of solutions for semilinear equations of the form
\begin{equation}\label{eqV78}
    \begin{cases}
    y^{\prime\prime} - \Delta y + f(y) + g(y^\prime) = 0 \hspace{0.2cm} & \text{in} \hspace{0.2cm} \Omega \times (0, \infty)
    \\
    y = 0 \hspace{0.2cm} & \text{on} \hspace{0.2cm} \partial \Omega \times (0, \infty)
    \\
    y(0) = y^0; y^\prime (0) = y^1\hspace{0.2cm} & \text{in} \hspace{0.2cm} \Omega
    \end{cases}
\end{equation}
where $ f \in C^1(\mathbb{R}) $ is an increasing function satisfying certain growth conditions at infinity and such that $ f(0)=0 $ (cf. \cite{haraux1987semi}, \cite{zuazua1988stability}).
$ \hfill \square $

\item {\bf Localized damping.} The situation is considerably more complicated when we consider wave equations with localized dissipation
\begin{equation}\label{eqV79}
    \begin{cases}
    y^{\prime\prime} - \Delta y + f(y) + a(x) y^\prime = 0 \hspace{0.2cm} & \text{in} \hspace{0.2cm} \Omega \times (0, \infty)
    \\
    y = 0 \hspace{0.2cm} & \text{on} \hspace{0.2cm} \partial \Omega \times (0, \infty)
    \\
    y(0) = y^0; y^\prime (0) = y^1\hspace{0.2cm} & \text{in} \hspace{0.2cm} \Omega
    \end{cases}
\end{equation}
with $ a \in L^\infty_+ (\Omega) $ such that
\begin{equation}\label{eqV80}
    a(x) \geq a_0 > 0 \hspace{0.6cm} \forall x \in \omega
\end{equation}
where $ \omega $ is a nonempty open subset of $ \Omega $.

In the linear case $ f \equiv 0 $, it is shown that, if we have the estimate
\begin{equation}\label{eqV81}
    \| \phi^0 \|^2_{H^1_0(\Omega)} + \| \phi^1 \|^2_{L^2(\Omega)} \leq C \int^T_0 \int_\omega | \phi^\prime |^2 dx dt
\end{equation}
with $ T > 0 $ independent of the solution $ \phi = \phi (x,t) $ of the conservative equation
\begin{equation}\label{eqV82}
    \begin{cases}
    \phi^{\prime\prime} - \Delta \phi = 0 \hspace{0.2cm} & \text{in} \hspace{0.2cm} \Omega \times (0, \infty)
    \\
    \phi = 0 \hspace{0.2cm} & \text{on} \hspace{0.2cm} \partial \Omega \times (0, \infty)
    \\
    \phi(0) = \phi^0, \phi^\prime (0) = \phi^1\hspace{0.2cm} & \text{in} \hspace{0.2cm} \Omega,
    \end{cases}
\end{equation}
then the solutions of (\ref{eqV79}), under hypothesis (\ref{eqV80}), decay exponentially in $ H^1_0(\Omega) \times L^2(\Omega) $ (cf. \cite{haraux1989remarque}). We know, from Theorem \ref{theorem:III.3.1} of Chapter \ref{chapter03}, that (\ref{eqV81}) holds if $ \omega $ is a neighborhood of a subset of the boundary of the form $ \overline{\Gamma(x^0)} $. On the other hand, from \cite{bardos1992sharp}, we know that, if $ \Omega $ is of class $ C^\infty $, (\ref{eqV81}) holds if $ \omega $ satisfies the geometric control condition.

In the particular case in which $ \omega $ is a neighborhood of the entire boundary $ \Gamma $, these results are extended to the semilinear problem under suitable conditions on $f$ (cf.  \cite{zuazua1990exponential}). In  \cite{zuazua1991exponential} analogous results are shown for unbounded domains $ \Omega $ when $ \omega $ is the union of a neighborhood of infinity in $ \mathbb{R}^n $ and of a neighborhood of $ \partial \Omega $ in $ \Omega $.

J. Da Silva Ferreira in \cite{dasilva1994exponential}  obtained similar results for systems of semilinear wave equations with localized dissipation.

 Obtaining estimates of the decay rate of solutions of wave equations with localized nonlinear dissipation of the type
\begin{equation}\label{eqV85}
    \begin{cases}
    y^{\prime\prime} - \Delta y + \chi_\omega |y^\prime|^{p-1} y^\prime = 0 \hspace{0.2cm} & \text{in} \hspace{0.2cm} \Omega \times (0, \infty)
    \\
    y = 0 \hspace{0.2cm} & \text{on} \hspace{0.2cm} \partial \Omega \times (0, \infty)
    \\
    y(0) = y^0, y^\prime (0) = y^1 \hspace{0.2cm} & \text{in} \hspace{0.2cm} \Omega
    \end{cases}
\end{equation}
with $ p>1 $ and $ \omega $ a neighborhood of the boundary, requires significant further developments (see \cite{martinez1999method}, \cite{alabau2011sharp}, \cite{lasiecka2006energy} and \cite{tebou1998stabilization}).
$ \hfill \square $

\item {\bf  Thin support damping.} Let $ \omega_\varepsilon $ be a neighborhood of width $ \varepsilon > 0 $ of $ \partial \Omega $ and consider the wave equation  with localized dissipation in $ \omega_\varepsilon $:
\begin{equation}\label{eqV83}
    \begin{cases}
    y^{\prime\prime} - \Delta y + c_\varepsilon \chi_{\omega_\varepsilon} y^\prime = 0 \hspace{0.2cm} & \text{in} \hspace{0.2cm} \Omega \times (0, \infty)
    \\
    y = 0 \hspace{0.2cm} & \text{on} \hspace{0.2cm} \partial \Omega \times (0, \infty)
    \\
    y(0) = y^0; y^\prime (0) = y^1 \hspace{0.2cm} & \text{in} \hspace{0.2cm} \Omega
    \end{cases}
\end{equation}
with $ c_\varepsilon> 0 $.

It would be interesting to study the possibility of choosing $ c_\varepsilon $ so that (\ref{eqV83}) converges (in a precise sense) to a wave equation with dissipative boundary conditions, preserving the exponential decay.

This question is related to Comment 5 of Section \ref{sec:III.5}.
$ \hfill \square $

 \item {\bf Models involving non-local terms.} The techniques in this chapter also allow to study nonlinear wave equations of the form
\begin{equation}\label{eqV84}
    \begin{cases}
    y^{\prime\prime} - M (  \int_\Omega |\nabla y (x,t) |^2 dx) \Delta y + (-\Delta)^\alpha y^\prime = 0 \hspace{0.2cm} & \text{in} \hspace{0.2cm} \Omega \times (0, \infty)
    \\
    y = 0 \hspace{0.2cm} & \text{on} \hspace{0.2cm} \partial \Omega \times (0, \infty)
    \\
    y(0) = y^0; y^\prime (0) = y^1\hspace{0.2cm} & \text{in} \hspace{0.2cm} \Omega
    \end{cases}
\end{equation}
with $ \alpha  \geq 0 $ and $ M \in C(\mathbb{R}) $ such that $ M(s) \geq m_0 > 0$ for all $s \in \mathbb{R} $.

In L. A. Medeiros and M. Milla Miranda \cite{medeiros1990remarks} it is shown that, when $ \alpha > 0 $,  (\ref{eqV84}) has global, exponentially decaying solutions. $ \hfill \square $

\item {\bf Other references.} The literature on the problems we have discussed here and on related topics is very broad. Interested readers are referred for instance to \cite{alabau2012advances}, \cite{kafnemer2021weak}, \cite{martinez1999new} and \cite{martinez2000exponential} and the bibliographies therein. Among  pioneering works on this topic,  those by M. Nakao \cite{nakao1983energy} are also worth mentioning.

\end{enumerate}


%
%

\chapter{Boundary Stabilization of the Wave Equation}
\label{chapter06} 

\abstract{
In this chapter we study the boundary stabilization of the wave equation. We apply LaSalle’s Invariance Principle and show the decay of finite energy solutions under general conditions on the partition of the boundary and on the nonlinearity. We also prove explicit decay rates under suitable assumptions on the damping term. In addition, we present a second proof of the exponential decay for the linear dissipation, of interest when dealing with model perturbations.
}

\section{Problem Formulation}
\label{sec:VI.1}
In this chapter we study the rate of decay of the energy of solutions of wave equations with dissipative boundary conditions.

Let $ \Omega $ be a bounded domain of $ \mathbb{R}^n $ with boundary $ \Gamma $ of class $ C^2 $. Let $ \{ \Gamma_0, \Gamma_1 \} $ be a partition of the boundary $ \Gamma $ such that

\begin{equation}\label{eqVI01}
    \operatorname{int} (\Gamma_i) \neq \emptyset, \hspace{0.6cm} i = 0,1
\end{equation}
i.e. $ \Gamma_i $ contains a nonempty open subset of $ \Gamma $ for $ i=0,1 $.

Let $ g \in \mathbb{R} \to \mathbb{R} $ be a function such that
\begin{equation}\label{eqVI02}
  g \in C(\mathbb{R}); g(0)=0; g \text{ is nondecreasing}
\end{equation}
and
\begin{equation}\label{eqVI03}
  \alpha = \alpha (x) \in L^\infty_+ (\Gamma_0).
\end{equation}

We consider the wave equation
\begin{equation}\label{eqVI04}
    \begin{cases}
    y^{\prime\prime} - \Delta y = 0 \hspace{0.2cm} & \text{in} \hspace{0.2cm} \Omega \times (0, \infty)
    \\
    \frac{\partial y}{\partial \nu} = - \alpha (x) g(y^\prime) \hspace{0.2cm} & \text{on} \hspace{0.2cm} \Gamma_0 \times (0, \infty)
    \\ y = 0 \hspace{0.2cm} & \text{on} \hspace{0.2cm} \Gamma_1 \times (0, \infty)
    \\
    y(0) = y^0 \in V; y^\prime (0) = y^1 \in L^2(\Omega)& \text{in} \hspace{0.2cm} \Omega
    \end{cases}
\end{equation}
where
\begin{equation}\label{eqVI05}
  V = \{ v \in H^1(\Omega): v_{|\Gamma_1} = 0 \} =: H^1_{\Gamma_1} (\Omega).
\end{equation}

\begin{remark}
\label{remark:VI.6.1}
The condition (\ref{eqVI01}) for $ i=1 $ ensures that $ \| \nabla v \|_{L^2(\Omega)} $ defines a norm in $V$, equivalent to the one induced by $ H^1(\Omega) $.

On the other hand, if
\begin{equation}\label{eqVI06}
  g(s)s>0,  \, \forall s \in \mathbb{R}-\{0\}; \, 
  \alpha (x) \geq \alpha_0 > 0,  \, \text{for a.e.} \, x \in \Gamma^\prime_0,
\end{equation}
with $ \Gamma^\prime_0 $ an open and non-empty subset of $ \Gamma_0 $, we can ensure that the only equilibrium or stationary solution of the system is $ y \equiv 0 $.
$ \hfill \square $

We define the energy of the system
\begin{equation}\label{eqVI07}
E(t) = \frac{1}{2} \int_\Omega \left[ | \nabla y(x,t) |^2 + | y^\prime (x,t) |^2 \right] dx.
\end{equation}

Multiplying equation (\ref{eqVI04}) by $ y^\prime $ and integrating by parts we obtain (this argument is, by now, formal, but will be rigorously justified in the next section):
\begin{equation}\label{eqVI08}
\frac{dE(t)}{dt} = E^\prime(t) = - \int_{\Gamma_0} \alpha (x) g(y^\prime(x,t)) y^\prime (x,t) d\Gamma.
\end{equation}

From (\ref{eqVI02}), (\ref{eqVI03}) and (\ref{eqVI08}) it follows that $ E(t) $ is decreasing. Therefore, the boundary conditions of (\ref{eqVI04}) are dissipative.

Under hypothesis (\ref{eqVI06}) it can also be expected that
\begin{equation}\label{eqVI09}
E(t) \to 0 \hspace{0.3cm} \text{ when } \hspace{0.3cm} t \rightarrow \infty
\end{equation}
for any finite energy solution, or equivalently,
\begin{equation}\label{eqVI10}
\{ y(t), y^\prime(t) \} \to \{0, 0\} \hspace{0.3cm} \text{ in } \hspace{0.3cm} V \times L^2(\Omega) \hspace{0.3cm} \text{ when } \hspace{0.3cm} t \rightarrow +\infty.
\end{equation}

In Section \ref{sec:VI.3} we prove, using LaSalle's Invariance Principle, that, indeed, (\ref{eqVI09}) and (\ref{eqVI10}) are fulfilled under hypothesis (\ref{eqVI06}).

Once (\ref{eqVI09})-(\ref{eqVI10}) is shown, the following question arises naturally: Under suitable conditions on the partition $ \{ \Gamma_0, \Gamma_1 \} $ and the nonlinearity $ g $, can the rate of decay of the energy be estimated?

In particular, when the damping function $ g $  is linear, a uniform exponential decay of the energy is expected. However, by a classical argument due to D. L. Russell \cite{russell1978control}, we know that if the uniform exponential decay of the energy holds, with dissipation supported on $ \Gamma_0 $, then the wave equation is exactly controllable in the energy space with controls at $ \Gamma_0 \times (0,T) $ for $ T>0 $ large enough.
On the other hand, as described  in Chapter \ref{chapter02},  for the wave equation to be exactly controllable in the energy space with controls supported in $ \Gamma_0 \times (0,T) $, $ \Gamma_0 $ must be a sufficiently large part of $ \Gamma $ satisfying, roughly, the geometric control condition.

Thus, while (\ref{eqVI09})-(\ref{eqVI10}) can be expected to hold for a broad class of partitions $ \{ \Gamma_0, \Gamma_1 \} $, obtaining decay rates will require important geometric conditions on the partition $ \{ \Gamma_0, \Gamma_1 \} $.

Inspired by the exact controllability results of the previous chapters, obtained through multiplier techniques, the first case to consider is when
\begin{equation}\label{eqVI11}
\begin{cases}
  \Gamma_0 = \Gamma (x^0) = \{ x \in \Gamma: ( x - x^0 ) \cdot \nu(x) > 0 \}
  \\ \Gamma_1 = \Gamma_\star (x^0) = \Gamma \backslash \Gamma(x^0).
\end{cases}
\end{equation}

In Section \ref{sec:VI.4} we combine multiplier techniques and the methods of Chapter \ref{chapter05} to get decay rates for the partition (\ref{eqVI11}). As we will see, the choice of the potential $ \alpha (x) $ will play a fundamental role. Indeed, we will choose
\begin{equation}\label{eqVI12}
 \alpha (x) = ( x - x^0 ) \cdot \nu(x),
\end{equation}
which is a positive potential on $ \Gamma (x^0) $, vanishing at the interface points $ x \in \overline{\Gamma (x^0)} \cap \overline{\Gamma_\star(x^0)} $.

At first glance, this choice of $ \alpha $ seems to contradict the (erroneous) intuition that the higher the $ \alpha $, the larger the decay rate should be. Indeed, note that if $ \alpha = \text{constant} \to \infty $, the boundary conditions of (\ref{eqVI04}) tend, in a singular limit,  at least formally, towards the boundary conditions
\begin{equation*}
    g(y^\prime) = 0 \hspace{0.3cm} \text{ in } \hspace{0.3cm} \Gamma(x^0) \times (0, \infty); \hspace{0.3cm} y = 0 \hspace{0.3cm} \text{ in } \hspace{0.3cm} \Gamma_\star (x^0) \times (0, \infty).
\end{equation*}
In the particular case in which $ g(s)s >0 $ for all $s \neq 0 $ and $ y^0 \in H^1_0(\Omega) $, these conditions coincide with the homogeneous Dirichlet boundary conditions. In the limit system, therefore, the energy of the system is conserved and the dissipativity is lost. 

This is a manifestation of the classical overdamping phenomenon. Our intuition, according to which, the decay rate should increase when the damping term grows, is correct for small amplitude dampings, but fails beyond a critical threshold of the damping term.

Therefore, it is not surprising that the degeneration of the damping potential at the interface points $ \overline{\Gamma(x^0)} \cap \overline{\Gamma_\star(x^0)} $,  not being uniformly positive, might contribute to ensure the decay.

Due to the Dirichlet-Neumann mixed type boundary conditions, the solutions of (\ref{eqVI04}) develop singularities at the interface $ x \in \overline{\Gamma_0} \cap \overline{\Gamma_1} $ where the change of boundary conditions occurs, even if the initial data are regular. For this reason, it is necessary to be cautious in the application of multiplier techniques, since some of the integrations by parts that are required cannot be justified. P. Grisvard in \cite{grisvard1987controlabilite305} and \cite{grisvard1989controlabilite} made a detailed study of these singularities and proved that some of the identities that would be desirable in the application of multipliers are not true. However, in \cite{grisvard1987controlabilite305}, \cite{grisvard1989controlabilite} it is shown that for the partition (\ref{eqVI11}), certain inequalities (and not identities!) are satisfied, which is sufficient to obtain precise bounds. P. Grisvard's study is limited to dimensions $ n \leq 3 $ and therefore our results will also refer only to cases $ n = 1,2,3 $.

The constraint $ n \leq 3 $ is most likely purely technical, and the results in this chapter are probably valid in any dimension. In fact, the arguments that we use (except for the use of P. Grisvard's inequalities) do not depend on the dimension. But, so far, they are only fully justified if $ n \leq 3 $.

Of course, solutions do not develop singularities if there are no interface points, i.e. if:
\begin{equation}\label{eqVI13}
   \overline{\Gamma(x^0)} \cap \overline{\Gamma_\star(x^0)} = \emptyset.
\end{equation}

In this case the previous remarks on the difficulties generated by the singularities at the interface do not apply, and the analysis can be carried out in any dimension $ n \geq 1 $. However (\ref{eqVI13}) imposes very restrictive conditions on the geometry of the domain and is essentially satisfied if
$
    \Omega = \Omega_0 \backslash \overline{\Omega_1}
$
with $\Omega_0$ and $\Omega_1$ star-shaped domains with respect to $ x^0 $ such that $ \overline{\Omega_1} \subset \Omega_0 $.
$ \hfill \square $

It is convenient to make the following remark about condition (\ref{eqVI01}) for the partition (\ref{eqVI11}). The condition
$
    \operatorname{int} ( \Gamma_*(x^0) ) \neq \emptyset
$
must hold for the solution  to decay. It ensures, for instance, that Poincar\'e inequality, which plays a key role in the estimates leading to the decay, holds.  If $ \Gamma_\star(x^0) = \emptyset $ and $ \Gamma(x^0) = \Gamma $, the boundary conditions of (\ref{eqVI04}) are of the Neumann type
\begin{equation*}
    \frac{\partial y}{\partial \nu} = -\alpha (x) g (y^\prime) \hspace{0.3cm} \text{ in } \hspace{0.3cm} \Gamma \times (0, \infty).
\end{equation*}

In this case, every constant function is a solution of the system and, although it is of null energy, it does not converge to the null equilibrium. In this case $ V = H^1(\Omega) $ and the semi-norm $ \| \nabla \cdot \|_{L^2(\Omega)} $ does not define a norm in $ H^1(\Omega) $.

In this setting, in order to eliminate these non-trivial stationary solutions and to ensure the decay of any trajectory $ \{ y(t), y^\prime(t) \} $ in $ H^1(\Omega) \times L^2(\Omega) $ towards $ \{0, 0 \} $, it is convenient to consider boundary conditions of the type
\begin{equation*}
    \frac{\partial y}{\partial \nu} + \beta (x) y = -\alpha (x) g (y^\prime) \hspace{0.3cm} \text{ in } \hspace{0.3cm} \Gamma \times (0, \infty)
\end{equation*}
with $ \beta \in L^\infty_+ (\Gamma)$, non trivial.

This chapter is organized as follows. In Section \ref{sec:VI.2} we study the initial value problem (\ref{eqVI04}). In Section \ref{sec:VI.3} we apply LaSalle's Invariance Principle and show the decay (without estimates on the decay rate) of all finite-energy solutions under very general conditions on the partition of the boundary and the nonlinearity $ g $. In Section \ref{sec:VI.4} we obtain estimates of the energy decay rate for the partition (\ref{eqVI11}). In Section \ref{sec:VI.5} we give a second proof of the exponential decay when the dissipation is linear. We conclude this chapter with a section devoted to comment on some extensions and  open problems.
\end{remark}

\section{The Initial Value Problem}
\label{sec:VI.2}
Let $ V = H^1_{\Gamma_1} (\Omega) $ (as in (\ref{eqVI05})) and let us define the operator $ A \in L(V, V^\prime) $ by
\begin{equation}\label{eqVI14}
    \langle Av, w \rangle = \int_\Omega \nabla v(x) \cdot \nabla w(x) dx, \hspace{0.6cm} \forall v, w \in V.
\end{equation}

In (\ref{eqVI14}), $ \langle \cdot , \cdot \rangle $ denotes the duality product between $ V^\prime $ and $ V $.

Suppose that $ g \in C( \mathbb{R} ) $, in addition to  (\ref{eqVI02}), satisfies
\begin{equation}\label{eqVI15}
\exists C > 0: | g(s) | \leq C (1 + |s|), \hspace{0.6cm} \forall s \in \mathbb{R}.
\end{equation}

We define the following nonlinear operator $ B: V \to V^\prime $:
\begin{equation}\label{eqVI16}
\langle B v, w \rangle = \int_{\Gamma_0} \alpha (x) g ( v(x) ) w(x) d\Gamma \hspace{0.3cm} \forall v, w \in V.
\end{equation}

From (\ref{eqVI15}) it follows that $ v \rightarrow g(v) $ is continuous from $ L^2(\Gamma_0) \rightarrow L^2(\Gamma_0) $. On the other hand, since $ \alpha \in L^\infty (\Gamma_0) $ and the trace $ \gamma: H^1_{\Gamma_1} (\Omega) \rightarrow L^2(\Gamma_0) $ is continuous, $ Bv $ defines a single element of $ V^\prime $ for each $ v \in V $.
Since $ \alpha \geq 0 $ and $ g $ is nondecreasing, the operator $ B $ is monotone since
\begin{equation*}
\langle Bv - Bw, v - w \rangle = \int_{\Gamma_0} \alpha (x) ( g(v(x) ) - g( w(x))) (v(x) - w(x)) d\Gamma \geq 0 \hspace{0.3cm} \forall v, w \in V.
\end{equation*}

We introduce the vector variable
\begin{equation}\label{eqVI17}
    u(t) = 
    \begin{bmatrix}
    y(t)\\ y^\prime(t)
    \end{bmatrix}
\end{equation}
(which, to simplify the notation, we will simply write as $ u(t) = \{ y(t), y^\prime (t) \} $ and the matrix operator)
\begin{equation}\label{eqVI18}
    \mathcal{A} = 
    \begin{bmatrix}
    0 & -I \\ A & B
    \end{bmatrix}
\end{equation}
with domain
\begin{equation}\label{eqVI19}
\mathcal{D}(\mathcal{A}) = \{ \{v,w\} \in V \times V \mid Av + Bw \in L^2(\Omega) \}.
\end{equation}

With these definitions, the system (\ref{eqVI04}) can be written as
\begin{equation}\label{eqVI20}
 \begin{cases}
   u^\prime (t) + \mathcal{A} u(t) = 0 \hspace{0.6cm} \text{ for } \hspace{0.3cm} t \in (0, \infty)
   \\ u(0) = \{ y^0, y^1 \} \in V \times L^2(\Omega)
   \\ u(t) \in V \times L^2(\Omega) \hspace{0.6cm} \text{ for } \hspace{0.3cm} t \in (0, \infty).
 \end{cases}
\end{equation}

Let us prove that the operator $ (\mathcal{A}, \mathcal{D}(\mathcal{A})) $ is maximal monotone in $ V \times L^2(\Omega) $.

Given $ u_1 = \{ v_1, w_1 \}, u_2 = \{ v_2, w_2\} \in \mathcal{D}(\mathcal{A}) $ one has
\begin{equation}\label{eqVI21}
\begin{split}
& (\mathcal{A}u_1 - \mathcal{A}u_2,u_1 - u_2 )_{V \times L^2(\Omega)}
\\ & = ( \{ w_2 - w_1, \mathcal{A} v_1 - \mathcal{A} v_2 + Bw_1 - Bw_2 \}, \{ v_1 - v_2, w_1 - w_2 \} )_{V \times L^2(\Omega)}
\\ & = \int_\Omega \nabla (w_2 - w_1) \cdot \nabla (v_1 - v_2) dx + \langle Av_1 - Av_2, w_1 - w_2 \rangle
\\ & \, \, \, \, \,  + \langle Bw_1 - Bw_2, w_1 - w_2 \rangle = \langle B w_1 - B w_2, w_1 - w_2 \rangle \geq 0
\end{split}
\end{equation}
and therefore $ (\mathcal{A}, \mathcal{D}(\mathcal{A})) $ is monotone in $ V \times L^2(\Omega) $.

In (\ref{eqVI21}) $ (\cdot,\cdot)_{V \times L^2(\Omega)} $ denotes the scalar product in $ V \times L^2(\Omega) $
\begin{equation}\label{eqVI22}
(\{ v_1, w_1 \}, \{ v_2, w_2\})_{V \times L^2(\Omega)} = 
\int_\Omega \nabla v_1 \cdot \nabla v_2 dx + \int_\Omega w_1 w_2 dx,
\end{equation}
for all $ \{ v_1, w_1 \}, \{ v_2, w_2\} \in V \times L^2(\Omega).$

Let us see that $ R(I + \mathcal{A}) $, the range of $ (I + \mathcal{A})$, is the full space $V \times L^2(\Omega) $. To do it, given $ \{ f,g \} \in V \times L^2(\Omega) $, we consider the equation
\begin{equation}\label{eqVI23}
\{ v, w \} + \mathcal{A} \{ v, w \} = \{ f, g \},
\end{equation}
that is
\begin{equation}\label{eqVI24}
v - w = f \in V; \, w + Av + Bw = g \in L^2(\Omega).
\end{equation}

Writing $ v = f + w $ in the second equation we get
\begin{equation}\label{eqVI25}
Aw + Bw + w = g - Af \in V^\prime.
\end{equation}

On the other hand, $A$ is strictly coercive in $V$, since
\begin{equation}\label{eqVI26}
 \langle Av, v \rangle = \int_\Omega | \nabla v |^2 dx
\end{equation}
and $ \langle Bw, w \rangle \geq 0 $ in $ V $. Consequently $ ( A + B + I )^{-1} $, which is well defined in $ R( A + B + I ) $, maps  bounded sets of $ V^\prime $ in bounded sets of $V$. Since $A+B$ is monotone and continuous from $V$ to $V^\prime$, it follows that  $ R( A + B + I ) = V^\prime $ (see F. Browder [6], Chapter I).

If $ w \in V $ is the only solution to (\ref{eqVI25}), then $ \{ w + f, w \} \in V \times V $ is the only solution of (\ref{eqVI23}). On the other hand $ Bw + Av = g - w \in L^2(\Omega)$ (with  $v = w + f$) and therefore $ \{ v, w \} \in \mathcal{D}(\mathcal{A}) $.

Since $ \mathcal{A} $ is maximal monotone,  it generates a continuous semigroup of nonlinear contractions in $ \overline{\mathcal{D}(\mathcal{A})} $. On the other hand, since $ g(0)=0 $,
\begin{equation}\label{eqVI27}
 \{ v \in V:A v \in L^2(\Omega), \frac{\partial v}{\partial v} \bigg|_{\Gamma_0} = 0 \} \times V \subset \mathcal{D}(\mathcal{A}).
\end{equation}

We deduce that $ \mathcal{D}(\mathcal{A}) $ is dense in $ V \times L^2(\Omega)$.

Applying the results of H. Brezis \cite{brezis1983analyse}, Chapter III, we obtain (see also A. Haraux \cite{haraux1981nonlinear}, p. 70-79) the following result on the  existence, uniqueness and stability of solutions.

\begin{theorem}
\label{theorem:VI.6.1}
Suppose that (\ref{eqVI01})-(\ref{eqVI03}) and (\ref{eqVI15}) are satisfied.  Then the following results hold
\label{sec:VI.2-a2}
\begin{enumerate}[label=\roman*)]
\item Strong solutions.
If $ \{ y^0, y^1 \} \in \mathcal{D}(\mathcal{A}) $, there exists a unique strong solution of (\ref{eqVI04}) in the class
\begin{equation}\label{eqVI28}
    \begin{cases}
      y \in W^{1,\infty} (0, \infty; V) \cap W^{2, \infty} (0, \infty; L^2(\Omega))
      \\ Ay(t) + B y^\prime(t) \in L^2(\Omega), \hspace{0.3cm} \forall t>0
    \end{cases}
\end{equation}
and the associated energy satisfies
\begin{equation}\label{eqVI29}
\frac{dE(t)}{dt} = E^\prime(t) = - \int_{\Gamma_0} \alpha (x) g( y^\prime (x,t) ) y^\prime(x,t) d\Gamma, \hspace{0.3cm} \forall t > 0.
\end{equation}

\label{sec:VI.2-b}
\item Weak solutions.
If $ \{ y^0, y^1 \} \in V \times L^2(\Omega) $, the system (\ref{eqVI04}) admits a unique weak finite energy solution such that
\begin{equation}\label{eqVI30}
    y \in C ([0, \infty); V) \cap C^1 ([0, \infty); L^2(\Omega)).
\end{equation}

\label{sec:VI.2-c}
\item Stability. 
Given two weak solutions $y$, $\hat{y}$ we have
\begin{equation}\label{eqVI31}
\begin{split}
\frac{1}{2} \int_\Omega \left[ | \nabla y (x,t) - \nabla \hat{y} (x,t) |^2 + | y^\prime (x,t) - \hat{y}^\prime (x,t) |^2 \right] dx
\\ \leq \frac{1}{2} \int_\Omega \left[ | \nabla y^0 (x,t) - \nabla \hat{y}^0 (x) |^2 + | y^1 (x) - \hat{y}^1 (x) |^2 \right] dx.
\end{split}
\end{equation}
\end{enumerate}
\end{theorem}

\begin{remark}
\label{remark:VI.6.2}
Existence, uniqueness and stability results can be proved when $g$ is a multivalued function with superlinear growth at infinity (cf. I. Lasiecka \cite{lasiecka1989stabilization}, \cite{lasiecka1988stabilization}).
$ \hfill \square $ 
\end{remark}

\section{Consequences of LaSalle's Invariance Principle}
\label{sec:VI.3}
We have the following result:

\begin{theorem}
\label{theorem:VI.6.2}
Suppose that conditions (\ref{eqVI01})-(\ref{eqVI03}), and (\ref{eqVI15}) are satisfied.
Suppose furthermore that
\begin{equation}\label{eqVI32}
 \alpha (x) \geq \alpha_0 > 0 \hspace{0.3cm} \text{ for a.e. } \hspace{0.3cm} x \in \Gamma^\prime_0 \text{ non-empty open subset of } \Gamma_0,
\end{equation}
\begin{equation}\label{eqVI33}
g(s) s > 0 \hspace{0.6cm} \forall s \in \mathbb{R}.
\end{equation}

Then, for all $ \{ y^0, y^1 \} \in V \times L^2(\Omega) $ the corresponding weak solution of (\ref{eqVI04}) satisfies
\begin{equation}\label{eqVI34}
E(t) \to 0 \hspace{0.3cm} \text{ when } \hspace{0.3cm} t \to \infty.
\end{equation}
\end{theorem}

\begin{remark}
\label{remark:VI.6.3}
Conditions (\ref{eqVI32}) and (\ref{eqVI33}) ensure that the dissipation is effective in $ \Gamma_0^\prime $ whatever the value of $ y^\prime \neq 0 $.
$ \hfill \square $
\end{remark}

\begin{proof}
\label{proof:VI.theorem.6.3}
As we observed in the proof of Theorem \ref{theorem:V.5.2}  in Chapter \ref{chapter05}, thanks to the stability property (\ref{eqVI31}) and the density of $ \mathcal{D}(\mathcal{A})$ in $ V \times L^2(\Omega) $, it is enough to consider the case of the initial data $ \{ y^0, y^1 \} \in \mathcal{D}(\mathcal{A}) $.

Combining (\ref{eqVI28}) and equation (\ref{eqVI04}) we deduce that
\begin{equation}\label{eqVI35}
\{ y(t), y^\prime(t) \}_{t \geq 0} \hspace{0.3cm} \text{ is relatively compact in } \hspace{0.3cm} V \times L^2(\Omega).
\end{equation}

Therefore, it is sufficient to prove that the only point of accumulation of the trajectory in $ V \times L^2(\Omega) $ is $ \{ 0,0 \} $.

Let $ t_n \rightarrow \infty $ be such that
\begin{equation}\label{eqVI36}
\{ y(t_n), y^\prime(t_n) \} \to \{ z^0, z^1 \} \hspace{0.3cm} \text{ in } \hspace{0.3cm} V \times L^2(\Omega).
\end{equation}


We define the translated solutions
\begin{equation}\label{eqVI37}
z_n (x,t) = y (x, t + t_n).
\end{equation}

Given any $ T>0 $
\begin{equation}\label{eqVI38}
\{ z_n \} \hspace{0.1cm} \text{ is bounded in } \hspace{0.1cm} W^{1, \infty} (0,T;V) \cap W^{2, \infty} (0,T;L^2(\Omega))
\end{equation}
and it is also a solution of (\ref{eqVI04}) for $ t \in [0,T] $. Using (\ref{eqVI04}) and (\ref{eqVI28}) we deduce that
\begin{equation}\label{eqVI39}
\{ z_n \} \hspace{0.1cm} \text{ is relatively compact in } \hspace{0.1cm} C( [0,T];V) \cap C^1 ([0,T];L^2(\Omega)).
\end{equation}

Thus, extracting a subsequence (which we will continue to denote by $ \{ z_n \} $) we have
\begin{equation}\label{eqVI40}
z_n \to z \hspace{0.3cm} \text{ in } \hspace{0.3cm} C( [0,T];V) \cap C^1 ([0,T];L^2(\Omega))
\end{equation}
with
\begin{equation}\label{eqVI41}
z(0) = z^0, z^\prime(0) = z^1.
\end{equation}

On the other hand, we have
\begin{equation}\label{eqVI42}
z^{\prime\prime} - \Delta z = 0 \hspace{0.2cm} \text{ in } \hspace{0.2cm} \Omega \times (0,T).
\end{equation}

From (\ref{eqVI29}) we deduce that the energy $E(t)$ associated with $ y = y(x,t) $ is decreasing and, accordingly, the following limit exists
\begin{equation}\label{eqVI43}
\exists \lim_{t \rightarrow \infty} E(t) = \ell.
\end{equation}

Now since
\begin{equation}\label{eqVI44}
\begin{split}
    \int^T_0 \int_{\Gamma_0} \alpha (x) g (z_n^\prime (x,t)) z_n^\prime (x,t) d\Gamma dt
    & = \int^{t_n + T}_{t_n} \int_{\Gamma_0} \alpha (x) g (y^\prime) y^\prime d\Gamma dt 
    \\ &= E(t_n) - E(t_n + T)
\end{split}
\end{equation}
we deduce that
\begin{equation}\label{eqVI45}
    \lim_{n \rightarrow \infty}
    \int^T_0 \int_{\Gamma_0} \alpha (x) g (z_n^\prime (x,t)) z_n^\prime (x,t) d\Gamma dt = 0.
\end{equation}

From (\ref{eqVI45}) it follows that
\begin{equation*}
    \alpha (x) g (z^\prime (x,t)) z^\prime (x,t) = 0 \hspace{0.3cm} \text{ in } \hspace{0.3cm} \Gamma_0 \times (0,T) 
\end{equation*} 
and by (\ref{eqVI32})-(\ref{eqVI33}) we deduce
\begin{equation}\label{eqVI46}
    z^\prime (x,t) = 0 \hspace{0.6cm} \forall (x,t) \in \Gamma_0^\prime \times (0,T).
\end{equation}

But then, passing to the limit in the boundary condition satisfied by $ z_n $, we also get
\begin{equation}\label{eqVI47}
    \frac{\partial z}{\partial v} = 0 \hspace{0.3cm} \text{ in } \hspace{0.3cm} \Gamma_0 \times (0,T).
\end{equation}

Therefore, the limit $ z = z(x,t) $ satisfies
\begin{equation}\label{eqVI48}
    \begin{cases}
      z^{\prime\prime} - \Delta z = 0 & \text{ in } \hspace{0.3cm} \Omega \times (0,T)
      \\ \partial z/\partial v = 0 \hspace{0.3cm} & \text{ in } \hspace{0.3cm} \Gamma_0 \times (0,T)\\
       z=0 \hspace{0.3cm} &\text{ in } \hspace{0.3cm} \Gamma_1 \times (0,T)
      \\ z^\prime = 0 \hspace{0.3cm} & \text{ in } \hspace{0.3cm} \Gamma_0^\prime \times (0,T)
      \\ z(0) = z^0, z^\prime(0) = z^1\hspace{0.3cm} & \text{ in } \hspace{0.3cm} \Omega.
    \end{cases}
\end{equation}

Then the function $ \eta = z^\prime $ satisfies
\begin{equation}\label{eqVI49}
    \begin{cases}
      \eta^{\prime\prime} - \Delta \eta = 0 & \text{ in } \hspace{0.3cm} \Omega \times (0,T)
      \\ \frac{\partial \eta}{\partial v} = \eta = 0 \hspace{0.3cm} & \text{ in } \hspace{0.3cm} \Gamma_0^\prime \times (0,T).
    \end{cases}
\end{equation}

But then, by Holmgren's Theorem, we know that, if $ T > 0 $ is large enough (larger than $T_0$ = twice  the maximal geodesic distance from $\Gamma_0^\prime$ to all other points in $\Omega$), from (\ref{eqVI49}) it follows (with $ 2 \varepsilon = T - T_0$)
\begin{equation*}
    \eta \equiv 0 \hspace{0.3cm} \text{ in } \hspace{0.3cm} \Omega \times \left( \frac{T}{2} - \varepsilon, \frac{T}{2} + \varepsilon \right)
\end{equation*}
and therefore
\begin{equation*}
    z = z(x) \hspace{0.3cm} \text{ in } \hspace{0.3cm} \Omega \times \left( \frac{T}{2} - \varepsilon, \frac{T}{2} + \varepsilon \right).
\end{equation*}

Then, $ z $ satisfies
\begin{equation*}
    \begin{cases}
      -\Delta z = 0 & \text{ in } \hspace{0.3cm} \Omega \times \left( \frac{T}{2} - \varepsilon, \frac{T}{2} + \varepsilon \right)
      \\ z = 0 & \text{ in } \hspace{0.3cm} \Gamma_1
      \\ \partial z/\partial v = 0 \hspace{0.3cm} & \text{ in } \hspace{0.3cm} \Gamma_0 \times \left( \frac{T}{2} - \varepsilon, \frac{T}{2} + \varepsilon \right)
    \end{cases}
\end{equation*}
and we deduce that $ z \equiv 0 $ in $ \Omega \times (T/2 - \varepsilon, T/2 + \varepsilon) $ and therefore $ z \equiv 0 $ in $ \Omega \times (0,T) $. In particular $ z^0 = z^1 = 0 $.

Theorem \ref{theorem:VI.6.2} has been proven.
\end{proof}

\begin{remark}
\label{remark:VI.6.4}
In G. Chen and H. K. Wang \cite{wang1989asymptotic} and I. Lasiecka \cite{lasiecka1988stabilization} more general results are shown for multivalued functions $g$ growing superlinearly at infinity.
Here we did not pursue to get the optimal growth conditions on $g$, since in the next section we will require that $g$ is sublinear at infinity to get explicit decay rates.
$ \hfill \square $
\end{remark}

\section{Decay Rates}
\label{sec:VI.4}
Let $ \Omega $ be a bounded domain of $ \mathbb{R}^n $, $n\ge 1$,  of class $C^2$. Given a point $ x^0 \in \mathbb{R}^n $ such that
\begin{equation}\label{eqVI50}
    \operatorname{int } ( \Gamma_\star (x^0)) \neq \emptyset
\end{equation}
we define
\begin{equation}\label{eqVI51}
    \alpha (x) = ( x - x^0 ) \cdot \nu(x)
\end{equation}
where $\alpha(x)$ is the outward unit normal vector at $x \in \partial \Omega$, and consider the system
\begin{equation}\label{eqVI52}
    \begin{cases}
      y^{\prime\prime} -\Delta y = 0 & \text{ in } \hspace{0.3cm} \Omega \times (0,\infty)
      \\ y = 0 & \text{ in } \hspace{0.3cm} \Gamma_\star (x^0) \times (0,\infty)
      \\ \partial y/\partial \nu = - \alpha (x) g (y^\prime) \hspace{0.3cm} & \text{ in } \hspace{0.3cm} \Gamma (x^0) \times (0,\infty)
      \\ y(0) = y^0 \in V, y^\prime(0) = y^1 \in L^2(\Omega) &\text{ in } \hspace{0.3cm} \Omega
    \end{cases}
\end{equation}
with $ V = \{ v \in H^1(\Omega): v_{| \Gamma_\star(x^0)} = 0 \} $.
The following holds.

\begin{theorem}
\label{theorem:VI.6.3}
Suppose the space dimension $ n \leq 3 $, that (\ref{eqVI50}) holds, and $ g \in C(\mathbb{R}) $ is a non-decreasing function satisfying
\begin{equation}\label{eqVI53}
 \exists \lambda > 0, C > 0: |g(s)| \leq C|s|^\lambda, \, \forall s \in [-1, 1]
\end{equation}
\begin{equation}\label{eqVI54}
 \exists C > 0:|g(s)| \leq C|s|, \,  \forall s \in \mathbb{R}: |s| \geq 1
\end{equation}
\begin{equation}\label{eqVI55}
 \exists p>0, c > 0: g(s)s \geq c|s|^{p+1}, \, \forall s \in [-1,1]
\end{equation}
\begin{equation}\label{eqVI56}
 \exists c > 0: g(s)s \geq c|s|^2, \,  \forall s \in \mathbb{R}: |s| \geq 1.
\end{equation}

Then,
\begin{enumerate}[label=\roman*)]
\item If $ \lambda = p = 1 $, there are constants $ M > 1, \gamma > 0 $, independent of the initial data, such that
\begin{equation}\label{eqVI57}
 E(t) \leq ME(0)e^{-\gamma t} \hspace{0.6cm} \forall t \geq 0
\end{equation}
for all solutions of (\ref{eqVI52}).

\item If $ \lambda < 1 $, there is a constant $ M $ depending on $ E(0) $ such that
\begin{equation}\label{eqVI58}
 E(t) \leq 4 \left[ Mt + (E(0))^{-(p + 1 - 2\lambda)/2 \lambda} \right]^{-2\lambda/(p+1-2\lambda)}  \hspace{0.6cm} \forall t \geq 0.
\end{equation}

\item If $ \lambda\geq 1 $ and $ p > 1 $, there is a constant $ M>0 $ depending on $ E(0) $ such that
\begin{equation}\label{eqVI59}
 E(t) \leq 4 \left[ Mt + (E(0))^{-(p - 1)/2} \right]^{-2/(p-1)}  \hspace{0.6cm} \forall t \geq 0.
\end{equation}
\end{enumerate}
\end{theorem}

\begin{remark}
\label{remark:VI.6.5}
The decay rates are analogous to those obtained in Theorem \ref{theorem:V.5.3} of Chapter \ref{chapter05}. Note that conditions (\ref{eqVI54}) and (\ref{eqVI56}) at infinity impose restrictions on the growth of $g$ from below and above. However, the order of decay of the energy depends solely on the properties of $g$ at the origin.
$ \hfill \square $
\end{remark}

\begin{remark}
\label{remark:VI.6.6}
The restriction $ n \leq 3 $ on the dimension is due, as mentioned in Section \ref{sec:VI.1}, to the singularities that the solutions of (\ref{eqVI52}) develop at the interface points $ x \in \overline{\Gamma(x^0)} \cap \overline{\Gamma_\star(x^0)} $. The result is probably true for any space dimension $n\ge 1$.

Of course, in the particular case where
$
   \overline{\Gamma(x^0)} \cap \overline{\Gamma_\star(x^0)} = \emptyset,
$
the solutions do not present singularities and the previous result holds without any restriction on the dimension $n\ge 1$.
$ \hfill \square $
\end{remark}

\begin{remark}
\label{remark:VI.6.7}
This result was proved in  \cite{zuazua1990stabilization} for $ \lambda \geq 1 $. The method  of proof when $ \lambda < 1 $ is inspired by the energy perturbation arguments in \cite{lagnese1991uniform}.
$ \hfill \square $
\end{remark}

\begin{proof}
\label{proof:VI.remark.6.7}
Thanks to the stability property (\ref{eqVI31}) and to the fact that all the constants  in the following estimates  depend continuously on the energy $ E(0) $, it is enough for us to consider initial data $\{ y^0, y^1 \} \in \mathcal{D}(\mathcal{A}).$

We consider the function
\begin{equation}\label{eqVI60}
    \rho (t) = 2 \int_\Omega y^\prime (x,t) m(x) \cdot \nabla y (x,t) dx + (n-1) \int_\Omega y^\prime (x,t) y (x,t) dx
\end{equation}
with $ m(x) = x - x^0 $.

Note that
\begin{equation}\label{eqVI61}
    | \rho (t) | \leq C_1 E(t)
\end{equation}
with $ C_1 = 2R + (n-1) \mu $ where $ R = \| x - x^0 \|_{L^\infty (\Omega)} $ and $ \mu > 0 $ is the smallest constant such that
\begin{equation}\label{eqVI62}
    \| v \|_{L^2(\Omega)} \leq \mu \| \nabla v \|_{L^2(\Omega)} \hspace{0.6cm} \forall v \in V.
\end{equation}

From now on,  to simplify the notation, we will omit the variables $x$ and $t$ in the integration signs.

We compute the derivative of $ \rho $ and obtain
\begin{equation}\label{eqVI63}
    \begin{split}
        \frac{d \rho (t)}{dt} = \rho^{\prime}(t) & = 2 \int_\Omega y^{\prime\prime} m \cdot \nabla y dx + 2 \int_\Omega y^\prime m \cdot \nabla y^\prime dx
        \\ & \quad + (n-1) \int_\Omega |y^\prime|^2 dx + (n-1) \int_\Omega y^{\prime\prime} y dx
        \\ & = 2 \int_\Omega \Delta ym \cdot \nabla y dx + 2 \int_\Omega m \cdot \nabla (|y^\prime|^2) dx
        \\ &  \quad+ (n-1) \int_\Omega |y^\prime|^2 dx + (n-1) \int_\Omega \Delta y y dx.
    \end{split}
\end{equation}

We have
\begin{equation}\label{eqVI64}
    \begin{split}
    \int_\Omega m \cdot \nabla(|y^\prime|^2) dx & = -n \int_\Omega |y^\prime|^2 dx + \int_\Gamma (m \cdot \nu) |y^\prime|^2 d\Gamma
    \\ & = -n \int_\Omega |y^\prime|^2 dx + \int_{\Gamma (x^0)} (m \cdot \nu) |y^\prime|^2 d\Gamma
    \end{split}
\end{equation}
and
\begin{equation}\label{eqVI65}
    \begin{split}
    \int_\Omega \Delta y ydx = - \int_\Omega |\nabla y|^2 dx + \int_\Gamma \frac{\partial y}{\partial \nu} y d\Gamma = - \int_\Omega | \nabla y |^2 dx - \int_{\Gamma (x^0)} \alpha (x) g(y^\prime) y d\Gamma.
    \end{split}
\end{equation}

We now need the following technical result, derived from those in \cite{grisvard1987controlabilite305} and \cite{grisvard1989controlabilite}.

\begin{lemma}
\label{lemma:VI.6.1}
For every $ y \in V $ such that $ \Delta y \in L^2(\Omega) $ and $ \partial y/ \partial \nu = (m \cdot \nu) v $ with $ v \in L^2(\Gamma (x^0)) $ the following inequality is satisfied
\begin{equation}\label{eqVI66}
    2 \int_\Omega \Delta y m \cdot \nabla y dx \leq (n - 2) \int_\Omega |\nabla y|^2 dx + R^2 \int_{\Gamma(x^0)} 
    (m \cdot \nu) |v|^2 d\Gamma.
\end{equation}
\end{lemma}

\begin{proof}
\label{proof:VI.lemma.6.1}
 In \cite{grisvard1987controlabilite}, \cite{grisvard1989controlabilite} it was proved that the following inequality holds when  $ v=0 $:
\begin{equation}\label{eqVI67}
\begin{split}
    2 \int_\Omega \Delta y m \cdot \nabla y dx & \leq (n - 2) \int_\Omega |\nabla y|^2 dx 
    \\ &  \quad+ \int_{\Gamma_\star(x^0)} (m \cdot \nu) \left| \frac{\partial y}{\partial \nu} \right|^2 d\Gamma
     - \int_{\Gamma_(X^0)} (m \cdot \nu) \left| \nabla_\sigma y \right|^2 d\Gamma
\end{split}
\end{equation}
when $ n \leq 3 $, where $ \nabla_\sigma $ denotes the tangential gradient.

Now since $ m \cdot \nu \geq 0 $ on $ \Gamma(x^0) $ and $ m \cdot \nu \leq 0 $ or $ \Gamma_\star(x^0) $, we get
\begin{equation}\label{eqVI68}
    2 \int_\Omega \Delta y m \cdot \nabla y dx \leq (n - 2) \int_\Omega |\nabla y|^2 dx.
\end{equation}

Next, in V. Komornik and E. Zuazua \cite{komornik1990direct} it was proved that, if $ v \in H^{1/2} (\Gamma (x^0)) $, the following holds
\begin{equation}\label{eqVI69}
\begin{split}
    2 \int_\Omega \Delta y m \cdot \nabla y dx
     \leq (n - 2) \int_\Omega |\nabla y|^2 dx 
+ 2 \int_\Gamma \frac{\partial y}{\partial \nu} m \cdot \nabla y d\Gamma - \int_\Gamma (m \cdot \nu) |\nabla y|^2 d\Gamma.
\end{split}
\end{equation}

On the other hand,
\begin{equation}\label{eqVI70}
\begin{split}
&2 \int_\Gamma \frac{\partial y}{\partial v} m \cdot \nabla y d\Gamma - \int_\Gamma (m \cdot \nu) |\nabla y|^2 d\Gamma
\\& = 2 \int_{\Gamma(x^0)} \frac{\partial y}{\partial \nu} m \cdot \nabla y d\Gamma
 - \int_{\Gamma(x^0)} (m \cdot \nu) |\nabla y|^2 d\Gamma 
 + \int_{\Gamma_{\star}(X^0)} (m \cdot \nu) |\nabla y|^2 d\Gamma
\\& \leq 2 \int_{\Gamma(x^0)} \frac{\partial y}{\partial \nu} m \cdot \nabla y d\Gamma - \int_{\Gamma(x^0)} (m \cdot \nu) |\nabla y|^2 d\Gamma
\\& \leq R^2 \int_{\Gamma(x^0)} \frac{1}{(m \cdot \nu)} \left| \frac{\partial y}{\partial \nu} \right|^2 d\Gamma = R^2 \int_{\Gamma(x^0)} (m \cdot \nu) |v|^2 d\Gamma.
\end{split}
\end{equation}
Combining (\ref{eqVI69}) and (\ref{eqVI70}) we get (\ref{eqVI66}).

Finally, by density, inequality (\ref{eqVI66}) is extended to the case where $ v \in L^2(\Gamma (x^0)) $. \quad \quad
\end{proof}

Combining (\ref{eqVI63}), (\ref{eqVI64}), (\ref{eqVI65}) and (\ref{eqVI66}) we get
\begin{equation}\label{eqVI71}
    \rho^{\prime}(t) \leq -2E(t) + \int_{\Gamma(x^0)} (m \cdot \nu) \left[ R^2 |g(y^\prime)|^2 - (n-1) g(y^\prime) y + |y^\prime|^2 \right] d\Gamma.
\end{equation}
At this level we distinguish the cases $ \lambda \geq 1 $ and $ \lambda < 1 $.
\smallskip

{ \bf $ \text{Case 1. } \lambda \geq 1.$ }
\label{case.1:proof:VI.lemma.6.1}
Let $ \beta >0 $ be the largest constant such that
\begin{equation}\label{eqVI72}
    \beta \int_{\Gamma(x^0)} (m \cdot \nu) |v|^2 d\Gamma \leq \int_\Omega |\nabla v|^2 dx \hspace{0.6cm} \forall v \in V.
\end{equation}

We have
\begin{equation*}
    (n-1) \left| \int_{\Gamma(x^0)} (m \cdot \nu) g(y^\prime) y d\Gamma \right| \leq E(t) + \frac{(n-1)^2}{2\beta} \int_{\Gamma(x^0)} (m \cdot \nu) |g(y^\prime)|^2 d\Gamma
\end{equation*}
and therefore, by (\ref{eqVI53})-(\ref{eqVI54}),
\begin{equation}\label{eqVI73}
    (n-1) \left| \int_{\Gamma(x^0)} (m \cdot \nu) g(y^\prime) y d\Gamma \right| \leq E(t) + C \int_{\Gamma(x^0)} (m \cdot \nu) |y^\prime|^2 d\Gamma.
\end{equation}

Combining (\ref{eqVI71}) and (\ref{eqVI73}) we get
\begin{equation}\label{eqVI74}
    \rho^\prime(t) \leq -E(t) + C \int_{\Gamma(x^0)} (m \cdot \nu) |y^\prime|^2 d\Gamma.
\end{equation}

Let
\begin{equation}\label{eqVI75}
    \phi(t) = \left[ E(t) \right]^{(p-1)/2} \rho (t).
\end{equation}

We have
\begin{equation}\label{eqVI76}
\begin{split}
    \phi^\prime(t) & = \frac{p-1}{2} \left[ E(t) \right]^{\frac{p-3}{2}} \rho (t) E^\prime(t) + \left[ E(t) \right]^{\frac{p-1}{2}} \rho^\prime (t)
    \\ & \leq -C_2 E^\prime(t) - \left[ E(t) \right]^{\frac{p+1}{2}}  + C \left[ E(t) \right]^{\frac{p-1}{2}} \int_{\Gamma(x^0)} (m \cdot \nu) |y^\prime|^2d\Gamma
\end{split}
\end{equation}
with $ C_2 = (p-1) C (E(0))^{(p-1)/2}/2.$

Given $ \varepsilon > 0 $ we define the functional
\begin{equation}\label{eqVI77}
    E_\varepsilon (t) = (1 + \varepsilon C_2) E(t) + \varepsilon \left[ E(t) \right]^{(p-1)/2} \rho (t).
\end{equation}

Since $ E(t) \leq E(0) $, taking $ \varepsilon \leq \varepsilon_0 (E_0) $ we have
\begin{equation*}
    \frac{1}{2} \left[ E_\varepsilon(t) \right]^{\frac{(p+1)}{2}} \leq \left[ E(t) \right]^{\frac{(p+1)}{2}} \leq 2 \left[ E_\varepsilon(t) \right]^{\frac{(p+1)}{2}}.
\end{equation*}

From (\ref{eqVI29}) and (\ref{eqVI76}) we obtain
\begin{equation}\label{eqVI78}
\begin{split}
    E_\varepsilon^\prime(t) \leq & - \int_{\Gamma(x^0)} (m \cdot \nu) g (y^\prime) y^\prime d\Gamma
    \\ & - \varepsilon \left[ E_\varepsilon(t) \right]^{\frac{(p+1)}{2}} + \varepsilon C \left[ E_\varepsilon(t) \right]^{\frac{(p-1)}{2}} \int_{\Gamma(x^0)} (m \cdot \nu) |y^\prime|^2 d\Gamma.
\end{split}
\end{equation}

From (\ref{eqVI56}) it follows that, if $ \varepsilon \leq \varepsilon_1 (E(0)) $,
\begin{equation}\label{eqVI79}
\begin{split}
\varepsilon C \left[ E(t) \right]^{\frac{(p-1)}{2}} \int_{\Gamma(x^0) \cap \{ |y^\prime| \geq 1 \} } (m \cdot \nu) |y^\prime|^2 d\Gamma
 \leq \int_{\Gamma(x^0) \cap \{ |y^\prime| \geq 1 \} } (m \cdot \nu) g(y^\prime) y^\prime d\Gamma.
\end{split}
\end{equation}

On the other hand, if $ p > 1 $,
\begin{equation}\label{eqVI80}
\begin{split}
&C \left[ E(t) \right]^{\frac{(p-1)}{2}} \int_{\Gamma(x^0) \cap \{ |y^\prime| \leq 1 \}} (m \cdot \nu) |y^\prime|^2 d\Gamma
\\& \quad\quad\quad\quad\quad\leq \frac{1}{2} \left[ E(t) \right]^{\frac{(p+1)}{2}} + C \left[ \int_{\Gamma(x^0) \cap \{ |y^\prime| \leq 1 \} } (m \cdot \nu) |y^\prime|^2 d\Gamma \right]^{\frac{p+1}{2}}
\end{split}
\end{equation}
but, by (\ref{eqVI55}),
\begin{equation}\label{eqVI81}
\begin{split}
&\left[ \int_{\Gamma(x^0) \cap \{ |y^\prime| \leq 1 \} } (m \cdot \nu) |y^\prime|^2 d\Gamma \right]^{\frac{p+1}{2}}
\\&\leq C\int_{\Gamma(x^0) \cap \{ |y^\prime| \leq 1 \} } (m \cdot \nu) |y^\prime|^{p+1} d\Gamma \leq C \int_{\Gamma(x^0)} (m \cdot \nu) g(y^\prime) y^\prime d\Gamma.
\end{split}
\end{equation}

Combining (\ref{eqVI78})-(\ref{eqVI81}) with $ \varepsilon >0 $, possibly smaller, it follows that
\begin{equation}\label{eqVI82}
    E_\varepsilon^\prime(t) \leq - \frac{\varepsilon}{2} \left[ E(t) \right]^{\frac{(p+1)}{2}} \leq - \frac{\varepsilon}{4} \left[ E_\varepsilon(t) \right]^{\frac{(p+1)}{2}},
\end{equation}
which implies (\ref{eqVI59}).

When $ p=1 $,  from (\ref{eqVI78}) the following holds:
\begin{equation}\label{eqVI83}
    E_\varepsilon^\prime(t) \leq - \varepsilon E(t) \leq - \frac{\varepsilon}{2} E_\varepsilon(t)
\end{equation}
which implies (\ref{eqVI57}).
\smallskip

{\bf Case 2. $\lambda \leq 1.$} We define
\label{case.2:proof:VI.lemma.6.1}
\begin{equation}\label{eqVI84}
    \phi (t) = \left[ E(t) \right]^{\frac{(p+1-2\lambda)}{2 \lambda}} \rho (t).
\end{equation}

We have
\begin{equation}\label{eqVI85}
\begin{split}
   \phi^\prime (t)  &= \frac{p+1-2\lambda}{2\lambda} \left[ E(t) \right]^{\frac{(p+1-4\lambda)}{2 \lambda}} \rho (t) E^\prime (t)
     + \left[ E(t) \right]^{\frac{(p+1-2\lambda)}{2 \lambda}} \rho^\prime(t)
    \\ & \leq -C E^\prime (t) - 2 \left[ E(t) \right]^{(p+1)/2 \lambda}
    \\ &\quad+ \left[ E(t) \right]^{\frac{(p+1-2\lambda)}{2 \lambda}} \int_{\Gamma(x^0)} (m \cdot \nu) \left[ R^2 | g(y^\prime) |^2 - (n-1) g (y^\prime) y + | y^\prime |^2 \right] d\Gamma.
\end{split}
\end{equation}

We observe that, thanks to (\ref{eqVI53}) and (\ref{eqVI55}):
\begin{equation}\label{eqVI86}
\begin{split}
&\left| \int_{\Gamma(x^0) \cap \{ |y^\prime| \leq 1 \} } (m \cdot \nu) g(y^\prime) y d\Gamma \right| 
\leq C \int_{\Gamma(x^0) \cap \{ |y^\prime| \leq 1 \} } (m \cdot \nu) |y^\prime|^\lambda |y| d\Gamma
\\&\leq C \int_{\Gamma(x^0) \cap \{ |y^\prime| \leq 1 \} } (m \cdot \nu) (g(y^\prime) y^\prime)^{\frac{\lambda}{(p+1)}} |y| d\Gamma
\\&\leq C \left[ \int_{\Gamma(x^0)} (m \cdot \nu) g(y^\prime) y^\prime d\Gamma \right]^{\frac{\lambda}{(p+1)}} \left[ \int_{\Gamma(x^0)} (m \cdot \nu) |y|^{\frac{p+1}{(p+1-\lambda)}} d\Gamma \right]^{\frac{p+1-\lambda}{(p+1)}}
\\&\leq C \left|E(t)\right|^{\frac{1}{2}} \left[ \int_{\Gamma(x^0)} (m \cdot \nu) g(y^\prime) y^\prime d\Gamma \right]^{\frac{\lambda}{(p+1)}}
\end{split}
\end{equation}
since 
$
(p+1)/(p+1-\lambda) \leq 2.
$

On the other hand, we have
\begin{equation}\label{eqVI87}
\begin{split}
&(n-1) \left| \int_{\Gamma(x^0) \cap \{ |y^\prime| \geq 1 \}} (m \cdot \nu) g(y^\prime) y d\Gamma \right|
\\&\leq C \left( \int_{\Gamma(x^0)} (m \cdot \nu) |y|^2 d\Gamma \right)^{\frac{1}{2}} \left( \int_{\Gamma(x^0) \cap \{ |y^\prime| \geq 1 \}} (m \cdot \nu) |g(y^\prime)|^2 d\Gamma \right)^{1/2}
\\&\leq C \left[E(t)\right]^{\frac{1}{2}} \left( \int_{\Gamma(x^0)} (m \cdot \nu) g(y^\prime) y^\prime d\Gamma \right)^{\frac{1}{2}}
\end{split}
\end{equation}
and, furthermore,
\begin{equation}\label{eqVI88}
\begin{split}
\int_{\Gamma(x^0)} (m \cdot \nu) |g(y^\prime)|^2 d\Gamma &\leq C
\int_{\Gamma(x^0) \cap \{ |y^\prime| \geq 1 \}} (m \cdot \nu) g(y^\prime) y^\prime d\Gamma
\\&\quad + C \int_{\Gamma(x^0) \cap \{ |y^\prime| \leq 1 \}} (m \cdot \nu) [g (y^\prime) y^\prime]^{\frac{2\lambda}{(1+\lambda)}} d\Gamma
\\&\leq -C E^\prime(t) + C \left[ \int_{\Gamma(x^0)} (m \cdot \nu) g(y^\prime) y^\prime d\Gamma \right]^{\frac{2\lambda}{1+\lambda}}.
\end{split}
\end{equation}

Finally
\begin{equation}\label{eqVI89}
\begin{split}
\int_{\Gamma(x^0)} (m \cdot \nu)& |y^\prime|^2 d\Gamma \leq C
\int_{\Gamma(x^0)} (m \cdot \nu) g(y^\prime) y^\prime d\Gamma 
\\&\quad\quad\quad\quad+ C \int_{\Gamma(x^0) \cap \{ |y^\prime| \geq 1 \}} (m \cdot \nu) [g (y^\prime) y^\prime]^{\frac{2}{(p+1)}} d\Gamma    \\
&\leq
    \begin{cases}
    -CE^\prime(t) & \hspace{0.1cm} \text{ if } \hspace{0.1cm} p \leq 1
    \\ -CE^\prime(t) + C \left[ \int_{\Gamma(x^0)} (m \cdot \nu) g(y^\prime) y^\prime d\Gamma \right]^{\frac{2}{(p+1)}} & \hspace{0.1cm} \text{ if } \hspace{0.1cm} p > 1.
    \end{cases}
    \end{split}
\end{equation}

Combining (\ref{eqVI85})-(\ref{eqVI89}) we get
\begin{equation}\label{eqVI90}
\begin{split}
\phi^\prime(t) & \leq -CE^\prime(t) - 2 \left[ E(t) \right]^{\frac{(p+1)}{2\lambda}}
 + \left[ E(t) \right]^{\frac{(p+1-\lambda)}{2\lambda}} 
\Bigg[ \left[ \int_{\Gamma(x^0)} (m \cdot \nu) g (y^\prime) y^\prime d\Gamma  \right]^{\frac{\lambda}{(p+1)}} 
\\ &\quad+ \left[ \int_{\Gamma(x^0)} (m \cdot \nu) g (y^\prime) y^\prime d\Gamma  \right]^{\frac{1}{2}} \Bigg]
 + \left[ E(t) \right]^{\frac{(p+1-2\lambda)}{2\lambda}}
 \left[ \int_{\Gamma(x^0)} (m \cdot \nu) g(y^\prime) y^\prime d\Gamma \right]^{\frac{2 \lambda}{(1+ \lambda)}}
\\ & \quad\left(+ \left[ E(t) \right]^{\frac{(p+1-2\lambda)}{2\lambda}}
 \left[ \int_{\Gamma(x^0)} (m \cdot \nu) g(y^\prime) y^\prime d\Gamma \right]^{\frac{2}{(p+1)}} \hspace{0.3cm} \text{ if } p>1 \right).
\end{split}
\end{equation}

Applying Young's inequality in (\ref{eqVI90}) we have
\begin{equation}\label{eqVI91}
\phi^\prime(t) \leq -C_3 E^\prime(t) - \left[ E(t) \right]^{\frac{(p+1)}{2\lambda}} 
\end{equation}

with $C_3 >0$ continuously depending on $E(0)$.

Given $ \varepsilon >0 $ we define the functional
\begin{equation}\label{eqVI92}
E_\varepsilon(t) = (1 + \varepsilon C_3) E(t) + \varepsilon \phi (t).
\end{equation}

As $ E(t) \leq E(0) $, taking $ \varepsilon \leq \varepsilon_0 (E(0)) $ we have
\begin{equation}\label{eqVI93}
\frac{1}{2} \left[ E_\varepsilon(t) \right]^\frac{(p+1)}{2 \lambda} \leq \left[ E(t) \right]^\frac{(p+1)}{2 \lambda} \leq 2 \left[ E_\varepsilon(t) \right]^\frac{(p+1)}{2 \lambda}.
\end{equation}

On the other hand, from (\ref{eqVI29}) and (\ref{eqVI91}), it follows
\begin{equation}\label{eqVI94}
E_\varepsilon^\prime(t) \leq E^\prime(t) - \varepsilon \left[ E(t) \right]^\frac{(p+1)}{2 \lambda} \leq - \frac{\varepsilon}{2} \left[ E_\varepsilon(t) \right]^\frac{(p+1)}{2 \lambda}
\end{equation}
from where we may conclude (\ref{eqVI58}).
\end{proof}

\section{Second Proof of the Exponential Decay}
\label{sec:VI.5}
In this section we are going to present a second proof of the exponential decay of the energy of the solutions of (\ref{eqVI52}) with linear boundary conditions. It has the advantage of being applicable in situations where we do not know how to build energy functionals allowing, as in the previous section, to obtain decay rates. Its biggest drawback is that it mainly applies with linear boundary conditions (or  growing linearly both at the origin and at infinity).

Let us consider the system
\begin{equation}\label{eqVI95}
    \begin{cases}
    y^{\prime\prime} - \Delta y = 0 & \text{in} \hspace{0.2cm} \Omega \times  (0, \infty)
    \\
    \frac{\partial y}{\partial v} + by = - [ (x-x^0) \cdot \nu(x) ] y^\prime \hspace{0.2cm} & \text{in} \hspace{0.2cm} \Gamma (x^0) \times (0, \infty)
    \\ y = 0 \hspace{0.2cm} & \text{in} \hspace{0.2cm} \Gamma_\star (x^0) \times (0, \infty)
    \\
    y(0) = y^0 \in V, y^\prime (0) = y^1 \in L^2(\Omega) & \text{in} \hspace{0.2cm} \Omega.
    \end{cases}
\end{equation}
We assume that
\begin{equation}\label{eqVI96}
b > -\beta
\end{equation}
where
\begin{equation}\label{eqVI97}
\beta = \beta (x^0) = \inf_{v \in V} \frac{\int_\Omega |\nabla v|^2 dx}{\int_{\Gamma(x^0)} |v|^2 d\Gamma}
\end{equation}

When $ \operatorname{int} (\Gamma_\star(x^0)) \neq \emptyset $ we have $ \beta > 0 $ and therefore the case $ b=0 $ is included in (\ref{eqVI96}). It would correspond to (\ref{eqVI52}) with $ g (s) = s$.

When $ \operatorname{int} (\Gamma_\star(x^0)) = \emptyset, \beta = 0 $ and therefore (\ref{eqVI96}) corresponds to $ b>0 $.

The energy associated with this system is
\begin{equation}\label{eqVI98}
E_b(t) = \frac{1}{2} \int_\Omega \left[ |\nabla y (x,t)|^2 + |y^\prime(x,t)|^2 \right] dx + \frac{b}{2} \int_{\Gamma(x^0)} |y(x,t)|^2 d\Gamma
\end{equation}
and thanks to (\ref{eqVI96}), we can ensure that $ [ E_b ]^{1/2} $ defines a norm in  $ V \times L^2(\Omega) $, equivalent to the one induced by $ H^1(\Omega) \times L^2(\Omega) $.

On the other hand,
\begin{equation}\label{eqVI99}
\frac{dE_b(t)}{dt} = E_b^\prime(t) = - \int_{\Gamma(x^0)} \left[ (x-x^0) \cdot \nu(x) \right] |y^\prime(x,t)|^2  d\Gamma.
\end{equation}

Therefore, the equilibrium points, i.e. the non dissipative solutions, satisfy
\begin{equation}\label{eqVI100}
    \begin{cases}
    y^{\prime\prime} - \Delta y = 0 \hspace{0.2cm} & \text{in} \hspace{0.2cm} \Omega \times (0, \infty)
    \\
    \frac{\partial y}{\partial v} + by = y^\prime = 0 \hspace{0.2cm} & \text{in} \hspace{0.2cm} \Gamma (x^0) \times (0, \infty)
    \\ y = 0 \hspace{0.2cm} & \text{in} \hspace{0.2cm} \Gamma_\star (x^0) \times (0, \infty).
    \end{cases}
\end{equation}

But, applying Holmgren's Theorem, it follows that $ y^\prime \equiv 0 $, i.e. $ y= y(x) $. But then:
\begin{equation}\label{eqVI101}
    \begin{cases}
    - \Delta y = 0 \hspace{0.2cm} & \text{in} \hspace{0.2cm} \Omega
    \\
\partial y/\partial v + by = 0 \hspace{0.2cm} & \text{in} \hspace{0.2cm} \Gamma (x^0)\\
     y=0 \hspace{0.2cm} &\text{in} \hspace{0.2cm} \Gamma_\star (x^0).
\end{cases}
\end{equation}

Multiplying in (\ref{eqVI101}) by $y$ and integrating by parts we obtain
\begin{equation*}
    \int_\Omega | \nabla y (x) |^2 dx + b \int_{\Gamma(x^0)} |y(x)|^2 d\Gamma = 0
\end{equation*}
which, combined with (\ref{eqVI96}), implies $ y \equiv 0 $.

These remarks suggest the exponential decay of the energy $ E_b(t) $ when $ t \rightarrow \infty $, or, what is the same, the exponential convergence of all trajectories $ \{ y(t), y^\prime(t) \} $ to the only point of equilibrium $ \{ 0, 0 \} $ in $ H^1(\Omega) \times L^2(\Omega) $ when $ t \rightarrow + \infty $.

We have the following theorem (see \cite{zuazua1989remarks}).

\begin{theorem}
\label{theorem:VI.6.4}
Let $ \Omega $ be a bounded domain of $ \mathbb{R}^n $ with a boundary of class $ C^2 $ and $n\leq3$. Let $ x^0 \in \mathbb{R}^n $ and $ b > -\beta (x^0) $.

Then, there exist constants $ C>1 $ and $ \gamma >0 $ such that
\begin{equation}\label{eqVI102}
E_b(t) \leq C E_b(0) e^{-\gamma t}, \, \forall t \geq 0
\end{equation}
for every solution of (\ref{eqVI95}), where $E_b(t) $ is the energy defined in (\ref{eqVI98}).
\end{theorem}

\begin{remark}
\label{remark:VI.6.8}
Again, if $ \overline{\Gamma(x^0)} \cap \overline{\Gamma_\star(x^0)} = \emptyset $, the theorem is valid without any restriction on the dimension $n\ge 1$.
$ \hfill \square $
\end{remark}


\begin{remark}
\label{remark:VI.6.9}
The results of existence, uniqueness, regularity and stability of solutions in Sections 6.2 and 6.3 extend to system (\ref{eqVI95}).
$ \hfill \square $
\end{remark}

\begin{proof}
\label{proof:VI.remark.6.9}
Note that it is enough to prove the existence of a time $ T>0 $ and a constant $ C_0 > 0 $ such that
\begin{equation}\label{eqVI103}
    E_b(T) \leq C_0 \int^T_0 \int_{\Gamma(x^0)} (m \cdot \nu)|y^\prime|^2 d\Gamma dt
\end{equation}
for every solution of (\ref{eqVI95}). In (\ref{eqVI103}), $ m = m(x) =x - x^0 $.

Indeed, combining (\ref{eqVI103}) with (\ref{eqVI99}), which implies
\begin{equation}\label{eqVI104}
    E_b(T) = E_b(0) - \int^T_0 \int_{\Gamma(x^0)} (m \cdot \nu)|y^\prime|^2 d\Gamma dt,
\end{equation}
we get
\begin{equation}\label{eqVI105}
    E_b(T) \leq \left( \frac{C_0}{1+C_0} \right)  E_b(0).
\end{equation}

Combining (\ref{eqVI105}) with the decreasing character of the energy and the semigroup property, it follows that
\begin{equation}\label{eqVI106}
    E_b(t) \leq \left( \frac{C_0}{1+C_0} \right)^{\frac{t}{T} - 1}  E_b(0)
\end{equation}
from where we get (\ref{eqVI102}) with
\begin{equation}\label{eqVI107}
C = \frac{1+C_0}{C_0}; \hspace{0.3cm} \gamma = \frac{1}{T} \operatorname{log} \left[ \frac{1+C_0}{C_0} \right].
\end{equation}

By density, it is enough to consider solutions with regular initial data $ \{ y^0, y^1 \} \in V \times V $ such that $ -\Delta y^0 \in L^2(\Omega) $, satisfying the compatibility condition
\begin{equation*}
    \frac{\partial y^0 }{\partial v} + b y^0 = -(m \cdot \nu) y^1 \hspace{0.3cm} \text{ in } \hspace{0.3cm} \Gamma(x^0).
\end{equation*}

In this case the solution $ y = y(x,t) $ of (\ref{eqVI95}) satisfies
\begin{equation*}
 y \in C^1([0, \infty); V), \hspace{0.3cm} \Delta y\in C([0,\infty); L^2(\Omega)).
\end{equation*}

We multiply equation (\ref{eqVI95}) by $ (x- x^0) \cdot \nabla y $ and integrate by parts in $ \Omega \times (0,T) $ (where $ T>0 $ will be fixed later). Applying similar arguments as in Lemma 6.1, which is possible since $ n \leq 3 $, we obtain (hereinafter, to simplify the notation, we write $ E(t) $ instead of $ E_b(t) $):
\begin{equation}\label{eqVI108}
\begin{split}
    & \int_\Omega y^\prime m \cdot \nabla y dx \bigg|^T_0 + \int^T_0 E(t) dt + \frac{(n-1)}{2} \int^T_0 \int_\Omega \left[ |y^\prime|^2 - | \nabla y |^2  \right] dx dt
    \\ & \leq \frac{b}{2} \int^T_0 \int_{\Gamma(x^0)} |y|^2 d\Gamma dt + \frac{1}{2} \int^T_0 \int_{\Gamma(x^0)} ( m \cdot \nu ) |y^\prime|^2 d\Gamma dt
    \\ & \quad+ \int^T_0 \int_\Gamma \frac{\partial y}{\partial \nu}  m \cdot \nabla y d\Gamma dt - \frac{1}{2} \int^T_0 \int_\Gamma (m \cdot \nu) |\nabla y|^2 d\Gamma dt
    \\ & \leq \frac{b}{2} \int^T_0 \int_{\Gamma(x^0)} |y|^2 d\Gamma dt + \frac{1}{2} \int^T_0 \int_{\Gamma(x^0)} ( m \cdot \nu ) |y^\prime|^2 d\Gamma dt
    \\ & \quad+ \frac{1}{2} \int^T_0 \int_{\Gamma_\star(x^0)} ( m \cdot \nu ) \left| \frac{\partial y}{\partial \nu} \right|^2 d\Gamma dt + \frac{1}{2} \int^T_0 \int_{\Gamma(x^0)} ( m \cdot \nu ) \left| \frac{\partial y}{\partial \nu} \right|^2 d\Gamma dt
    \\ &  \quad- b \int^T_0 \int_{\Gamma(x^0)} y m \cdot \nabla_\sigma y d\Gamma dt - \int^T_0 \int_{\Gamma(x^0)} (m \cdot \nu) y^\prime m \cdot \nabla_\sigma y d\Gamma dt
    \\ &  \quad- \frac{1}{2} \int^T_0 \int_{\Gamma(x^0)} (m \cdot \nu) \nabla_\sigma |y|^2 d\Gamma dt
    \\ & \leq C \int^T_0 \int_{\Gamma(x^0)} |y|^2 d\Gamma dt + C \int^T_0 \int_{\Gamma(x^0)} (m \cdot \nu) |y^\prime|^2 d\Gamma dt
\end{split}
\end{equation}
for $ C >0 $ depending only on $ m $ and $ b $, since $ m \cdot \nu \leq 0 $ in $ \Gamma_\star (x^0) $.

We have used the identity
\begin{equation*}
    \int_{\Gamma(x^0)} y m \cdot \nabla_\sigma y d\Gamma =   \int_{\Gamma(x^0)} \frac{m}{2} \cdot \nabla_\sigma ( |y|^2 ) d\Gamma = 0
\end{equation*}
and $ \| \partial y / \partial \nu \|^2 \leq 2b^2 |y|^2 + 2 ( m \cdot \nu )^2 |y^\prime|^2 $ in $ \Gamma(x^0) $. We have also used the inequality
\begin{equation*}
    ( m \cdot \nu ) | y^\prime m \cdot \nabla_\sigma y | \leq \frac{1}{2} ( m \cdot \nu ) | \nabla_\sigma y |^2 + \frac{R^2 ( m \cdot \nu )}{2} |y^\prime|^2 \hspace{0.2cm} \text{ in } \hspace{0.2cm} \Gamma(x^0)
\end{equation*}
with $ R = \| m \|_{L^\infty(\Omega)}  $.

On the other hand, multiplying equation (\ref{eqVI95}) by $ y $ and integrating by parts, we obtain
\begin{equation}\label{eqVI109}
\begin{split}
\int^T_0 \int_\Omega \left[ | y^\prime |^2 - | \nabla y |^2 \right] dx dt & = \int_\Omega y^\prime y dx \bigg|^T_0 + b \int^T_0 \int_{\Gamma(x^0)} |y|^2 d\Gamma dt
\\ &  \quad+ \int^T_0 \int_{\Gamma(x^0)} ( m \cdot \nu ) y y^\prime d\Gamma dt.
\end{split}
\end{equation}

Combining (\ref{eqVI108}) and (\ref{eqVI109}) we deduce
\begin{equation}\label{eqVI110}
\begin{split}
 \int_\Omega y^\prime \bigg[ m \cdot \nabla y &+ \frac{(n-1) y}{2} \bigg] dx \bigg|^T_0 + \int^T_0 E(t) dt
\\ & \leq C \int^T_0 \int_{\Gamma(x^0)} |y|^2 d\Gamma dt + C \int^T_0 \int_{\Gamma(x^0)} (m \cdot \nu) |y^\prime|^2 d\Gamma dt.
\end{split}
\end{equation}

Now observe that there exists $ C_1 > 0 $ (which depends only on $ \Omega $),  such that
\begin{equation*}
    \left| \int_\Omega y^\prime \left[ m \cdot \nabla y + \frac{(n-1)}{2} y \right] dx \right| \leq C_1 E(t)
\end{equation*}
and therefore
\begin{equation}\label{eqVI111}
\begin{split}
    \bigg| \int_\Omega y^\prime \bigg[ m \cdot \nabla y &+ \frac{(n-1)}{2} y \bigg] dx \bigg|^T_0 \bigg|  \leq C_1 ( E(0) + E(T) )
    \\ & \leq 2 C_1 E(T) + C_1 \int^T_0 \int_{\Gamma(x^0)} ( m \cdot \nu ) |y^\prime|^2 d\Gamma dt.
\end{split}
\end{equation}

From (\ref{eqVI110}), (\ref{eqVI111}) and the fact that
$ T E(T) \leq \int^T_0 E(t) dt $ we obtain
\begin{equation}\label{eqVI112}
\begin{split}
    (T - 2C_1) E(T) \leq C \int^T_0 \int_{\Gamma(x^0)} |y|^2 d\Gamma dt + C \int^T_0 \int_{\Gamma(x^0)} ( m \cdot \nu ) |y^\prime|^2 d\Gamma dt.
\end{split}
\end{equation}

Therefore, it is enough to prove that, if $ T > 2C_1 $, there is a constant $ C>0 $ such that
\begin{equation}\label{eqVI113}
\begin{split}
    \int^T_0 \int_{\Gamma(x^0)} |y|^2 d\Gamma dt \leq C \int^T_0 \int_{\Gamma(x^0)} ( m \cdot \nu ) |y^\prime|^2 d\Gamma dt.
\end{split}
\end{equation}

We argue by contradiction. If (\ref{eqVI113}) does not hold, there is a sequence $ \{ y_n \} $ of solutions of (\ref{eqVI95}) such that
\begin{equation}\label{eqVI114}
\begin{cases}
    \int^T_0 \int_{\Gamma(x^0)} |y_n|^2 d\Gamma dt = 1 \hspace{0.6cm} \forall n \in \mathbb{N}
    \\ 
    \\ \int^T_0 \int_{\Gamma(x^0)} ( m \cdot \nu ) |y_n^\prime|^2 d\Gamma dt \rightarrow 0.
\end{cases}
\end{equation}

Combining (\ref{eqVI114}) with (\ref{eqVI112}) and taking into account that $ T > 2 C_1 $ we obtain that
\begin{equation*}
 \{ y_n \} \hspace{0.2cm} \text{ is bounded in } \hspace{0.2cm} C( [0,T];V ) \cap C^1 ( [0,T];L^2(\Omega) )
\end{equation*}
and therefore
\begin{equation*}
 \left\{ y_{n} |_{\Gamma \times (0,T) }\right\}\hspace{0.2cm} \text{ is relatively compact in } \hspace{0.2cm} L^2(\Gamma \times (0,T)).
\end{equation*}

Extracting a subsequence and passing to the limit, we obtain a solution $ y = y(x,t) $ of (\ref{eqVI95}), thanks to (\ref{eqVI114}), satisfying:
\begin{equation}\label{eqVI115}
    y^\prime = 0 \hspace{0.3cm} \text{ in } \hspace{0.3cm} \Gamma(x^0) \times (0,T)
\end{equation}
\begin{equation}\label{eqVI116}
    \int^T_0 \int_{\Gamma(x^0)} |y|^2 d\Gamma dt = 1.
\end{equation}

But, as we observed at the beginning of this section, by Holmgren's Theorem, if $ T>0 $ is large enough, the only solution of (\ref{eqVI95}) satisfying (\ref{eqVI115}) is $ y \equiv 0 $, which contradicts (\ref{eqVI116}). This concludes the proof of Theorem \ref{theorem:VI.6.4}.
\end{proof}

\label{sec:VI.5-f}
\begin{remark}
The same result is obtained if in (\ref{eqVI95}) we replace the constant $ b $ by a function $ b = b(x) \in C^1 (\Gamma(x^0)) $ such that
$
    b \geq b_0 > 0
$
in an open and non-empty subset of $ \Gamma(x^0) $.

Theorem \ref{theorem:VI.6.4} can be extended to systems of type (\ref{eqVI95}) with boundary conditions
\begin{equation}\label{eqVI117}
    \frac{\partial y}{\partial \nu} + b(x)y = -\left[ (x-x^0) \cdot \nu(x) \right] g(y^\prime) \hspace{0.2cm} \text{ in } \hspace{0.2cm} \Gamma(x^0) \times (0, \infty)
\end{equation}
with $ g \in C(\mathbb{R}) $ increasing and satisfying
\begin{equation*}
    c|s| \leq |g(s)| \leq C|s| \hspace{0.6cm} \forall s \in \mathbb{R}.
\end{equation*}
We do not know if these arguments can be extended to the setting of Theorem \ref{theorem:VI.6.3}.

In \cite{zuazua1990stabilization} Theorem \ref{theorem:VI.6.3} is extended to boundary conditions of the type (\ref{eqVI117}), but with
$
    b(x) = k [ (x - x^0) \cdot \nu(x) ]
$
where $ k>0 $ is a sufficiently small constant.
\end{remark}

\section{Comments}
\label{sec:VI.6}
\begin{enumerate}
\item {\bf Variable coefficients.} The results of Section \ref{sec:VI.3}, based on LaSalle's Invariance Principle and energy methods, can be easily generalized to wave equations with variable coefficients.

However, the adaptation of the results of Sections \ref{sec:VI.4} and \ref{sec:VI.5} to the variable coefficients case is much more delicate since in their proof we have used multiplier techniques.
In fact, the extension of Theorem \ref{theorem:VI.6.3} to equations with variable coefficients (but without imposing the technical restrictions mentioned in the previous chapters that are due to the multiplier technique itself) is an open problem.

Theorem 6.4 can be easily extended to wave equations with variable $BV$ coefficients in one space dimension $ n=1 $.

On the other hand, the microlocal analysis techniques of \cite{bardos1992sharp} allow obtaining stabilization results like Theorem \ref{theorem:VI.6.4} for equations with variable coefficients and for partitions of the boundary $ \{ \Gamma_0, \Gamma_1 \} $ satisfying suitable geometric control conditions, more general than $ \{ \Gamma(x^0), \Gamma_\star(x^0) \} $. The case where $ \overline{\Gamma_0} \cap \overline{\Gamma_1} \neq \emptyset $, i.e. when the boundary conditions change type, needs further developments from a microlocal perspective.

\item  {\bf Semilinear wave equations.}  Theorem \ref{theorem:VI.6.3} can be extended to semilinear wave equations
$
    y^{\prime\prime} - \Delta y + f(y) = 0
$
with dissipative boundary conditions (cf. \cite{zuazua1989remarks}).

\item {\bf Other models.}  The techniques that we have developed in this chapter also allow obtaining stabilization results for the system of elasticity (cf. J. Lagnese \cite{lagnese1983boundary}) and for the Schrödinger equation (cf. E. Machtyngier \cite{machtyngier1990control}).

\item {\bf Stabilisation in weaker energy norms.} The problem of stabilization of the wave equation arises also in other functional settings and with other boundary conditions.

Let $ G = (-\Delta)^{-1} : H^{-1}(\Omega) \to H^1_0(\Omega) $ and  consider the system
\begin{equation}\label{eqVI118}
    \begin{cases}
    y^{\prime\prime} - \Delta y = 0 \hspace{0.2cm} & \text{in} \hspace{0.2cm} \Omega \times (0, \infty)
    \\
    y = \partial (G y^\prime)/\partial \nu \hspace{0.2cm} & \text{in} \hspace{0.2cm} \Gamma_0 \times (0, \infty)
    \\ y = 0 \hspace{0.2cm} & \text{in} \hspace{0.2cm} \Gamma_1 \times (0, \infty)
    \\ y(0) = y^0 \in L^2(\Omega), y^\prime (0) = y^1 \in H^{-1}(\Omega)
    & \text{in} \hspace{0.2cm} \Omega
    \end{cases}
\end{equation}
where $ \{ \Gamma_0, \Gamma_1 \} $ is a partition of the boundary $ \Gamma $.

In I. Lasiecka and R. Triggiani \cite{lasiecka1987uniform} it is shown that (\ref{eqVI118}) admits a single solution in the class
\begin{equation}\label{eqVI119}
    y \in C( [0,\infty); L^2(\Omega) ) \cap C^1 ( [0, \infty) ; H^{-1}(\Omega) ).
\end{equation}

The energy associated with the system is
\begin{equation}\label{eqVI120}
    E(t) = \frac{1}{2} \int_\Omega | y(x,t) |^2 dx + \frac{1}{2} \| y^\prime(t) \|^2_{H^{-1}(\Omega)}
\end{equation}
with $ \| u \|^2_{H^{-1}(\Omega)} = \langle G u, u \rangle $, where $ \langle \cdot , \cdot \rangle $ denotes the duality product between $ H^1_0(\Omega) $ and $H^{-1}(\Omega)$.
It fulfils
\begin{equation}\label{eqVI121}
    \frac{dE(t)}{dt} = - \int_{\Gamma_0} | y^\prime |^2 d\Gamma \hspace{0.6cm} \forall t \geq 0
\end{equation}
which allows to prove, applying LaSalle's Invariance Principle, that
\begin{equation}\label{eqVI122}
    E(t) \to 0 \hspace{0.3cm} \text{ when } \hspace{0.3cm} t \rightarrow \infty
\end{equation}
for every solution of (\ref{eqVI118}), as long as
\begin{equation}\label{eqVI123}
    \operatorname{int} ( \Gamma_0 ) \neq \emptyset.
\end{equation}

In \cite{lasiecka1987uniform}, under suitable geometric conditions on the domain $ \Omega $, it is proved that if $ \Gamma_0 = \Gamma (x^0) $ and $ \Gamma_1 = \Gamma_\star(x^0) $, the uniform exponential decay of the energy holds, i.e.
\begin{equation}\label{eqVI124}
    E(t) \leq C E(0) e^{-\gamma t} \hspace{0.6cm} \forall t \geq 0.
\end{equation}

The results of \cite{lasiecka1987uniform} apply essentially to the case where $ \Omega = \Omega_0 \backslash \overline{\Omega_1} $ with $ \overline{\Omega_1}\subset \Omega_0 $ and such as $ \Gamma(x^0) = \partial \Omega_0, \Gamma_\star (x^0) = \partial \Omega_1 $.

In \cite{bardos1988stabilisation} using Microlocal Analysis techniques, it was proved that if $ \Omega $ is regular (without any geometric condition on $ \Omega $) and $ \Gamma_0 $ geometrically controls the domain $ \Omega $, then  (\ref{eqVI124}) holds as long as $ \overline{\Gamma_0} \cap \overline{\Gamma_1} = \emptyset $.
The case where $ \overline{\Gamma_0} \cap \overline{\Gamma_1} \neq \emptyset $ is more delicate due to the possible presence of singularities at the interface points where boundary conditions change type.

In \cite{bardos1988stabilisation} stabilization results were also proved for other types of dissipative boundary conditions. For example, for absorbent or transparent boundary conditions (cf. B. Engquist and L. Halpern \cite{engquist1988far}).

\item {\bf Slow decay.}  We refer to \cite{lebeau1997stabilisation} for the analysis of the slow logarithmic decay of the smooth solutions of the linearly damped wave equation, when the damping is localized in an arbitrarily small subset of the boundary. Note that, in order to avoid difficulties with the possible singularities on the interface points where the boundary conditions change, in \cite{lebeau1997stabilisation} authors consider boundary conditions of the form
\begin{equation*}
    \frac{\partial y}{\partial v} + a(x) y^\prime = 0 \hspace{0.3cm} \text{on} \hspace{0.3cm} \Gamma \times (0, \infty)
\end{equation*}
with $ a \in C^\infty(\Gamma) $, and, $a \geq 0$, vanishing on a subset of the boundary.

\item {\bf Backstepping.} The methods developed in this paper do not directly apply to unstable wave-like equations. They are rather designed to force the decay of solutions of conservative models by means of added damping terms. Backstepping techniques allow to design non-local mechanisms that can effectively stabilize unstable systems. The literature on this topic is very broad. See for instance \cite{krstic2008boundary}. Note however that the design of backstepping feedback operators requires  the introduction of suitable kernels, whose existence is often guaranteed only for particular geometries on the domains in which waves propagate.
\end{enumerate}


%
%

\chapter{Further Reading}

\abstract{The topics addressed in these Notes have significantly evolved since they were written for the first time back in 1989. Here we enclose some other important references that complement the synthetic presentation of these Notes. }

\section{Texts on ODE and/or PDE Control}
\label{sec:odepde}
\begin{itemize}
    \item {Bastin, G., \& Coron, J. M. (2016). \emph{Stability and Boundary Stabilization of 1-d Hyperbolic Systems} (Vol. 88). Basel: Birkhäuser.}
    \item {Coron, J. M. (2007). \emph{Control and Nonlinearity} (No. 136). American Mathematical Society.}
    \item {Gugat, M. (2015). Optimal Boundary Control and Boundary Stabilization of Hyperbolic Systems. Birkhäuser.}
    \item {Komornik, V. (1994). \emph{Exact Controllability and Stabilization: the Multiplier Method} (Vol. 39, p. 351). Chichester: Wiley.}
    \item {Trélat, E. (2005). \emph{Contrôle optimal: théorie \& applications} (Vol. 36). Paris: Vuibert.}
    \item {Tucsnak, M., \& Weiss, G. (2009). \emph{Observation and Control for Operator Semigroups}. Springer Science \& Business Media.}
\end{itemize}

\section{Survey Articles}
\label{sec:survey}
\begin{itemize}
    \item {Zuazua, E. (2002). Controllability of partial differential equations and its semi-discrete approximations. \emph{Discrete and continuous dynamical systems}, 8(2), 469-517.}
    \item {Zuazua, E. (2007). Controllability and observability of partial differential equations: some results and open problems. In \emph{Handbook of differential equations: evolutionary equations} (Vol. 3, pp. 527-621). North-Holland.}
\end{itemize}

\section{Control of Fractional PDE}
\label{sec:fracpde}
\begin{itemize}
    \item {Biccari, U., Warma, M., \& Zuazua, E. (2022). Control and numerical approximation of fractional diffusion equations. \emph{Handbook of Numerical Analysis XXIII. Elsevier. Numerical Control: Part A}, 23, 1-58.}
\end{itemize}

\section{Numerical Analysis of Wave Control}
\label{sec:numanal}
\begin{itemize}
    \item {Ervedoza, S., \& Zuazua, E. (2013). \emph{Numerical Approximation of Exact Controls for Waves}. SpringerBriefs in Mathematics. Springer, New York.}
    \item {Marica, A., \& Zuazua, E. (2014). \emph{Symmetric Discontinuous Galerkin Methods for 1-D Waves: Fourier Analysis, Propagation, Observability and Applications}. Springer, New York.}
    \item {Zuazua, E. (2005). Propagation, observation, and control of waves approximated by finite difference methods. \emph{SIAM Review}, 47(2), 197-243.}
\end{itemize}

\section{Control of PDEs on Networks}
\label{sec:pdenets}
\begin{itemize}
    \item {Dáger, R., \& Zuazua, E. \emph{Wave Propagation, Observation and Control in $1-d$ Flexible Multi-structures}. Springer Science \& Business Media, 2006.}
    \item {Lagnese, J. E., Leugering G., and Schmidt, E.J.P.G. \emph{Modeling, Analysis and Control of Dynamic Elastic Multi-link Structures}. Springer Science \& Business Media, 2012.}
\end{itemize}

\backmatter

\printindex

\bibliographystyle{plain}
\nocite{*}
\bibliography{biblio}

\end{document}